\documentclass[12pt]{article}
\usepackage{sectsty}
\usepackage{theorem}
\usepackage{amssymb}
\usepackage{graphicx}
\usepackage[mathscr]{eucal}
\usepackage{amsbsy}
\usepackage{amsmath}
\usepackage{multirow}
\newtheorem{prop}{}[section]

{\theorembodyfont{\upshape} \newtheorem{rema}[prop]{}}
\newcommand{\boma}[1]{{\mbox{\boldmath $#1$} }}
\newcommand{\Di}[1]{D^{{\,#1}}}
\newcommand{\Lap}[1]{\sqrt{-\Delta}^{\,#1}}

\newcommand{\WW}[2]{ \dot{\LL}^{#1 , #2}}
\newcommand{\ZZ}[3]{ \LL^{#1 , #2 , #3}}
\newcommand{\Hh}[1]{ \dot{H}^{#1}}
\newcommand{\HH}[1]{ H^{#1} }
\begin{document}
\def\L{\mathfrak{P}}
\def\param{\mathfrak{p}}
\def\alf{\Upsilon}
\def\Ge{\mathfrak{G}}
\def\Se{\mathfrak{S}}
\def\AZ{Z}
\def\N{N}
\def\GNA{\mathscr{A}}
\def\GA{A}
\def\GNB{\mathscr{B}}
\def\GB{B}
\def\GN{$\mathscr{G}\mathscr{N}~$}
\def\GNE{\hat{\mathscr{G}}}
\def\SO{$\mathscr{S}~$}
\def\HL{\mathscr{H}}
\def\EF{\mathfrak{F}}
\def\Is{\mathfrak{I}}
\def\del{\delta}
\def\Sp{S_{+}}
\def\Spp{S_{++}}
\def\Sm{S_{-}}
\def\Smm{S_{--}}
\def\Gp{G_{+}}
\def\Gpp{G_{++}}
\def\Gm{G_{-}}
\def\Gmm{G_{--}}
\def\PP{{\mathcal P}}
\def\PPp{Pol}
\def\t{{r'}}
\def\ep{\varepsilon}
\def\CC{{\mathcal C}}
\def\DD{\mathcal{D}}
\def\DDd{Dir}
\def\Sob{\mathfrak{S}}
\def\Gag{\mathfrak{G N}}
\def\Gagg{\overset{\circ}{\mathfrak{G N}}}
\def\Gg{\overset{\circ}{G}}
\def\si{\sigma}
\def\al{\alpha}
\def\rs{\hat{r}}
\def\So{\mathcal{S}}
\def\Ga{\mathcal{G}}
\def\Sa{\mathcal{S}}
\def\pb{\hskip -0.1cm \phantom{$\Bigg|$}}
\def\pbb{\hskip -0.1cm \phantom{$\Big|$}}
\def\RR{\mathscr R}
\def\LL{\mathscr L}
\def\BB{\mathscr B}
\def\intl{\int_{-i \infty}^{+ i \infty} \hspace{-0.1cm}}
\def\intc{\int_{c -i \infty}^{c+ i \infty} \hspace{-0.1cm}}
\def\intcn{\int_{{c \over N} -i \infty}^{{c \over N} + i \infty} \hspace{-0.1cm}}
\def\bs{b^{*}}
\def\as{a^{*}}
\def\Bs{B^{*}}
\def\As{A^{*}}
\def\ns{{n^{*}}}
\def\ms{{m^{*}}}
\def\T{{\mathfrak S}}
\def\te{\vartheta}
\def\leqs{\leqslant}
\def\geqs{\geqslant}
\def\Mm{{\mathfrak M}}
\def\Dd{{\mathfrak D}}
\def\s{\scriptstyle}
\def\ss{\scriptscriptstyle}
\def\dd{\displaystyle}
\def\Hn{ H^{n}(\reali^d, \complessi) }
\def\Ld{\LL^{2}(\reali^d, \complessi)}
\def\Lpi{L^{p}(\reali^d, \complessi)}
\def\Lq{L^{q}(\reali^d, \complessi)}
\def\Lr{L^{r}(\reali^d, \complessi)}
\def\k{\mbox{{\tt k}}}
\def\x{\mbox{{\tt x}}}
\def\lgraffa{ \mbox{\Large $\{$ } \hskip -0.2cm}
\def\rgraffa{ \mbox{\Large $\}$ } \hskip -0.2cm}
\def\sc{ {\scriptstyle{\bullet} }}
\def\la{\lambda}
\def\complessi{{\mathbb{C}}}
\def\reali{{\mathbb{R}}}
\def\interi{{\mathbb{Z}}}
\def\naturali{{\mathbb{N}}}
\def\FF{{\mathcal F}}
\def\SS{{\mathcal S}}
\def\circa{\thickapprox}
\def\vain{\rightarrow}
\def\parn{\par \noindent}
\def\salto{\vskip 0.2truecm \noindent}
\def\saltino{\vskip 0.1truecm \noindent}
\def\spazio{\vskip 0.5truecm \noindent}
\def\vs1{\vskip 1cm \noindent}
\def\fine{\hfill $\square$ \vskip 0.2cm \noindent}
\newcommand{\rref}[1]{(\ref{#1})}
\def\beq{\begin{equation}}
\def\feq{\end{equation}}
\def\beqq{\begin{eqnarray}}
\def\feqq{\end{eqnarray}}
\def\barray{\begin{array}}
\def\farray{\end{array}}
%%%%%%%%% THIS NUMBERS EQUATIONS BY SECTIONS %%%%%%%%%%%%%
\makeatletter
\@addtoreset{equation}{section}
\renewcommand{\theequation}{\thesection.\arabic{equation}}
%\thesection instead of \arabic{section} for correct equation numbering
% in appendices
\makeatother
%%%%%%%%%%%%%%%%%%%%%%%%%%%INTESTAZIONE%%%%%%%%%%%%%%%%%%%%%%%%%%%%%%%

%%%%%THIS DEFINES THE STYLE SUBSECTIONS, USING PACKAGE \usepackage{sectsty}%%%%%%%%%%%%%%%%
\subsectionfont{\normalsize}
%%%%%%%%%%%%%%%%%%%%%%%%%%%%%%%%%%%%%%%%%%%%%%%%%%%%%%%%%%%%%%%

\begin{center}
{\Large \textbf{On the constants for some fractional \\ ~\\ Gagliardo-Nirenberg and Sobolev inequalities}}
\end{center}
\vspace{0.4truecm}
\begin{center}
{\large
Carlo Morosi${}^1$, Livio Pizzocchero${}^2$} \\
\vspace{0.5truecm}
{\footnotesize
${}^1$ Dipartimento di Matematica, Politecnico di
Milano, \\ P.za L. da Vinci 32, I-20133 Milano, Italy \\
e--mail: carlo.morosi@polimi.it \\
${}^2$ Dipartimento di Matematica, Universit\`a di Milano\\
Via C. Saldini 50, I-20133 Milano, Italy\\
and Istituto Nazionale di Fisica Nucleare, Sezione di Milano, Italy \\
e--mail: livio.pizzocchero@unimi.it \\
}
\end{center}
\begin{abstract} \noindent
We consider the inequalities of Gagliardo-Nirenberg and Sobolev
in $\reali^d$, formulated in terms of the Laplacian $\Delta$
and of the fractional powers $\Di{n} := \sqrt{-\Delta}^{~n}$ with real
$n \geqs 0$; we review known facts and present novel results
in this area. After illustrating the equivalence between
these two inequalities and the relations between the corresponding
sharp constants and maximizers, we focus the attention on the $\LL^2$
case where, for all sufficiently regular $f : \reali^d
\vain \complessi$, the norm $\| \Di{j} f \|_{\LL^r}$
is bounded in terms of $\| f \|_{\LL^2}$ and $\| \Di{n} f \|_{\LL^2}$,
for $1/r = 1/2 - (\te n - j)/d$, and suitable
values of $j,n,\te$ (with $j,n$ possibly noninteger).
In the special cases $\te = 1$ and $\te = j/n + d/2 n$ (i.e., $r = + \infty$),
related to previous results of Lieb
and Ilyin, the sharp constants
and the maximizers can be found explicitly; we
point out that the maximizers can be expressed in terms
of hypergeometric, Fox and Meijer functions.
For the general $\LL^2$ case, we present two kinds of upper bounds
on the sharp constants: the first kind is suggested by the literature,
the second one is an alternative proposal of ours, often
more precise than the first one. We also
derive two kinds of lower bounds.
Combining all the available upper and lower bounds,
the Gagliardo-Nirenberg and Sobolev sharp constants
are confined to quite narrow intervals.
Several examples are given.
\end{abstract}
\textbf{Key words:} fractional Gagliardo-Nirenberg and Sobolev inequalities.
\vskip 0.2cm \noindent
\textbf{AMS 2000 Subject classifications:} 46E35, 26D10, 26A33.
\vskip 0.2cm \noindent
\vfill \eject \noindent
\section{Introduction}
\label{intro}
In this paper we work in $\reali^d$, using the Laplacian $\Delta$,
the operator $D := \sqrt{-\Delta}$ and its powers
$\Di{n}$ with real exponent $n \geqs 0$.
In the sequel $f$ stands for a complex-valued
function on $\reali^d$, with suitable regularity properties.
\par
We consider the embedding inequalities of
Gagliardo \cite{Gag}, Nirenberg \cite{Nir} and Sobolev \cite{Sob}.
The terms \textsl{Gagliardo-Nirenberg inequality} and \textsl{Sobolev inequality} are used
to indicate, respectively, the inequalities
({\footnote{The association
of the cited authors to either \rref{igag} or \rref{isob} is to some
extent conventional; in particular, the cited papers of Gagliardo consider
mainly the inequality \rref{isob}. However, these historical aspects are not relevant for our
purposes.}})
\beq \| \Di{j} f \|_{\LL^r} \leqs G \, \| f \|^{1 - \te}_{\LL^p} \| \Di{n} f \|^\te_{\LL^q}
\qquad \Big(\,{1 \over r} =  {1 -\te \over p} +
{\te \over q} - {\te n - j \over d} \Big), \label{igag} \feq
\beq \| \Di{j} f \|_{\LL^r} \leqs S \, ( \| f \|^t_{\LL^p} + \| \Di{n} f \|^t_{\LL^q})^{1/t}
\qquad (r~ \mbox{as in \rref{igag}})~,
\label{isob} \feq
holding if the parameters $p, q, j, n, \te, t$ fulfill appropriate conditions.
Here and in the sequel, $\LL^p$ is the usual space $L^p \equiv L^p(\reali^d, \complessi)$
for $p \in [1,+\infty)$, while $\LL^\infty$ is the subspace of $C(\reali^d, \complessi)$
made of the functions vanishing at infinity, with the usual sup norm (see the forthcoming Eq.\,\rref{linfo},
and the related comments). \par
The inequalities \rref{igag} \rref{isob} are found to be equivalent via appropriate
scaling considerations. We are especially interested in their $\LL^2$ versions
which are obtained setting $p=q=t=2$, and read
\beq \| \Di{j} f \|_{\LL^r} \leqs G \, \| f \|^{1 - \te}_{\LL^2} \| \Di{n} f \|^\te_{\LL^2}
\qquad \Big({1 \over r} = {1 \over 2} - {\te n - j\over d}  \Big)~,
\label{igagg} \feq
\beq \| \Di{j} f \|_{\LL^r} \leqs S \, \sqrt{\| f \|^2_{\LL^2} + \| \Di{n} f \|^2_{\LL^2}}
\qquad (r~\mbox{as in \rref{igagg}})~.
 \label{isobb} \feq
They hold under suitable conditions on $j,n,\te$, given in the forthcoming
Eq.\,\rref{unoo} and here anticipated: \beq {~} \hspace{-0.4cm}
0 \leqs \te \leqs 1\,, \quad 0 \leqs n,j < + \infty\,, \quad 0 \leqs \te n - j \leqs {d \over 2}\,,
\quad \te \neq 1\, \mbox{if $n = j + \dd{d \over 2}$}\,; \label{unooi} \feq
we write $G(j,n,\te)$ and $S(j,n,\te)$ for the
sharp constants of \rref{igagg} and \rref{isobb}, respectively. The aims of
this paper are as follows.
\begin{itemize}
\item To summarize some basic facts related to Eqs.\,\rref{igag} \rref{isob},
from the definition of the fractional Laplacian via Lizorkin distributions
and Fourier transform to the derivation of the equivalence between \rref{igag} and \rref{isob}
via scaling considerations. These generalities prepare the analysis
of the $\LL^2$ inequalities \rref{igagg} \rref{isobb}.
\item To review some special cases in which the sharp constants
and some maximizers for either \rref{igagg} or \rref{isobb} have been
determined in the literature; by the general equivalence mentioned
before, any result of this kind for one of the two inequalities can
be converted into a result for the other one.
\item To point out a fact unnoticed in the literature, i.e., that
the maximizers of the special cases mentioned before can be expressed in terms
of hypergeometric functions, Fox $H$-functions or Meijer $G$-functions.
\item To derive by more or less known strategies some
upper bounds for the sharp constants of the general inequalities \rref{igagg} and \rref{isobb}.
\item To propose some lower bounds for these sharp constants
and compare
them with the above mentioned upper bounds. This analysis is
performed in a fully quantitative way and ultimately
confines the unknown sharp constants of \rref{igagg} and \rref{isobb}
to rather small intervals, as shown by several numerical examples.
\end{itemize}
Let us describe with more details the organization of the paper.
In Section \ref{prelim} we fix our standards for some basic spaces of functions
or distributions in $\reali^d$; moreover, we recall how to define
the fractional Laplacians $\Di{n}$ via Fourier transform, in the framework of Lizorkin distributions. \par
In Section \ref{secgen} (and in the related Appendix \ref{appega})
we present formally the Gagliardo-Nirenberg and Sobolev inequalities \rref{igag} \rref{isob}
in their natural functional setting.
Moreover, using some scaling considerations we derive the equivalence between
\rref{igag} and \rref{isob}, and make explicit the relations
between the corresponding sharp constants and maximizers. This idea has
been used in the literature in many special cases \cite{Agu} \cite{Lie4};
our general formulation is contained in Proposition
\ref{soga}. \par
The results of Section \ref{secgen} are stated
on purely logical grounds, independently of the actual
conditions on the parameters $p, q, j, n, \te, t$ for
the validity of \rref{igag} \rref{isob}; these conditions
are reported in Section \ref{secfour},
with an appropriate bibliographical support.
\par
The remaining Sections \ref{secelledue}, \ref{seces} (and the
related Appendix \ref{appe}) form the longest part of the paper,
entirely devoted to the $\LL^2$ inequalities \rref{igagg} \rref{isobb};
hereafter we describe in detail the contents of these sections that provide, amongst else,
simple and self-consistent proofs of \rref{igagg}
or \rref{isobb} (for $j,n,\te$ as in Eq.\,\rref{unooi}),
independently of the general results on the validity of \rref{igag} \rref{isob}.
\begin{itemize}
\item[(i)] In subsection \ref{jten} we begin our discussion of the $\LL^2$
inequalities analyzing the elementary case $j = \te n$, i.e., $r=2$.
In this case the Gagliardo-Nirenberg inequality
\rref{igagg} is reduced (via Fourier transform) to the H\"older
inequality: the Gagliardo-Nirenberg sharp constant is readily
found to be $1$, and by the equivalence between
\rref{igagg} and \rref{isobb} one also obtains the Sobolev sharp constant.
\item[(ii)] In subsection \ref{subsuno} we pass to
the (much less elementary) case $\te=1$. Eq.\,\rref{igagg} (and, more generally, Eq.\,\rref{igag})
with $\te=1$ is equivalent to the so-called
Hardy-Littlewood-Sobolev inequality concerning
convolution with a power of the radius $|x|$
($x \in \reali^d$) \cite{HL} \cite{Sob} (see, e.g., \cite{Mizo}
for a more up-to-date reference on these issues).
The sharp constant and the
maximizers of the Hardy-Littlewood-Sobolev
inequality in the $\LL^2$ case have been
found by Lieb \cite{Lie0}, and these
results can be used in relation
to the $\te=1$ case of \rref{igagg}; this was done in
explicit terms in \cite{Cot} (with the additional assumption $j=0$). A review of these facts
is presented (extending the considerations
of \cite{Cot} to nonzero values of $j$). The maximizer
derived along these lines is a radial function and
can be represented as the inverse Fourier transforms
of a Bessel-type function; when expressed in terms of the ordinary space
variables, it is an elementary function
for $j=0$ and a hypergeometric-type function for suitable values of $j,n$.
\item[(iii)] In subsection \ref{subsdue}, following a path
indicated by \cite{Haj}, we combine the previous results on
the cases $j = \te n$ and $\te=1$ to derive the Gagliardo-Nirenberg inequality
\rref{igagg} and to obtain upper bounds $\Gp(j,n,\te)$
on its sharp constants in an ``almost general case''; this covers
all choices of $j, n, \te$ in Eq.\,\rref{unooi} except
$\te = j/n + d/2 n$ (i.e., $r = + \infty$).
The upper bounds obtained by this strategy coincide
with the sharp constant if $\te = 1$, and diverge
if $\te$ approaches the exceptional value
$j/n + d/2 n$.
By the equivalence between the Gagliardo-Nirenberg
and the Sobolev inequalities, the above results can be
rephrased making reference to \rref{isobb}; in particular we have
upper bounds $\Sp(j,n,\te)$ on the corresponding sharp constants
(again, coinciding with the sharp constant for $\te=1$
and diverging for $\te \vain j/n + d/2 n$).
\item[(iv)] Subsection \ref{substre} focuses on the previously excluded
case $\te = j/n + d/2 n$ (i.e., $r = + \infty$), requiring
a different strategy. For $j=0$ this case was analyzed
by Ilyin \cite{Ily} who
derived via Fourier transform the Sobolev inequality \rref{isobb}, found the sharp constant
and the maximizer, and then pointed out the implications of these results on
the Gagliardo-Nirenberg inequality \rref{igagg}.
In this subsection we propose a similar derivation for \rref{isobb},
giving the sharp constant and a maximizer
holding for arbitrary $j$; the consequences for \rref{igagg} are indicated.
Our maximizer is radial and agrees for $j=0$ with the one of Ilyin;
as in \cite{Ily}, it can be represented as the inverse Fourier transform
of an elementary function. We also derive its expressions in terms of
the space variables, using the Fox $H$-function or
the Meijer $G$-function (whose definitions are reviewed in Appendix \ref{appe}).
\item[(v)] \label{panti}In subsection \ref{subsqua} we propose
a second, ``almost general'' approach to the $\LL^2$ inequalities
\rref{igagg} \rref{isobb}, alternative to the one of item (iii);
this covers all choices of $j,n,\te$ in Eq.\,\rref{unooi}, except the special
case $\te=1$ of item (ii). We use again Fourier transform arguments,
inspired by a previous work of ours \cite{imb} on a variant
of \rref{isobb}; these yield an alternative proof of the Sobolev
inequality \rref{isobb}, giving upper bounds
$\Spp(j,n,\te)$ for its sharp constants. These
bounds coincide with the Sobolev sharp constant
for $\te = j/n + d/2 n$, and diverge for $\te \vain 1$. Due to the equivalence between
the Sobolev and Gagliardo-Nirenberg inequalities,
these results on \rref{isobb} yield upper bounds
$\Gpp(j,n,\te)$ on the sharp constants
of \rref{igagg}, with an analogous behavior for $\te =j/n + d/2 n$
and for $\te \vain 1$. \parn
\item[(vi)] As emphasized in subsection
\ref{subsqui}, the upper bounds $\Gp, \Sp$ of item (iii) and
$\Gpp,\Spp$ of item (v) cover jointly the general $\LL^2$ inequalities
\rref{igagg} \rref{isobb}. The upper bounds
$\Gp, \Sp$ are expected to be better (i.e., smaller) for $\te$ close to $1$,
while $\Gpp, \Spp$ are expected to be better for
$\te$ close to $j/n + d/2 n$. To be more specific one
can make a direct comparison of the numerical values of these
bounds, an issue that is treated in Section \ref{seces}
for $d=1, 2, 3$ and some test values of $j, n, \te$.
\item[(vii)] As a final step in the theoretical
investigation of the $\LL^2$ inequalities
\rref{igagg} \rref{isobb}, in subsection \ref{subssei} we derive
lower bounds for their sharp constants.
We obtain two types of lower bounds
$\Gm(j,n,\te)$,
$\Sm(j,n,\te)$
and $\Gmm(j,n,\te)$,
$\Smm(j,n,\te)$,
derived substituting
``trial functions'' of two kinds
for $f$ in \rref{igagg} or in \rref{isobb}.
The $-$ bounds hold under the general $\LL^2$ conditions \rref{unooi},
while the $--$ bounds require some limitations
for $j, n, \te$ (see Eq.\,\rref{unoi}).
\item[(viii)] In Section \ref{seces}, as examples we
write explicitly the sharp constants
and maximizers of \rref{igagg} \rref{isobb}
for $d=1,2,3$, $\te = j/n + d/2n$ (see item (iv))
and some choices of $j,n$: see Tables I, II.
Moreover, we give the numerical values
of the upper bounds $G_{+}, G_{++}$
and of the lower bounds $G_{-}, G_{--}$
on the sharp constants of \rref{igagg},
for $d=1,2,3$ and several choices of
$j,n,\te$: see Table III. (In these examples
the $+$ bounds are better than the $++$ bounds for
$\te \simeq 1$, as expected; however,
this occurs only for $\te$
\textsl{very} close to $1$.)
In these numerical tests the best lower and upper bounds
are close together, thus confining the sharp constants
$G(j,n,\te)$ to narrow intervals. By the equivalence
between the Gagliardo-Nirenberg and the Sobolev
inequality, it would be easy to produce similar
numerical results for the bounds $\Sp, \Spp, \Sm, \Smm$
on the sharp constants of \rref{isobb}.
\end{itemize}
In the above description of the contents of the paper,
we have indicated their connections with the existing
literature to the best of our knowledge; more details
are given in the sequel. We hope that the survey
of known results presented here in a unified language,
and our contributions mentioned in items (v)-(viii), will
allow a more complete understanding of the $\LL^2$
Gagliardo-Nirenberg-Sobolev inequalities.
\section{Some preliminaries}
\label{prelim}
\textbf{Notations.}
Throughout the paper we intend
\beq 0^0 := (+\infty)^0 := 1 ~; \label{nota1} \feq
\beq (a^t + b^t)^{1/t} := \max(a,b) \qquad \mbox{for $a, b \in [0, +\infty)$
and $t = + \infty$}~. \label{nota2}
\feq
We work in the Euclidean space $\reali^d$ for a fixed
space dimension $d \in \{1,2,3,...\}$. We often write $\x$
for the identity map $\reali^d \vain \reali^d$, $x \mapsto x$,
and $|\x|$ for the map $\reali^d \vain [0,+\infty)$,
$x \mapsto |x|$. When dealing with the Fourier transform,
a typical wave vector in $\reali^d$ is indicated
with $k$; in this framework, the identity map and the
Euclidean norm of $\reali^d$ are indicated with $\k$ and $|\k|$.
\par
Given two complex topological vector spaces $E, F$ we
say that $E$ is continuously embedded in $F$, and
write $E \hookrightarrow F$, if $E$ is a vector
subspace of $F$ and the natural inclusion of $E$ into
$F$ is continuous.
\vskip 0.2cm\noindent
\textbf{Standard terminology for Banach space inequalities.}
Let us consider two complex Banach spaces $X,Y$ with norms $\|~\|_X$,
$\|~\|_Y$. Moreover, assume we are given a linear map
$\EF: X \vain Y'$
where $Y'$ is vector space containing $Y$ as
a vector subspace. We often consider statements with the following
structure:
\beq \EF X \subset Y~, \qquad \|\EF f \|_Y \leqs C \|f \|_X
\qquad \mbox{for all $f \in X$ and some $C \in [0,+\infty)$}.
\label{gineq} \feq
Any such statement is referred to as an \textsl{inequality};
of course \rref{gineq} indicates that $\EF$ is continuous from $X$ to $Y$ and that,
if $\EF$ is the identity, there is a continuous embedding $X \hookrightarrow Y$.
The \textsl{sharp constant} $C_s$ for the inequality \rref{gineq} is
the inf of the set of the constants $C \in [0,+\infty)$ which fulfill it;
this inf is in fact a minimum. A \textsl{maximizer}
for \rref{gineq}, if it exists, is a nonzero element $f \in X$ such that
$\|\EF f \|_Y = C_s \|f \|_X$. \par
It is clear that \rref{gineq} holds if and only if
$\EF X \subset Y$ and the ratio $\|\EF f \|_Y/\|f \|_X$
is bounded for $f$ ranging in $X \setminus \{0 \}$.
If this happens, the sharp constant can be expressed as
$C_s = \sup_{f \in X, f \neq 0} {\| \EF f \|_Y / \| f \|_X}$;
an element $f \in X \setminus \{0 \}$ is a maximizer
for \rref{gineq} if and only if it is a maximum point
for the above ratio.
Of course, any transformation leaving invariant this
ratio sends maximizers into maximizers; the simplest
example is the map $f \mapsto K f$ where $K \in \complessi
\setminus \{0 \}$.
\salto
\vskip 0.2cm\noindent
\textbf{Some spaces of functions and distributions on $\reali^d$.}
As usual, we say that a function $\phi: \reali^d \vain \complessi$
is rapidly decreasing if $(1 + |\x|)^N \phi$
is bounded for all $N \in \naturali$.
\par
We employ the standard symbol
$\SS(\reali^d, \complessi) \equiv \SS$ for the Schwartz space,
formed by the $C^\infty$ functions $\varphi : \reali^d \vain
\complessi$ rapidly decreasing with all
their derivatives; this space is equipped with the Fr\'echet topology induced by the seminorms
$p_{N i_1,...,i_m}(\varphi) := \sup_{x \in \reali^d}
(1 + |x|)^N |\partial_{i_1,...,i_m} \varphi (x)|$,
where $N,m \in \naturali$ and $i_1,...,i_m \in \{1,...,d\}$. \par
We consider the Fourier transform
\beq \FF : \SS \vain \SS~, \qquad \varphi \mapsto \FF \varphi~, \feq
normalized so that
$\FF \varphi(k)= (2 \pi)^{-d/2} $ $\int_{\reali^d} d x~e^{-i k \sc x}
\varphi(x)$; this is a linear homeomorphism
of $\SS$ into itself.
The dual space $\SS'(\reali^d, \complessi)
\equiv \SS'$ is the standard space of tempered distributions, and
one can extend the Fourier transform to a map of $\SS'$ into itself.
\par
Some less conventional spaces of
distributions and the related Fourier transforms are more interesting in relation to fractional
differential calculus \cite{Sam}. Their construction relies
on the Lizorkin-type spaces of test functions
\beq \Phi(\reali^d, \complessi) \equiv \Phi := \{ \phi \in \SS~|~\int_{\reali^d} d x \, x_{i_1} ... x_{i_m}
\phi(x) = 0~\mbox{for all $i_1,...,i_m$}\, \}~, \feq
\beq \Psi(\reali^d, \complessi) \equiv \Psi :=  \{ \psi \in \SS~|~\partial_{i_1,...,i_m} \psi(0) = 0
~\mbox{for all $i_1,...,i_m$} \}~.
\feq
(In both cases, ``for all $i_1,...,i_m$'' means
``for all $m \in \naturali$ and $i_1,...,i_m \in \{1,...,d\}$''.
Obviously enough, $x_{i_1} ... x_{i_m} := 1$ if $m=0$.
Functions $\psi \in \Psi$ are usually written as $\psi : \reali^d
\vain \complessi$, $k \mapsto \psi(k)$.)
\par
$\Phi, \Psi$ are closed vector subspaces of $\SS$, and thus
are Fr\'echet spaces with the induced topology.
One readily checks that
\beq \FF \Phi = \Psi~, \label{check} \feq
and that $\FF$ is a linear
homeomorphism between $\Phi$ and $\Psi$. The Lyzorkin-type
distribution spaces
$\Phi'(\reali^d) \equiv \Phi'$ and $\Psi'(\reali^d) \equiv \Psi'$
are the dual spaces of $\Phi$ and $\Psi$,
equipped with their weak topologies
({\footnote{Of course, there are continuous
linear maps $\SS' \vain \Phi'$, $f \mapsto f \restriction \Phi$
and $\SS' \vain \Psi'$, $f \mapsto f \restriction \Psi$;
the kernels of these maps are, respectively, the space $\PPp(\reali^d) \equiv \PPp$
of the polynomial functions on $\reali^d$ and the space
$\DDd(\reali^d) \equiv \DDd$
of the finite linear combinations of the Dirac delta at
the origin and of its derivatives. Therefore, there are
linear homeomorphisms $\SS'/\PPp \vain \Phi'$ and
$\SS'/\DDd \vain \Psi'$ \cite{Sam}.}}).
One can define a Fourier transform
\beq \FF : \Phi' \vain \Psi',
~~\mbox{$\langle \FF f, \psi \rangle := \langle f, \FF^{-1} \psi_R \rangle $ ~ for $f \in \Phi'$,
$\psi \in \Psi$, $\psi_R(k) := \psi(-k)$}~;
\label{lizo} \feq
this is a linear homeomorphism between
$\Phi'$ and $\Psi'$
({\footnote{Formally, one has $\FF^{-1} \psi_R = \FF \psi$. This
identity, even though correct, is somehow misleading;
it can be stated just because the spaces of wave vectors
$k$ and position vectors $x$ are both identified
with $\reali^d$.}}). \par
For $j\in \{1,...,d\}$ the map $\phi \mapsto \partial_j \phi$ is
linear and continuous from $\Phi$ into itself; we define the distributional
derivative $\partial_j : \Phi' \vain \Phi'$ by the usual procedure
\textsl{\`a la Schwartz}, i.e., setting
$\langle \partial_j f, \phi \rangle := - \langle f, \partial_j \phi \rangle$
for all $f \in \Phi', \phi \in \Phi$. \par
To go on, let us consider any function $\eta \in C^\infty(\reali^d \setminus \{0 \}, \complessi)$
such that, for all $n \in \naturali$ and $j_1,...,j_n \in \{1,...,d\}$ one has
$|\partial_{j_1,...,j_n} \eta(k)| = O(|k|^{-M})$ for $k \vain 0$
and $|\partial_{j_1,...,j_n} \eta(k)| = O(|k|^N)$ for $k \vain \infty$,
for suitable real exponents $N, M$ (depending on $n,j_1,...,j_n$).
For $\psi \in \Psi$
the product $\eta \psi$ is clearly $C^\infty$ on $\reali^d \setminus \{ 0 \}$
and possesses a unique $C^\infty$ extension to $\reali^d$, vanishing at
the origin and rapidly decreasing with all its derivatives:
in this sense, we have $\eta \psi \in \Psi$. The map
\beq \Psi \vain \Psi~, \qquad \psi \mapsto \eta \psi \label{mapeta} \feq
is linear and continuous. This fact can be used to introduce a continuous
linear map
\beq \Psi' \vain \Psi'~, \qquad g \mapsto \eta g~ \label{mapetaa} \feq
putting $\langle \eta \,g, \psi \rangle := \langle g, \eta \,\psi \rangle$
for all $\psi \in \Psi$. \textsl{A fortiori}, this construction works
in the special case of a function $\eta \in C^{\infty}(\reali^d, \complessi)$
with the previously mentioned behavior for $k \vain \infty$.
In particular, let us choose for $\eta$ the function
$\k_j : \reali^d \vain \reali$, $k \mapsto k_j$, for some
$j \in \{1,...,d\}$. This gives a map
$\Psi' \vain \Psi'$, $g \mapsto \k_j \,g$ and it turns out that
$\FF(\partial_j f) = i \k_j \, \FF f$ for all $f \in \Phi'$.
\par
For each $p \in [1,+\infty]$, we write $L^p$ for the usual space $L^p(\reali^d,
\complessi)$ and $\| ~\|_{L^p}$ for its
norm; moreover, we define as follows the Banach space $\LL^p$ and
its norm:
\beq \LL^p := L^p~,~~ \|~\|_{\LL^p}:= \|~ \|_{L^p}~~\mbox{if $p \in [1,+\infty)$}~,
\label{ellepi} \feq
\beq \LL^{\infty} := C_0(\reali^d, \complessi) = \{ f \in C(\reali^d, \complessi)~|~\lim_{x \vain \infty}
f(x) = 0 \}~,~ \| f \|_{\LL^\infty} := \sup_{x \in \reali^d} |f(x) |~. \label{linfo} \feq
$\LL^\infty$ is a closed subspace of $L^\infty$, and $\|~\|_{\LL^\infty}$ is just
the restriction to this subspace of $\|~\|_{L^\infty}$
({\footnote{Considering the family
\rref{ellepi} \rref{linfo} is not unusual, see e.g.
\cite{Dev} \cite{Fic} where similar families are
employed in different situations; admittedly,
the notation $\LL^p$ is not standard. For the purposes of the present
work, the space $\LL^\infty$ is much more interesting than $L^\infty$.}}).
For any $p \in [1,+\infty]$, we clearly have
\beq \Phi, \Psi
\hookrightarrow \LL^p~. \label{clear1} \feq
Again for $p \in [1,+\infty]$, given $f \in \LL^p$ we can define a continuous linear form
$\langle f,~\rangle$ on $\Phi$
setting $\langle f,\phi \rangle := \int_{\reali^d} d x \, f(x) \phi(x)$ for
$\phi \in \Phi$. The map $f \in \LL^p \mapsto \langle f,~\rangle \in \Phi'$
is a continuous linear injection ({\footnote{According to
the footnote before Eq.\,\rref{lizo}, the kernel of this map is the intersection of $\LL^p$ with the space
$\PPp$ of polynomials on $\reali^d$; clearly, $\LL^p \cap \PPp = \{ 0 \}$ for all
$p \in [1,+\infty]$. To obtain this result for $p = \infty$, it is crucial to
define $\LL^{\infty}$ in terms of functions vanishing at infinity.}}); thus, up to a natural identification (to be employed from
now on), for all $p \in [1,+\infty]$ we have
\beq \LL^p \hookrightarrow \Phi'~.
\label{clear2} \feq
In a similar way, for $p \in [1,+\infty]$ and $f \in \LL^p$
we can define a continuous linear form on $\Psi$ setting $\langle f,\psi
\rangle := \int_{\reali^d} d k \, f(k) \psi(k)$ for all $\psi \in \Psi$. We
have a continuous linear injection $f \in \LL^p \mapsto \langle f,~\rangle \in
\Psi'$, so that
\beq \LL^p \hookrightarrow \Psi'~. \label{clear2bis} \feq
Of course, many other topological vector spaces of complex, measurable functions on $\reali^d$
can be continuously embedded in $\Phi'$ or $\Psi'$ using the previous prescriptions
to identify a function $f$ with a continuous linear form on $\Phi$ or $\Psi$.
For example, we have the embeddings
\beq \Phi \hookrightarrow
\Phi', \quad \Psi \hookrightarrow \Psi' \label{clear3} \feq
(which can be seen as trivial consequences of Eqs.\,\rref{clear1}-\rref{clear2bis};
using the second of these relations, we can view
the map \rref{mapeta} as a restriction of the map \rref{mapetaa}).
\vskip0.2cm\noindent
\textbf{Hausdorff-Young inequality.} Let
\beq p \in [1,2]~,
\qquad p'~ \mbox{such that}~1/p + 1/p' = 1~; \label{pip} \feq
it is well known
that ({\footnote{Note that \rref{known} holds as well for $p=1$, $p' =
+\infty$, using our definition \rref{linfo} for $\LL^\infty$; in fact,
according to the Riemann-Lebesgue lemma \cite{Lie2}, the Fourier
transform of an $\LL^1$ function is a continuous function vanishing at
infinity.}})
\beq {~} \hspace{-0.5cm} f \in \LL^p~ \Rightarrow \FF f \in \LL^{p'}, \|
\FF f \|_{\LL^{p'}} \leqs C_p \| f \|_{\LL^p}\,, ~C_p := {1 \over (2 \pi)^{d/p
- d/2}}~ \Big[ {(1/p')^{1/p'} \over ({1/p})^{1/p}} \Big]^{d/2}. \label{known}
\feq
The inequality in \rref{known} is the familiar \textsl{Hausdorff-Young inequality}; the constant
$C_p \equiv C_{p d}$ is known to be sharp (see
\cite{Lie1}, \cite{Lie2} Chapter 5 and references therein; our expression for
$C_{p}$ differs by a factor from the one in \cite{Lie2} due to a different
normalization for the Fourier transform). A Hausdorff-Young inequality holds as
well for the inverse Fourier transform: with $p, p'$ as in \rref{pip} and $C_p$
as in \rref{known},
\beq {~} \hspace{-0.5cm} g \in \LL^p~ \Rightarrow \FF^{-1}
g \in \LL^{p'},\, \| \FF^{-1} g \|_{\LL^{p'}} \leqs C_p \, \| g \|_{\LL^p}~;
\label{knowy} \feq
the constant $C_p$ of \rref{known} is sharp as well for the
formulation \rref{knowy}.
\vskip 0.2cm \noindent
\textbf{Fractional Laplacian.}
Let $s \in \reali$, and consider the $C^\infty$ function $|\k|^s :\reali^d
\setminus \{0 \} \vain \reali$ (having a $C^\infty$ extension to $\reali^d$ if
and only if $s$ is a nonnegative, even integer). Making reference to
Eqs.\,\rref{mapeta} \rref{mapetaa} and to the related comments, we can define
two maps
\beq \Psi \vain \Psi, \quad \psi \mapsto |\k|^s \psi~; \qquad \Psi'
\vain \Psi', \quad g \mapsto |\k|^s g~; \label{kas} \feq
both of them are
linear homeomorphisms, with inverses corresponding to multiplication by
$|\k|^{-s}$. \par Composing the second map \rref{kas} with $\FF^{-1} : \Psi'
\vain \Phi'$ and $\FF : \Phi' \vain \Psi'$ we obtain the linear homeomorphism
\beq \Di{s} :\Phi' \vain \Phi'~, \qquad f \mapsto \Di{s} f := \FF^{-1}( |\k|^s
\FF f)~, \label{dis} \feq
with inverse $\Di{-s}$; of course, $\Di{s} \Di{t} =
\Di{s+t}$ for any real $s,t$. Taking into account
Eqs.\,\rref{check}\,\rref{kas} and the comments which accompany them, we see
that
\beq \Di{s} \Phi = \Phi \label{seeth} \feq
(identifying $\Phi$ with a
subspace of $\Phi'$, see \rref{clear3}) and that $\Di{s}$ is a homeomorphism of
$\Phi$ into itself. Noting that the usual Laplacian $\Delta$ fulfills $-
\Delta f = \FF^{-1}( |\k|^2 \FF f)$, we see that there
would be good reasons to write
\beq \Di{s} \equiv \Lap{s}~; \feq
however, the
symbol $\Lap{s}$ is never employed in the sequel. \salto
\textbf{Convolution, and representation of $\boma{\Di{-n}}$ for $\boma{n \in (0, d)}$}. Given two
measurable functions $g, f: \reali^d \vain \complessi$, the convolution $g * f
: \reali^d \vain \complessi$, $(g * f)(x) := \int_{\reali^d} d y \,g(x-y) f(y)$
can be defined if the previous integral exists for almost all $x$. This
happens, in particular, if $g := 1/|\x|^{d-n}$ and $f \in \LL^q$ with $0 < n <
d$ and $1 \leqs q < d/n$ (see \cite{Stein}, Chapter V, Theorem 1). For $n,q$
and $f$ as above, one has
\beq \Di{-n} f = {\AZ_{n} \over |\x|^{d-n}}* f~,
\qquad
 \AZ_{n} \equiv \AZ_{d n} := {\Gamma(d/2 - n/2) \over 2^n \pi^{d/2} \Gamma(n/2)}~.
\label{223} \feq
(see again \cite{Stein}, Chapter V, Lemma 1 or \cite{Mizo}, Chapters 2 and 7)
({\footnote{The proof of \rref{223} can be sketched as follows. By definition
\hbox{$\Di{-n} f$} \hbox{$= \FF^{-1}(|\k|^{-n} \FF f)$}. But
\hbox{$\FF^{-1}(u v)$} \hbox{$= (2 \pi)^{-d/2} \FF^{-1} u * \FF^{-1} v$}
under suitable conditions on the functions $u,v$;
applying this with $u = |\k|^{-n}$ and $v = \FF f$
one gets $\Di{-n} f = (2 \pi)^{-d/2} (\FF^{-1} |\k|^{-n})*f$.
Finally, one shows that $(2 \pi)^{-d/2} \FF^{-1} |\k|^{-n} =
\AZ_{n} \,|\x|^{-d+n}$. (These manipulations are related
to the following prescriptions: for each $s \in \reali$
and $a \in (-d,0)$, the functions $|\k|^s$ and $|\x|^a$
are identified with elements of $\Psi'$ and $\Phi'$,
respectively, setting $\langle |\k|^s, \psi \rangle
:= \int_{\reali^d} d k \, |k|^s \psi(k)$ and
$\langle |\x|^a, \phi \rangle := \int_{\reali^d} d x \, |x|^a \phi(x)$.
The conditions on $s, a$ and the features of the test
functions spaces $\Psi$, $\Phi$ ensure
that the previous integrals converge and
depend continuously on $\phi \in \Phi$, $\psi \in \Psi$).}}).
\salto
\textbf{Translations and scalings}. Let $a \in \reali^d$, $\lambda \in
(0,+\infty)$; we consider the translation $\reali^d \vain \reali^d$, $x \mapsto
x + a$ and the scaling transformation $\reali^d \vain \reali^d$, $x \mapsto
\lambda x$. \par If $\phi$ is a complex-valued function on $\reali^d$, the
$a$-translated function $\phi^a$ and the $\lambda$-rescaled function
$\phi_{\lambda}$ are the compositions of $\phi$ with the previously mentioned
transformations; thus
\beq \phi^a, \phi_{\lambda} : \reali^d \vain \complessi~,
\qquad \phi^a(x) := \phi(x+a), \quad \phi_{\lambda}(x) := \phi(\lambda x)~.
\label{scal} \feq
One proves that $\phi \in \Phi$ $\Rightarrow$ $\phi^a,
\phi_{\lambda} \in \Phi$. For $f \in \Phi'$, $a \in \reali^d$ and $\lambda \in
(0,+\infty)$ we define $f^a, f_{\lambda} \in \Phi'$ by
\beq \langle f^a, \phi
\rangle := \langle f, \phi^{-a} \rangle~, \qquad \langle f_{\lambda}, \phi
\rangle := \lambda^{-d} \langle f, \phi_{\lambda^{-1}} \rangle \qquad \mbox{for
all $\phi \in \Phi'$}~. \feq
Let $s \in \reali$, $p \in [1,+\infty]$. Clearly,
the operators $\Di{s}$ and the spaces $\LL^p$ are translation invariant ($f \in
\Phi' \Rightarrow \Di{s} f^a = (\Di{s} f)^a$; $f \in \LL^p \Rightarrow f^a \in
\LL^p$, $\| f^a \|_{\LL^p} = \| f \|_{\LL^p}$). As for the scaling by $\lambda
\in (0,+\infty)$, the following is proved by elementary means:
\beq f \in \Phi'
\Rightarrow \Di{s} f_{\lambda} = \lambda^s (\Di{s} f)_{\lambda}~; \feq
\beq f
\in \LL^p \Rightarrow f_{\lambda} \in \LL^p,\, \| f_{\lambda} \|_{\LL^p} =
\lambda^{-d/p} \| f \|_{\LL^p}~; \label{nonum} \feq
\beq f \in \Phi', \Di{s} f
\in \LL^p \Rightarrow \Di{s} f_{\lambda} \in \LL^p,\, \| \Di{s} f_{\lambda}
\|_{\LL^p} = \lambda^{s -d/p} \| \Di{s} f \|_{\LL^p} \label{19} \feq
(of course, the third statement is a consequence of the first two).
\salto
\textbf{Spaces of Riesz potentials.} Let $q \in [1, + \infty]$, $n \in [0,+
\infty)$. The \textsl{space of Riesz potentials} (or \textsl{homogeneous Sobolev
space}) of type  $q,n$ is
({\footnote{Most treatments of Riesz potential spaces do not consider
the case $q = + \infty$, see e.g. \cite{Gra2} \cite{Sam}; the same
comment can be done on the inhomogeneous Sobolev spaces introduced
hereafter. $\WW{q}{n}$
is often denoted with the alternative notation $I^n(\LL^q)$, in which $I^n$
stands for $\Di{-n}$, see \cite{Sam}.}})
\beq \WW{q}{n}(\reali^d) \equiv \WW{q}{n} := \{ f \in
\Phi'~|~ \Di{n} f \in \LL^q \}~. \label{riesz} \feq
By construction, we have a linear isomorphism $\WW{q}{n} \vain \LL^q$, $f \mapsto \Di{n} f$ with
inverse $\LL^q \vain \WW{q}{n}$, $h \mapsto \Di{-n} h$.  Due to this fact, $\WW{q}{n}$
is a Banach space with norm
\beq \| f \|_{\WW{q}{n}} := \| \Di{n} f
\|_{\LL^q}~. \feq
\salto
\textbf{Inhomogeneous Sobolev spaces.} Let $p,q \in [1, + \infty]$,
$n \in [0,+ \infty)$. The \textsl{inhomogeneous Sobolev space} of order $p,q,n$
is
\beq \ZZ{p}{q}{n}(\reali^d) \equiv \ZZ{p}{q}{n} := \LL^p \cap \WW{q}{n} = \{
f \in \LL^p~|~\Di{n} f \in \LL^q~ \}~. \label{wpqn} \feq
This is a Banach space
with respect to any one of the equivalent norms
\beq \| f \|_{\ZZ{p}{q}{n}\,|\, t} := ( \| f \|^t_{\LL^p} + \| \Di{n} f \|^t_{\LL^q})^{1/t}
\qquad  (t \in [1,+\infty])\,  \feq
(if $t=+\infty$, this definition must be understood following
\rref{nota2}) ({\footnote{The fact that $\ZZ{p}{q}{n}$ is a vector space and
the equivalence of all norms $\| ~ \|_{\ZZ{p}{q}{n} \,| \,t}$ are evident. To
prove completeness, let us consider a Cauchy sequence $(f_\ell)_{\ell \in
\naturali}$ in $\ZZ{p}{q}{n}$. Then $(f_\ell)_{\ell \in \naturali}$ and
$(\Di{n} f_\ell)_{\ell \in \naturali}$ are Cauchy sequences in $\LL^p$ and
$\LL^q$, respectively, so there are functions $g \in \LL^p$ and $h \in \LL^q$
such that $f_\ell \vain g$ in $\LL^p$ and $\Di{n} f_\ell \vain h$ in $\LL^q$.
But $\LL^p, \LL^q \hookrightarrow \Phi'$, so $f_\ell \vain g$ and $\Di{n}
f_\ell \vain h$ in $\Phi'$; by the continuity of $\Di{-n}$ on $\Phi'$ we also
infer $f_\ell \vain \Di{-n} h$ in $\Phi'$ so that $\Di{-n} h = g$, i.e., $h =
\Di{n} g$. In conclusion $f_\ell \vain g$ in $\LL^p$ and $\Di{n} f_\ell \vain
\Di{n} g$ in $\LL^q$, which means $f_\ell \vain g$ in $\ZZ{p}{q}{n}$.}}).
\salto
\textbf{Embedding and density statements.}
Let $p,q \in [1, + \infty]$, $n \in [0,+ \infty)$. It is evident that $\Phi \hookrightarrow
\ZZ{p}{q}{n} \hookrightarrow \WW{q}{n}$. Moreover $\Phi$ (and, consequently,
$\ZZ{p}{q}{n}$) is dense in $\WW{q}{n}$ if $q \neq 1,+\infty$ ({\footnote{To
prove this density statement, we consider an $f \in \WW{q}{n}$ and show that
$f$ is the limit of a sequence of functions $f_{\ell} \in \Phi$. Indeed,
$\Di{n} f \in \LL^q$ and it is known that $\Phi$ is dense in $\LL^q$ (see
\cite{Sam} page 41, Theorem 2.7),
so there is a sequence $(g_\ell)_{\ell \in \naturali}$ in $\Phi$ such that
$g_{\ell} \vain \Di{n} f$ in $\LL^q$.
Recalling that $\Di{n}$ is a one-to-one map of $\Phi$ into itself,
with inverse $\Di{-n}$, let us introduce the functions
$f_{\ell} := \Di{-n} g_{\ell} \in \Phi$;
then by construction $\Di{n} f_{\ell} = g_{\ell} \vain \Di{n} f$
in $\LL ^q$, which means $f_{\ell} \vain f$ in $\WW{q}{n}$.}}).
\section{The inequalities of Gagliardo-Nirenberg and So\-bo\-lev,
with their logical connections}
\label{secgen}
The aims of this section are: \parn
(i) to \emph{define} in formal terms the above inequalities; \parn
(ii) to point out their \emph{logical connections} (i.e., the fact that one of
them implies the other one). \par
This discussion is carried over independently of the \emph{validity conditions}
for the above inequalities, that are the subject of Section \ref{secfour}. We think
that the logical status of these inequalities has its own interest,
independently of the strategies that one can use to prove their
validity in more or less general situations.
\salto
\subsection{The Gagliardo-Nirenberg inequality}
\label{subgagliardo}
The definition given hereafter for this inequality involves certain
parameters $p$, $q$, $j$, $n$, $\te$; for the moment we put on the
parameters the minimal conditions ensuring
well definedness of both sides in the inequality and
certain scaling properties of general use.
The \emph{validity} of the inequality requires more
stringent conditions, see Section \ref{secfour}. \par
The minimal conditions are the following ones:
\beq 1 \leqs p,q \leqs + \infty~,
\quad 0 \leqs \te \leqs 1~, \quad 0 \leqs n,j < + \infty~, \label{min} \feq
$$ 0 \leqs {1 -\te \over p} +
{\te \over q} - {\te n - j \over d} \leqs 1~; $$
due to the last condition, there is a unique $r \equiv r(p,q;j,n,\te) \in [1,+\infty]$ such that
\beq {1 \over r} = {1 -\te \over p} +
{\te \over q} - {\te n - j \over d} ~. \label{due} \feq
\begin{prop}
\label{defgen}
\textbf{Definition.} Let $p,q,j,n,\te,r$ be as in \rref{min} \rref{due}.
The \textsl{Gagliardo-Nirenberg inequality} of order
$(p,q;j,n,\te)$
is the following statement:
\beq \ZZ{p}{q}{n} \subset \ZZ{p}{r}{j}~ \mbox{and} \label{gag} \feq
$$ \| \Di{j} f \|_{\LL^r} \leqs G \, \| f \|^{1 - \te}_{\LL^p} \| \Di{n} f \|^\te_{\LL^q}
\quad \mbox{for some $G \in [0,+\infty)$ and all $f \in \ZZ{p}{q}{n}$}. $$
Whenever this holds, the symbol $G(p,q;j,n,\te)$ indicates its
\textsl{sharp} constant.
\end{prop}
\begin{rema}
\textbf{Remarks.}
(i) The above inequality and its sharp
constant are related to the space
dimension $d$; so we should write, say,
$G_d(p,q;j,n,\te)$ for the sharp constant. For the sake of simplicity,
in the sequel $d$ is fixed and omitted from
most of our notations. \parn
(ii) If \rref{gag}
holds, for the related sharp constant we have the representation
\beq G(p,q;j,n,\te) = \sup_{f \in \ZZ{p}{q}{n} \setminus \{0 \}}
{\| \Di{j} f \|_{\LL^r} \over \| f \|^{1 - \te}_{\LL^p} \| \Di{n} f \|^\te_{\LL^q}}~.  \label{sec} \feq
(iii) The ratio in the right hand side
of Eq.\,\rref{sec} is invariant
under translations $f \mapsto f^a$ ($a \in \reali^d$, see Eq.\,\rref{scal}).
Thus, if $f$ is a maximizer for \rref{gag}, the same holds for each translated function $f^a$.
\fine
\end{rema}
The ratio in \rref{sec} has another, less trivial
invariance property described hereafter;
this depends crucially on the definition
\rref{due} of $r$.
\begin{prop}
\textbf{Proposition.}
For $p,q,j,n,\te,r$ as in \rref{min} \rref{due}, let us consider the scaling transformation
$f \mapsto f_{\lambda}$ (see Eq.\,\rref{scal}). For each $\lambda >0$, the following holds:
\beq f \in \ZZ{p}{q}{n} \Rightarrow f_{\lambda} \in \ZZ{p}{q}{n}~; \feq
\beq {\| \Di{j} f_{\lambda} \|_{\LL^r} \over
\| f_{\lambda}\|^{1 - \te}_{\LL^p} \| \Di{n} f_{\lambda} \|^{\te}_{\LL^q} }
= {\| \Di{j} f \|_{\LL^r} \over \| f \|^{1 - \te}_{\LL^p} \| \Di{n} f \|^\te_{\LL^q} } \quad
\mbox{for all $f \in \ZZ{p}{q}{n} \setminus \{0 \}$}~.
\label{scag} \feq
Therefore, if $f$ is a maximizer for the Gagliardo-Nirenberg inequality \rref{gag},
the same holds for each rescaled function $f_{\lambda}$.
\end{prop}
\textbf{Proof.} To derive \rref{scag} use the scaling relations \,\rref{nonum}
\rref{19}, together with Eq.\,\rref{due} for $r$.
\fine
\textbf{The case $\boma{\te=1}$
of the  Gagliardo-Nirenberg inequality.}
Let  $p,q \in [1,+\infty]$, $j, n \in [0,+\infty)$.
We assume the last condition in \rref{min} and
Eq.\,\rref{due} to be fulfilled with  $\te=1$, so that
\beq 0 \leqs {1 \over r} = {1 \over q} - {n - j \over d} \leqs 1~.\label{req1} \feq
If one considers Definition \ref{defgen} and applies it
mechanically with $\te=1$, one obtains the following statement:
\beq \ZZ{p}{q}{n} \subset \ZZ{p}{r}{j}~, \label{gag1}\feq
$$\| \Di{j} f \|_{\LL^r} \leqs G \| \Di{n} f \|_{\LL^q}
\quad\mbox{for some $G \in [0,+\infty)$ and all $f \in \ZZ{p}{q}{n}$}.$$
The related sharp constant is indicated with $G(p,q;j,n)$.
\par
It is natural is to consider an \textsl{extended}
inequality very similar to \rref{gag1} but making
no reference to $\LL^p$, namely:
\beq  \WW{q}{n} \subset \WW{r}{j}~, \label{gagn1} \feq
$$ \| \Di{j} f \|_{\LL^r} \leqs G \| \Di{n} f \|_{\LL^q}
\quad \mbox{for some $G \in [0,+\infty)$ and all $f \in \WW{q}{n}$}~. $$
The corresponding sharp constant is indicated with $G(q;j,n)$.
\begin{prop}
\label{lemga}
\textbf{Proposition.}
Let $q \neq 1, +\infty$. The inequality \rref{gag1} and the
extended inequality \rref{gagn1} are
equivalent; if they hold, $G(p,q;j,n) = G(q;j,n)$.
\end{prop}
\textbf{Proof.} An elementary density argument, reported for
completeness in \hbox{Appendix \ref{appega}}. \fine
\salto
\subsection{The Sobolev inequality}
\label{subsobo}
This inequality depends on a set of parameters $p$, $q$, $j$, $n$, $\te$
as in \rref{min} and on an additional parameter $t$, with
\beq 1 \leqs t \leqs +\infty~. \label{mint} \feq
\begin{prop}
\textbf{Definition.} Let $p,q,j,n,\te,r,t$ be as in \rref{min}  \rref{due} \rref{mint}.
The \textsl{Sobolev inequality} of order $(p,q;j,n,\te|t)$ is
the following statement:
\beq \ZZ{p}{q}{n} \subset \ZZ{p}{r}{j}~, \label{sob} \feq
$$ \| \Di{j} f \|_{\LL^r} \leqs S \, ( \| f \|^t_{\LL^p} + \| \Di{n} f \|^t_{\LL^q})^{1/t}
\quad \mbox{for some $S \in [0,+\infty)$ and all $f \in \ZZ{p}{q}{n}$}. $$
\parn
If this holds, the symbol $S(p,q;j,n,\te | t)$ indicates the sharp constant.
\end{prop}
\begin{rema}
\textbf{Remarks.} (i) Of course
\beq S(p,q;j,n,\te|t) = \sup_{f \in
\ZZ{p}{q}{n} \setminus \{0 \}} {\| \Di{j} f \|_{\LL^r} \over ( \| f
\|^t_{\LL^p} + \| \Di{n} f \|^t_{\LL^q})^{1/t}}~.  \label{fir} \feq
(ii) The ratio in the right hand side of \rref{fir} and, consequently, the set of
maximizers for \rref{sob} are invariant under translations $f \mapsto f^a$
($a \in \reali^d$, see Eq.\,\rref{scal}).
\fine
\end{rema}
\textbf{Scaling considerations.}
Let $p,q,j,n,\te, r,t$ be as in \rref{min} \rref{due} \rref{mint}. We consider the scaling transformation
$f \mapsto f_{\lambda}$, see again Eq.\,\rref{scal}, and its effect
on the ratio in \rref{fir}.
\begin{prop}
\textbf{Lemma.} For each $\lambda >0$, the following holds:
\beq {\| \Di{j} f_{\lambda} \|_{\LL^r} \over
( \| f_{\lambda}\|^t_{\LL^p} + \| \Di{n} f_{\lambda} \|^t_{\LL^q})^{1/t}} \label{scas} \feq
$$ =
{ \| \Di{j} f \|_{\LL^r} \over [ (\lambda^{ - \te (d/p -d/q + n)} \| f \|_{\LL^p})^t +
(\lambda^{ (1-\te) (d/p - d/q + n) } \| \Di{n} f \|_{\LL^q})^t ]^{1/t}}~~
\mbox{for all $f \in \ZZ{p}{q}{n} \setminus \{0 \}$} $$
(understanding both sides of \rref{scas} via \rref{nota2}, if $t=+\infty$).
\end{prop}
\textbf{Proof.} Use again Eqs.\,\rref{nonum} \rref{19} and Eq.\,\rref{due} for $r$.
\fine
\subsection{Connecting the Gagliardo-Nirenberg and Sobolev inequalities}
\label{subseccon}
Let $\te \in [0,1]$, $t \in [1,+\infty]$ and $a, b \in [0,+\infty)$; it is
well known that
\beq a^{1-\te} b^{\te} \leqs [\,(1 -\te)^{1 - \te} \te^\te\,]^{1/t} (a^t + b^t)^{1/t}~; \label{you} \feq
\beq a^{1-\te} b^{\te} = [\,(1 -\te)^{1 - \te} \te^\te\,]^{1/t} (a^t + b^t)^{1/t}~
\mbox{if $\te \neq 1$ for $t < + \infty$ and $b = \left(\te \over 1 - \te \right)^{1/t} \! a$}
\label{youeq} \feq
(recall the conventions \rref{nota1} \rref{nota2}; due
to \rref{nota1}, here and in the sequel $\left(\te \over 1 - \te \right)^0 = 1$
even for $\te=0,1$). Let us also mention the following variant of
\rref{youeq}, holding as well for $\te = 1$ and $t < + \infty$: if
$(a_{\lambda})$, $(b_{\lambda})$ are nets with values
in $(0,+\infty)$,
\beq \lim_{\lambda} {a_{\lambda}^{1-\te} b_{\lambda}^{\te} \over (a_{\lambda}^t + b_{\lambda}^t)^{1/t}}
= [\,(1 -\te)^{1 - \te} \te^\te\,]^{1/t} \quad
\mbox{if
$\lim_{\lambda} \dd{b_{\lambda} \over a_{\lambda}} = \left(\te \over 1 - \te \right)^{1/t}$}~.
~\label{youeqq} \feq
The  Gagliardo-Nirenberg and Sobolev inequalities can be
connected using scaling considerations, combined with the above elementary facts;
this was pointed out by several authors in special cases,
see, e.g., \cite{Agu} \cite{Lie4}. A general
formulation of these ideas is as follows.
\begin{prop}
\label{soga}
\textbf{Proposition.} Let $p,q,j,n,\te,r,t$ be as in \rref{min} \rref{due} \rref{mint};
consider the  Gagliardo-Nirenberg and Sobolev inequalities \rref{gag} and \rref{sob}.
Then, the following holds. \parn
(i) If the Gagliardo-Nirenberg inequality holds, the Sobolev inequality
holds as well and the corresponding sharp constants are related by
\beq {S(p,q;j,n,\te|t) \leqs [\,(1 -\te)^{1 - \te} \te^\te\,]^{1/t}}~G(p,q;j,n,\te)
~. \label{rela0} \feq
(ii) In addition to \rref{min} \rref{due} \rref{mint}, let
\beq {d \over p} - {d \over q} + n  \neq 0~; \label{let} \feq
then the  Gagliardo-Nirenberg and Sobolev inequalities
are equivalent. Whenever
these inequalities hold, their sharp constants are related by
\beq {S(p,q;j,n,\te|t) = [\,(1 -\te)^{1 - \te} \te^\te\,]^{1/t}}~G(p,q;j,n,\te)~.
\label{rela} \feq
(iii) With the condition \rref{let}, assume the  Gagliardo-Nirenberg
inequality to hold and possess a maximizer $f$; consider the rescaled functions
$f_{\lambda}$ ($\lambda >0$). The function $f_{\lambda}$ is a maximizer for
the Sobolev inequality if we put
\beq \lambda := \Big( \Big({\te \over 1 - \te} \Big)^{1/t} {\| f \|_{\LL^p} \over \| \Di{n} f \|_{\LL^q}}
\Big)^{1/(d/p - d/q + n)}.  \label{lamax} \feq
In the cases $\te=0$, $1$ and $t < + \infty$ (where \rref{lamax} would give formally
$\lambda = 0$ or $\lambda = + \infty$), the previous statement must be intended
in this limit sense: the fundamental ratio \rref{fir} for the Sobolev inequality
tends to the sharp constant if it is evaluated on $f_{\lambda}$
and the limit $\lambda \vain 0^{+}$ or $\lambda \vain + \infty$ is taken. \parn
(iv) With the condition \rref{let}, assume the Sobolev inequality
to hold and possess a maximizer $f$; then
$f$ is as well a Gagliardo-Nirenberg maximizer. \parn
\end{prop}
\textbf{Proof.} We proceed in several steps, using for the sharp
constants the shorthand notations
\beq G(p,q;j,n,\te) \equiv G~, \qquad
S(p,q;j,n,\te|t) \equiv S~. \feq
For the sake of brevity, we assume
\beq \te \neq 0,1~, \qquad t \neq + \infty~; \feq
the cases with $\te=0,1$ and/or $t = + \infty$ are treated
by simple variations of the considerations that follow
(taking possibly the limits $\lambda \vain 0^{+}$ or $\lambda \vain + \infty$
mentioned in item (iii)). \parn
\textsl{Step 1. If the  Gagliardo-Nirenberg inequality holds,
the Sobolev inequality holds
as well and its sharp constant fulfills}
\beq S \leqs G [\,(1 -\te)^{1 - \te} \te^\te\,]^{1/t} \label{thes2} \feq
\textsl{(so, statement (i) in the proposition is proved).}
Let $f \in \ZZ{p}{q}{n}$. The  inequality \rref{gag} with its sharp constant $G$
states that $f \in \ZZ{p}{r}{j}$, and
\beq \| \Di{j} f \|_{\LL^r} \leqs G \| f \|^{1 - \te}_{\LL^p} \| \Di{n} f \|^\te_{\LL^q}~; \label{310} \feq
on the other hand, Eq.\,\rref{you} with
$a= \| f \|_{\LL^p}$ and $b = \| \Di{n} f \|_{\LL^q}$ gives
\beq
\| f \|^{1 - \te}_{\LL^p} \| \Di{n} f \|^\te_{\LL^q} \leqs
[\,(1 -\te)^{1 - \te} \te^\te\,]^{1/t} \, ( \| f \|^t_{\LL^p} + \| \Di{n} f \|^t_{\LL^q})^{1/t}.
\label{you2} \feq
Due to \rref{310} \rref{you2}, the inequality \rref{sob} holds and its sharp constant $S$
fulfills \rref{thes2}. \parn
\textsl{Step 2. With the condition \rref{let}, assume the Sobolev inequality to hold;
then the  Gagliardo-Nirenberg inequality holds
as well, and its sharp constant fulfills}
\beq G \leqs {S \over [\,(1 -\te)^{1 - \te} \te^\te\,]^{1/t}}~. \label{thes1} \feq
To prove this, let us fix $f \in \ZZ{p}{q}{n}$; for each $\lambda > 0$,
the Sobolev inequality for $f_{\lambda}$ gives
\beq \| \Di{j} f_{\lambda} \|_{\LL^r} \leqs  S
( \| f_{\lambda}\|^t_{\LL^p} + \| \Di{n} f_{\lambda} \|^t_{\LL^q})^{1/t} \feq
which implies, due to \rref{scas},
\beq {~} \hspace{-0.4cm} \| \Di{j} f \|_{\LL^r} \leqs S \, (\lambda^{- \te (d/p - d/q + n) t} \| f \|^t_{\LL^p} +
\lambda^{(1-\te) (d/p - d/q + n) t} \| \Di{n} f \|^t_{\LL^q} )^{1/t}
\equiv  S \,F(\lambda). \feq
By the arbitrariness of $\lambda$, this gives
\beq \| \Di{j} f \|_{\LL^r} \leqs S \inf_{\lambda \in (0,+\infty)} F(\lambda), \label{35} \feq
and one obtains by elementary means that
\beq \inf_{\lambda \in (0,+\infty)} F(\lambda) =
{\| f \|^{1 - \te}_{\LL^p} \| \Di{n} f \|^\te_{\LL^q}
\over [\,(1 -\te)^{1 - \te} \te^\te\,]^{1/t} }
\label{36} \feq
(if $f \neq 0$, the $\inf$ of $F$ is attained for $\lambda$ as in \rref{lamax}).
Eqs.\,\rref{35} \rref{36} imply
\beq \| \Di{j} f \|_{\LL^r} \leqs {S \over [\,(1 -\te)^{1 - \te} \te^\te\,]^{1/t}}\,
\| f \|^{1 - \te}_{\LL^p} \| \Di{n} f \|^\te_{\LL^q} ~; \feq
this happens for all $f \in \ZZ{p}{q}{n}$, so the  Gagliardo-Nirenberg inequality
holds and its sharp constant $G$ fulfills \rref{thes1}. \parn
\textsl{Step 3. Proof of statement (ii) in the proposition.}
This follows immediately from Steps 1 and 2. \parn
\parn
\textsl{Step 4. With the condition \rref{let},
assume the  Gagliardo-Nirenberg inequality to hold
and admit a maximizer $f$; then the rescaled function $f_{\lambda}$ is a
Sobolev maximizer, if $\lambda$ is taken as in
\rref{lamax} (so, statement (iii) in the proposition is proved).}
According to our assumptions,
\beq f \in \ZZ{p}{q}{n} \setminus \{0 \}~, \qquad
\| \Di{j} f \|_{\LL^r} =  G \, \| f \|^{1 - \te}_{\LL^p} \| \Di{n} f \|^\te_{\LL^q}~. \feq
Since the ratio
$\| \Di{j} f \|_{\LL^r}/(\| f \|^{1 - \te}_{\LL^p} \| \Di{n} f \|^\te_{\LL^q})$
is scale invariant (recall \rref{scag}), we also have
\beq \| \Di{j} f_{\lambda} \|_{\LL^r} =
G \, \| f_{\lambda} \|^{1 - \te}_{\LL^p} \| \Di{n} f_{\lambda} \|^\te_{\LL^q} \label{eqq} \feq
for any $\lambda > 0$. On the other hand,
Eq.\,\rref{youeq} with $a = \| f_{\lambda} \|_{\LL^p}$ and
$b = \| \Di{n} f_{\lambda} \|_{\LL^q}$ gives
\beq \| f_{\lambda} \|^{1 - \te}_{\LL^p} \| \Di{n} f_{\lambda} \|^\te_{\LL^q} =
[\,(1 -\te)^{1 - \te} \te^\te\,]^{1/t} ( \| f_{\lambda} \|^t_{\LL^p} + \| \Di{n} f_{\lambda} \|^t_{\LL^q})^{1/t}
\label{gives2} \feq
under the condition
\beq \| \Di{n} f_{\lambda} \|_{\LL^q} =
\Big({\te \over 1 - \te}\Big)^{1/t} \,\| f_{\lambda} \|_{\LL^p} ~;\feq
due to the scaling properties \rref{nonum} \rref{19}, this condition is equivalent to
the equality
\beq
\la^{n - d/q} \| \Di{n} f \|_{\LL^q} =
\Big({\te \over 1 - \te}\Big)^{1/t} \,\lambda^{-d/p} \| f \|_{\LL^p}~,
\feq
which is fulfilled if $\lambda$ is chosen as in \rref{lamax}. With
this choice for $\lambda$, Eqs.\rref{eqq}\,\rref{gives2}
imply
\beq \| \Di{j} f_{\lambda} \|_{\LL^r} = G [\,(1 -\te)^{1 - \te} \te^\te\,]^{1/t}
( \| f_{\lambda} \|^t_{\LL^p} + \| \Di{n} f_{\lambda} \|^t_{\LL^q})^{1/t}; \feq
but $G [\,(1 -\te)^{1 - \te} \te^\te\,]^{1/t} = S$
due to \rref{rela}, so $f_{\lambda}$ is a Sobolev maximizer. \parn
\textsl{Step 5. With the condition \rref{let}, assume the Sobolev
inequality to hold and admit a maximizer $f$; then $f$ is a
Gagliardo-Nirenberg maximizer as well (so, statement (iv) in the proposition is proved).}
According to our assumptions,
\beq f \in \ZZ{p}{q}{n} \setminus \{0 \}~, \qquad
\| \Di{j} f \|_{\LL^r} =  S \, ( \| f \|^t_{\LL^p} + \| \Di{n} f \|^t_{\LL^q})^{1/t}; \feq
from here and from the inequality \rref{you2} we infer
\beq \| \Di{j} f \|_{\LL^r} \geqs  {S \over [\,(1 -\te)^{1 - \te} \te^\te\,]^{1/t}}
\, \| f \|^{1 - \te}_{\LL^p} \, \| \Di{n} f \|^\te_{\LL^q}~. \feq
But $S/[\,(1 -\te)^{1 - \te} \te^\te\,]^{1/t} = G$ due to Eq.\,\rref{rela}, so
\beq \| \Di{j} f \|_{\LL^r} \geqs  G \, \| f \|^{1 - \te}_{\LL^p} \, \| \Di{n} f \|^\te_{\LL^q}~. \feq
The reversed relation is just the Gagliardo-Nirenberg inequality, holding as well, so
\beq \| \Di{j} f \|_{\LL^r} =  G \, \| f \|^{1 - \te}_{\LL^p} \| \, \Di{n} f \|^\te_{\LL^q}~; \label{massf} \feq
this indicates that $f$ is a  Gagliardo-Nirenberg maximizer. \fine
\section{Validity conditions for the Gagliardo-Nirenberg
and Sobolev inequalities} \label{secfour} The original work of
Nirenberg \cite{Nir}
proves the inequality \rref{gag} for
\beq 1 \leqs p,q \leqs + \infty~, \quad 0 \leqs \te \leqs 1~, \quad
0 \leqs n,j < + \infty~, \label{uno} \feq
$$ 0 \leqs \te n - j \leqs d \left( {1 - \te \over p} + {\te \over q} \right) ~, \quad
\neg~ (1 < q < +\infty~,~n = j + {d \over q}~,~\te=1)~,  $$
with the additional assumption that $j, n$ be \textsl{integer}, and
defining $\| \Di{j} f \|_{\LL^r}$  (resp., $\| \Di{n} f \|_{\LL^q}$) in terms
of the $L^r$ (resp., $L^q$) norms of the partial derivatives
of order $j$ (resp., $n$) of $f$.
In the above $\neg$ indicates the logical negation
(the case negated in \rref{uno} is called ``exceptional'' in \cite{Nir}).
We note that the conditions \rref{uno} imply $1 \leqs r \leqs + \infty$, where
$r$ is defined by \rref{due}.
\par
To the best of our knowledge, the validity of
the  Gagliardo-Nirenberg inequality \rref{gag} with $j, n$ possibly noninteger
(and intending $\Di{j}, \Di{n}$, etc. as in the present paper) is nowadays known
under conditions slightly more restrictive than \rref{uno}, namely
\parn
\vbox{
\beq 1 < p,q < + \infty~,
\quad 0 \leqs \te \leqs 1~, \quad 0 \leqs n,j < + \infty~, \label{unofrac} \feq
$$ 0 \leqs \te n - j < d \left( {1 - \te \over p} + {\te \over q} \right)~; $$
}
see Corollary 1.5 of \cite{Haj0} for a proof based on the representation
of the homogeneous Sobolev spaces as special types of Triebel-Lizorkin spaces
\cite{Gra1} \cite{Gra2} \cite{Run}.
The conditions \rref{unofrac} imply $1 < r < + \infty$, with
$r$ as in \rref{due}.
\par
Let us recall that, according to Proposition \ref{soga}, for any
$t \in [1,+\infty]$ the Sobolev
inequality \rref{sob} is implied by \rref{gag}, and the two inequalities are equivalent
if $d/p - d/q + n \neq 0$. \par
In the next section we present \textsl{direct, autonomous} proofs for the validity of
the  inequalities \rref{gag} or \rref{sob} in all cases described by \rref{uno}
with $p=q=2$ (and $t=2$), for both integer
and noninteger values of $j, n$. Our analysis relies
on a collection of methods not requiring the Triebel-Lizorkin formalism of
\cite{Haj0}; these methods
give the sharp constants in some subcases, and accurate bounds
for them in the remaining ones. In
all the subcases analyzed, we exhibit direct proofs
and estimates on the related constants for
either \rref{gag} or \rref{sob},
according to convenience; next, we use Proposition \ref{soga}
to infer conclusions on the other inequality.
\section{Analysis of the $\boma{\LL^2}$ case}
\label{secelledue}
In this section we specialize the previous considerations to the case
\beq p=q=t=2~; \feq
so, our basic function space is $\LL^2$.
In the sequel, for each $n \in [0,+\infty)$, the spaces $\WW{2}{n}$
and $\ZZ{2}{2}{n}$ defined via Eqs.\,\rref{riesz}\rref{wpqn}
are indicated with $\Hh{n}$ and $\HH{n}$, respectively; thus
\beq \Hh{n}(\reali^d) \equiv \Hh{n} = \{ f \in \Phi'~|~ \Di{n} f \in \LL^2 \}~, \label{dewwn} \feq
\beq \HH{n}(\reali^d) \equiv \HH{n} = \{ f \in \LL^2~|~ \Di{n} f \in \LL^2 \}~. \label{dewn} \feq
Using the fact that the Fourier transform $\FF$ maps isometrically
$\LL^2$ into itself, we readily infer the following for $f \in \Phi'$:
\beq f \in \Hh{n} \Leftrightarrow | \k |^n \FF f \in \LL^2;\qquad f \in \Hh{n} \Rightarrow
\| \Di{n} f \|_{\LL^2} = \|\, |\k|^n \FF f \|_{\LL^2}~; \label{usec} \feq
\beq f \in \HH{n} \Leftrightarrow \FF f, \,|\k|^n \FF f \in \LL^2 \Leftrightarrow
\sqrt{1 + |\k|^{2 n}} \FF f \in \LL^2~;
\label{norm} \feq
$$ f \in \HH{n} \Rightarrow \sqrt{\| f \|^2_{\LL^2} + \| \Di{n} f \|^2_{\LL^2}} =
\| \sqrt{1 + |\k|^{2 n}} \FF f \|_{\LL^2}~. $$
One recognizes that $\HH{n}$ is the familiar $\LL^2$-based Sobolev space of Bessel potentials of order $n$
\cite{Aro} \cite{Maz}.
If $n$ is integer, we can describe $\Hh{n}, \HH{n}$ and $\| D^n f \|_{\LL^2}$ in terms
of the partial derivatives $\partial_j: \Phi' \vain \Phi'$ ($j=1,...,d$).
Indeed, by the Fourier representations of $\partial_j$ and $\Di{n}$
we have the following for $f \in \Phi'$ and $n \in \naturali$:
\parn
\vbox{
\beq f \in \Hh{n}
\Leftrightarrow \partial_{j_1...j_n} f \in \LL^2\quad\mbox{for $j_1,...,j_n
\in \{1,..,d\}$}~; \feq
\beq \mbox{if $f \in \Hh{n}$},~~
\| \Di{n} f \|^2_{\LL^2} = \sum_{j_1,...,j_n=1}^d \| \partial_{j_1...j_d} f \|_{\LL^2}^2~;
\label{concern} \feq
\beq f \in \HH{n} \Leftrightarrow f \in \LL^2 , \partial_{j_1...j_n} f \in \LL^2\quad\mbox{for $j_1,...,j_n
\in \{1,..,d\}$} \feq
$$ \Leftrightarrow \partial_{j_1...j_m} f \in \LL^2\quad\mbox{for $m \in \{0,...,n\}$ and
$j_1,...,j_m \in \{1,..,d\}$}  $$
}
({\footnote{To check some of the above statements, note the following:
\begin{itemize}
\item[(i)] $|\k|^{2 n} |\FF f|^2$ $= (\sum_{j=1}^ d \k^2_j)^n | \FF f|^2$
$= \sum_{j_1,..,j_n=1}^ d | \k_{j_1} ... \k_{j_n} \FF f |^2$;
\item[(ii)] if $m \in \{0,...,n\}$ and $j_1,...,j_m \in \{1,...,d\}$,
$| \k_{j_1} ... \k_{j_m} \FF f |^2 \leqs |\k|^{2 m} | \FF f |^2
\leqs (1 + |\k|^{2 m}) | \FF f |^2 \leqs C_{n m} (1 + |\k|^{2 n}) | \FF f |^2 $
for a suitable, positive constant $C_{n m}$.
\end{itemize} }}).
The parameters of the $\LL^2$ case are three real numbers $j, n, \te$;
the conditions \rref{uno} and the definition
\rref{due} for $r \equiv r(j,n,\te)$
take the form
\parn
\vbox{
\beq
0 \leqs \te \leqs 1\,, \quad 0 \leqs n,j < + \infty\,, \quad 0 \leqs \te n - j \leqs {d \over 2}\,,
\quad \te \neq 1\, \mbox{if $n = j + \dd{d \over 2}$}\,;
\label{unoo} \feq
$$ {1 \over r} = {1 \over 2} - {\te n - j\over d} ~.
$$}
We note that the above conditions imply $0 \leqs 1/r \leqs 1/2$, i.e.,
\beq 2 \leqs r \leqs + \infty.  \label{remr} \feq
The Gagliardo-Nirenberg and Sobolev inequalities \rref{gag} \rref{sob} read
\beq \HH{n} \subset \ZZ{2}{r}{j},\quad \| \Di{j} f \|_{\LL^r} \leqs
G \| f \|^{1 - \te}_{\LL^2} \| \Di{n} f \|^\te_{\LL^2}
\label{gagg} \feq
$$ \mbox{for some $G \in [0,+\infty)$ and all $f \in \HH{n}$}\,; $$
\beq \HH{n} \subset \ZZ{2}{r}{j},\quad\| \Di{j} f \|_{\LL^r} \leqs
S \sqrt{\| f \|^2_{\LL^2} + \| \Di{n} f \|^2_{\LL^2}}
\label{sobb} \feq
$$ \quad\mbox{for some $S \in [0,+\infty)$ and all $f \in \HH{n}$}\,. $$
From now on, the sharp constants of these inequalities are denoted
with $G(j,n,\te)$ and $S(j,n,\te)$,
respectively.
\par
As anticipated, in this section we present direct proofs
of \rref{gagg} and/or \rref{sobb}, fitting the $\LL^2$ framework,
for all $j,n,\te$ as in \rref{unoo}. Our analysis,
carried over in the following subsections, follows the
scheme (i)-(vii) already outlined in the Introduction.
\subsection{The elementary case $\boma{j = \theta n}$}
\label{jten}
If we put $j=\te n$ in Eq.\,\rref{unoo} we obtain the conditions
\beq 0 \leqs \te \leqs 1\,, \quad 0 \leqs n < + \infty\,, \feq
and the definition of $r$ written therein gives $r=2$.
The Gagliardo-Nirenberg inequality \rref{gagg}
for this case reads
\parn
\vbox{
\beq \HH{n} \subset \HH{\te n}, \qquad \| \Di{\te n} f \|_{\LL^2} \leqs
\Ge \| f \|^{1 - \te}_{\LL^2} \| \Di{n} f \|^{\te}_{\LL^2} \label{gagbi} \feq
$$ ~\mbox{for some $\Ge \in [0,+\infty)$ and all $f \in \HH{n}$}. $$
}
As a matter of fact, this inequality is obtained by manipulations based on
the H\"older inequality.
\begin{prop}
\label{propgagbi}
\textbf{Proposition.} For $0 \leqs \te \leqs 1$ and $0 \leqs n < +\infty$,
Eq.\,\rref{gagbi} holds with sharp constant $\Ge(n,\te)=1$.
\end{prop}
\textbf{Proof.} It is divided in two steps, whose combination yields the thesis.
\parn
\textsl{Step 1. Eq.\,\rref{gagbi} holds with sharp constant $\Ge(n,\te) \leqs 1$.}
We use the H\"older inequality
$\int u v \leqs (\int u^{p})^{1/p} (\int v^q)^{1/q}$
(for $u, v \geqs 0$ and $p, q \geqs 1$, $1/p + 1/q = 1$) and apply it
with $u = |\FF f|^{\,2(1-\te)}$, $v = |\k|^{\,2 \te n} |\FF f|^{\,2 \te}$,
$p=1/(1 - \te)$, $q = 1/\te$; the result
is the inequality $\| | \k|^{\te n} \FF f \|_{\LL^2} \leqs \| \FF f \|^{1 - \te}_{\LL^2}
\| |\k|^n \FF f \|^{\te}_{\LL^2}$ which is equivalent to
$$ \HH{n} \subset \HH{\te n}, \qquad  \| \Di{\te n} f \|_{\LL^2} \leqs
\| f \|^{1 - \te}_{\LL^2} \| \Di{n} f \|^{\te}_{\LL^2} \quad
\mbox{for all $f \in \HH{n}$}. $$
This gives the statement to be proved. \parn
\textsl{Step 2. The sharp constant of \rref{gagbi} is such that $\Ge(n,\te) \geqs 1$}.
Of course, for each $f \in \HH{n} \setminus \{0 \}$ we have
\beq \Ge(n,\te) \geqs { \| \Di{\te n} f \|_{\LL^2} \over \| f \|^{1 - \te}_{\LL^2} \| \Di{n} f \|^{\te}_{\LL^2}}~.
\label{comesic} \feq
To go on, for $a \in \reali^d$ and $\ep \in (0,+\infty)$ let us put
\beq \delta_{a \ep} : \reali^d \vain [0,+\infty)~, \qquad \delta_{a \ep}(k) :=
{\chi_{B_{a \ep}}(k) \over v_{\ep}}~, \feq
where $B_{a \ep}$ is the ball in $\reali^d$ of center $a$ and radius $\ep$,
$\chi_{B_{a \ep}}$ is the characteristic function of this ball and $v_{\ep} =
\pi^{d/2} \ep^ d/\Gamma(d/2+1)$ is the volume of the ball, so that
$\int_{\reali^d} dk \, \delta_{a \ep}(k) = 1$.
For fixed $a$ this function approaches the Dirac delta at $a$ in the
limit $\ep \vain 0$; more precisely, for each continuous function $g : \reali^d \vain \complessi$
we have
\beq \int_{\reali^d} d k \, g(k) \, \delta_{a \ep} \vain g(a) \qquad \mbox{for $\ep \vain 0$}~. \feq
For $a$, $\ep$ as above, let us introduce
the function
\beq f_{a \ep} := \FF^{-1} \sqrt{\delta_{a \ep}}~. \feq
Then $f_{a \ep} \in \LL^2$ and
\beq \| f_{a \ep} \|_{\LL^2} = \| \FF f_{a \ep} \|_{\LL^2} = \sqrt{\int_{\reali^d} d k \, \delta_{a \ep}(k) } = 1~; \feq
moreover, for each $m \in [0,+\infty)$ one has $\Di{m} f_{a \ep} \in \LL^2$ and
\beq \| \Di{m} f_{a \ep} \|_{\LL^2} = \| |\k|^m \FF f_{a \ep} \|_{\LL^2} =
\sqrt{\int_{\reali^d} d k \, |k|^{2 m} \delta_{a \ep}(k) } \vain |a|^m \quad \mbox{for $\ep \vain 0$}~. \feq
To conclude, let us fix $a \in \reali^d \setminus \{0 \}$ and
apply \rref{comesic} with $f = f_{a \ep}$; in this way we get
\beq \Ge(n,\te) \geqs
{\| \Di{\te n} f_{a \ep} \|_{\LL^2} \over \| f_{a \ep} \|^{1 - \te}_{\LL^2} \| \Di{n} f_{a \ep} \|^{\te}_{\LL^2}}
\vain { |a|^{\te n} \over (|a|^n)^{\te} } = 1 \quad \mbox{for $\ep \vain 0$}~. \feq
\fine
Let us pass to the Sobolev inequality \rref{sobb}, that in the present
case $j=\te n$ reads
\beq \HH{n} \subset \HH{\te n}, \qquad \| \Di{\te n} f \|_{\LL^2} \leqs
\Se \sqrt{\| f \|^2_{\LL^2} + \| \Di{n} f \|^2_{\LL^2}}~~\mbox{for all $f \in \HH{n}$}~.
\label{sobbi} \feq
From the previous result on the Gagliardo-Nirenberg inequality we obtain the following result.
\begin{prop}
\label{corbyeq}
\textbf{Corollary}. Let $0 \leqs \te \leqs 1$, $0 < n < +\infty$.
The inequality \rref{sobbi} holds with sharp constant
\beq \Se(n, \te) = \sqrt{(1-\te)^{1 - \te} \te^{\te}}~. \label{byeq} \feq
\end{prop}
\textbf{Proof.} Use Proposition \ref{propgagbi}, together
with Proposition \ref{soga} on the general relations
between the Gagliardo-Nirenberg and Sobolev inequalities
(especially, Eq.\,\rref{rela}). \fine
\begin{rema}
\label{remobvious}
\textbf{Remark.}
Of course, Eq.\,\rref{sobbi} holds as well for $0 \leqs \te \leqs 1$,
$n =0$ with sharp constant $\Se(0,\te) = 1/\sqrt{2}$ (in fact,
the inequality in \rref{sobbi} with $n=0$ holds as an equality for $\Se = 1/\sqrt{2}$
and any $f \in \LL^2$).
\fine
\end{rema}
\subsection{The case $\boma{\te=1}$}
\label{subsuno}
If we put $\te=1$ in the general conditions \rref{unoo} we obtain
\beq 0 \leqs j \leqs n < j + {d \over 2}~, \qquad
{1 \over r} = {1\over 2} - {n - j \over d}~; \label{assu1} \feq
note that $r \in [2,+ \infty)$.
The inequalities \rref{gagg} \rref{sobb} read
\parn
\vbox{
\beq \HH{n} \subset \ZZ{2}{r}{j}~, \quad \| \Di{j} f \|_{\LL^r} \leqs \Ga \| \Di{n} f \|_{\LL^2}
\label{gaggl} \feq
$$ \mbox{for some $\Ga \in [0,+\infty)$ and all $f \in \HH{n}$}; $$
}
\vbox{
\beq \HH{n} \subset \ZZ{2}{r}{j}~, \quad \| \Di{j} f \|_{\LL^r} \leqs
\Sa \sqrt{ \| f \|^2_{\LL^2} + \| \Di{n} f \|^2_{\LL^2}} \label{sobbl} \feq
$$ \mbox{for some $\Sa \in [0,+\infty)$ and all $f \in \HH{n}$}. $$
}
In this case, it is natural to consider as well the extended
Gagliardo-Nirenberg inequality, i.e., statement \rref{gagn1} with $p=q=2$; this
reads
\beq \Hh{n} \subset \WW{r}{j},~\| \Di{j} f \|_{\LL^r} \leqs \Ga \| \Di{n} f
\|_{\LL^2} ~\mbox{for some $\Ga \in [0,+\infty)$ and all $f \in \Hh{n}$},
\label{gagn21} \feq
and is equivalent to \rref{gaggl} due to
Proposition \ref{lemga}. The inequality
\rref{gagn1} and, in particular, its $\LL^2$ case
\rref{gagn21} are strictly connected
with the Hardy-Littlewood-Sobolev inequality \cite{HL} \cite{Sob}
concerning
convolution with a power of the radius $|\x|$; this connection has
a crucial role even in Sobolev's seminal paper \cite{Sob} and is presented, e.g., by
Mizohata \cite{Mizo} in a more up-to-date language.
The sharp constants and maximizers of the Hardy-Littlewood-Sobolev
inequality have been determined more recently
by Lieb \cite{Lie0} for some cases, including the
$\LL^2$ case; by the previously mentioned
equivalence, these results of Lieb can be used
to determine the sharp constant and maximizers
for \rref{gagn21}, a fact somehow suggested by \cite{Lie0}
and described more explicitly in \cite{Cot}. \par
The situation outlined above can be understood
starting from the subcase $j=0$ of \rref{gagn21};
this is treated in the following proposition
(and in the subsequent
Remarks \ref{rem55}),
very close to Theorem 1.1 of \cite{Cot}
({\footnote{See also the announcement of this
theorem in \cite{Cot0}.
The cited theorem of \cite{Cot0} \cite{Cot} contains some imprecision, since it does not
mention $\Hh{n}$ and always refers to $\HH{n}$; in particular, it seems to
indicate that the maximizer $f$ in Eq.\,\rref{fmax} is in $H^n$
for all $n$ as in \rref{letn}.}}).
\begin{prop}
\label{prolieb}
\textbf{Proposition.} Let
\beq 0 \leqs n < {d \over 2}~, \quad {1 \over \rs} = {1 \over 2} - {n \over d} \label{letn} \feq
(where the second equation is understood as the definition of
$\rs \in [2,  +\infty)$).
Then
\beq \Hh{n} \subset \LL^{\rs}~, \quad \| f \|_{\LL^{\rs}} \leqs \Ga \| \Di{n} f \|_{\LL^2}
\quad \mbox{for some $\Ga \in [0,+\infty)$ and all $f \in \Hh{n}$}. \label{gagg21} \feq
Moreover, the sharp constant in \rref{gagg21} is
\beq \Ga(n) = {1 \over (4 \pi)^{n/2}} \sqrt{\Gamma(d/2-n) \over \Gamma(d/2+n)}
\left( \Gamma(d) \over \Gamma(d/2) \right)^{n/d}\!. \label{sharplz} \feq
A maximizer for \rref{gagg21} is
\beq f := {1 \over (1 + |\x|^2)^{d/2 - n}} =
{1 \over 2^{d/2-n-1} \Gamma(d/2-n)} \, \FF^{-1} \left({K_{n}(|\k|) \over |\k|^n} \right)
\in \Hh{n}, \label{fmax} \feq
where $K_n$ denotes the modified Bessel function of the second kind (Macdonald function);
note that $K_{n}(|\k|)/|\k|^n \in \LL^1$.
The above function $f$ is in $\LL^2$ (and thus in $\HH{n}$) if and only if
$n < d/4$.
\end{prop}
\textbf{Proof.}
It is divided in some steps; the main point is Step 1,
reproducing a basic result of Lieb on the Hardy-Littlewood-Sobolev inequality. In the
sequel we use the convolution $*$ and some of its properties, reviewed
in Section \ref{prelim}. \parn
\textsl{Step 1 (a sharp Hardy-Littlewood-Sobolev inequality). Let
$n$, $\rs$ be as in \rref{letn} and, in addition, $n \neq 0$.
Then
\parn
\vbox{
\beq h \in \LL^2 \quad \Rightarrow \quad {1 \over |\x|^{d-n}} * h \in \LL^{\rs}~,~
\| {1 \over |\x|^{d-n}} * h \|_{\LL^{\rs}} \leqs \N_{n} \| h \|_{\LL^2}~, \label{hali} \feq
$$ \N_{n} \equiv \N_{n d} := \pi^{d/2 - n/2}
{\Gamma(n/2) \over \Gamma(d/2 - n/2)}
\sqrt{ {\Gamma(d/2-n) \over \Gamma(d/2 + n)}}
\left( \Gamma(d) \over \Gamma(d/2) \right)^{n/d}.$$
}
The inequality in \rref{hali} is fulfilled as an equality by the function
\beq h := {1 \over |\x|^{d-n}}*{1 \over (1 + |\x|^2)^{d/2 +n}} \in \LL^2~.
\label{mahali} \feq
}
For all these statements see \cite{Lie0}, Corollary 3.2, item (ii).
\parn
\textsl{Step 2 (essentially, a reformulation of Step 1 via the fractional
Laplacian). Let $n$, $\rs$ and $\Ga(n)$ be as in \rref{letn}
\rref{sharplz}. Then
\beq h \in \LL^2 \quad \Rightarrow \quad \Di{-n} h \in
\LL^{\rs}~,~ \| \Di{-n} h \|_{\LL^{\rs}} \leqs \Ga(n) \| h \|_{\LL^2}~,
\label{halir} \feq
with $\Ga(n)$ as in Eq.\,\rref{sharplz}. The inequality in
\rref{halir} if fulfilled as an equality by the function
\beq h := \Di{-n} {1
\over (1 + |\x|^2)^{d/2 + n}} \in \LL^2~. \label{mahalir} \feq
}
We first prove these statements for $n
\neq 0$ (so that $0 < n < {d/2}$). To get the thesis, it suffices to write down
the results of Step 1 and note that $|\x|^{-(d-n)} *... =
\AZ^{-1}_{n} \Di{-n}(...)$ with $\AZ_n$ as in Eq.\,\rref{223}. In particular,
the inequality $\| \, |\x|^{-(d-n)} * h \|_{\LL^{\rs}} \leqs \N_{n} \| h
\|_{\LL^2}$ of Step 1 becomes $\| \Di{-n} h \|_{\LL^{\rs}} \leqs \AZ_n \N_{n}
\| h \|_{\LL^2}$, and one readily checks that $\AZ_n \N_{n} = \Ga(n)$. \par
Let us pass to the case $n=0$. Then all statements to be proved hold trivially
since $\rs =2$ and $\Ga(0)=1$; of course, the inequality \rref{halir} is
fulfilled as an equality by any function in $\LL^2$, including the function in
\rref{halir} with $n=0$.
\parn
\textsl{Step 3. Proof of all statements in the proposition.} Let
again $n$, $\rs$ and $\Ga(n)$ be as in \rref{letn} and
\rref{sharplz}. By the very definition of $\Hh{n}$, the map $f \mapsto
\Di{n} f$ is one-to-one between $\Hh{n}$ and $\LL^2$. Therefore, applying the
results of Step 2 with $h = \Di{n} f$ ($f \in \Hh{n}$) we infer that
\beq f \in \Hh{n} \quad \Rightarrow \quad f \in \LL^{\rs}~,~ \| f \|_{\LL^{\rs}} \leqs
\Ga(n) \| \Di{n} f \|_{\LL^2} \label{halirr} \feq
and that the above
inequality is fulfilled as an equality by the function
\beq f = \Di{-2 n} {1
\over (1 + |\x|^2)^{d/2 + n}} \in \Hh{n} \label{mahalirr} \feq
(or by any multiple of it by a constant factor). Summing up, the extended
Gagliardo-Nirenberg inequality holds in the case under consideration with
$\Ga(n)$ as sharp constant and the function \rref{mahalirr} as a maximizer.
\par
Hereafter we show, via the related Fourier
representations, that the maximizer \rref{mahalirr} coincides, up to a constant
factor, with the function in Eq.\,\rref{fmax}. Indeed, by definition
\hbox{$\Di{-2 n} (1 + |\x|^2)^{-(d/2 + n)} = \FF^{-1}\big(|\k|^{-2 n} \FF (1 +
|\x|^2)^{-(d/2 + n)} \big)$}; moreover, \hbox{$\FF (1 + |\x|^2)^{-(d/2 + n)} =
2^{-(d/2 + n -1)} \Gamma(d/2+n)^{-1} |\k|^{n} K_{-n}(|\k|)$} by Lemma
\ref{lembes} in Appendix \ref{appe} (with $\mu=\sigma=-n$) and $K_{-n}=K_n$, as well known. Thus
\beq
\Di{-2 n} {1 \over (1 + |\x|^2)^{d/2 + n}} = {1 \over 2^{d/2 + n - 1}
\Gamma(d/2 + n)} \, \FF^{-1} \left({K_{n}(|\k|) \over |\k|^n} \right)\,; \feq
on the other hand, using again Lemma \ref{lembes} (now with $\mu=\sigma=n$) we
find that $K_{n}(|\k|)/|\k|^n \in \LL^1$, and
\beq {1 \over 2^{d/2-n-1} \Gamma(d/2-n)} \, \FF^{-1}
\left({K_{n}(|\k|) \over |\k|^n} \right) = {1 \over (1 + |\x|^2)^{d/2 - n}}~.
\feq
To conclude the proof, it remains to show that the maximizer $(1 +
|\x|^2)^{-(d/2 - n)}$ is in $\LL^2$ (and thus in $\HH{n})$ if and only if $n <
d/4$; the verification is trivial.
\fine
\begin{rema}
\label{rem55}
\textbf{Remarks.} (i) For $n=1$, the result of the above proposition was
obtained by Aubin \cite{Aub} and Talenti \cite{Tal} some years before \cite{Lie0}. \parn
(ii) The analysis of Lieb on the Hardy-Littlewood-Sobolev
inequality \rref{hali} shows as well that the function
$h$ in \rref{mahali} is the \textsl{unique} maximizer
up to translation, rescaling and multiplication by
a constant factor. Therefore, one can make a similar statement for the maximizer \rref{fmax}
of \rref{gagg21}; in the sequel we do not insist on
such uniqueness issues. \fine
\end{rema}
Proposition \ref{prolieb} has a straightforward generalization to the case \rref{assu1}.
\begin{prop}
\label{coronj}
\textbf{Corollary.} Let $j, n, r$ be as in \rref{assu1}, and
consider the extended Gagliardo-Nirenberg inequality \rref{gagn21}:
$$ \Hh{n} \subset \WW{r}{j}, \quad \| \Di{j} f \|_{\LL^r} \leqs \Ga \| \Di{n} f \|_{\LL^2}
\quad \mbox{for some $\Ga \in [0,+\infty)$ and all $f \in \Hh{n}$}. $$
This statement is true, and the sharp constant therein is
\beq \Ga(j,n) = {1 \over (4 \pi)^{(n - j)/2}} \sqrt{\Gamma(d/2-n + j) \over \Gamma(d/2+n -j)}
\left( \Gamma(d) \over \Gamma(d/2) \right)^{(n-j)/d}. \label{sharpo} \feq
A maximizer for \rref{gagn21} is
\beq f := \Di{-j} {1 \over (1 + |\x|^2)^{d/2 - n + j}} \in \Hh{n}~; \label{djef} \feq
this function can be written as
\beq f = {1 \over 2^{d/2 - n + j -1}\Gamma(d/2- n+ j)} \, \FF^{-1} \Big({K_{n-j}(|\k|) \over |\k|^{n}}\Big)
\label{prefhyp} \feq
(note that ${K_{n-j}(|\k|)/|\k|^{n}}$ makes sense in $\Psi'$ as the
product between $|\k|^{-j}$ and the $\LL^1$ function ${K_{n-j}(|\k|)/|\k|^{n-j}}$).\par
One has ${K_{n-j}(|\k|)/|\k|^{n}} \in \LL^1$ if and only if the stronger condition
$n < j/2 + d/2$ holds; in this case
$f \in \LL^\infty$, and we have the representation
\beq f = {\Gamma(d/2 - j/2) \Gamma(d/2 - n + j/2) \over 2^{j} \, \Gamma(d/2) \Gamma(d/2 - n + j)}\,\,
{}_2 F_1( d/2 - j/2, d/2 - n + j/2; d/2; - |\x|^2)~. \label{fhyp} \feq
$f$ is in $\LL^2$ (and thus in $\HH{n}$) if and only if
the even stronger condition $n < j/2 + d/4$ holds.
\end{prop}
\textbf{Proof.} Let $0 \leqs s < d/2$, and write the inequality \rref{gagg21}
with $n$ replaced by $s$ and with the sharp constant
therein. This states that, with $1/\rs = 1/2 - s/d$, $\Hh{s} \subset \LL^{\rs}$
and  $\| g \|_{\LL^{\rs}} \leqs \Ga(s) \| \Di{s} g \|_{\LL^2}$
for all $g \in \Hh{s}$; $\Ga(s)$ is as in \rref{sharplz}
with $n$ replaced by $s$, and the inequality
holds as an equality if $g =1/(1 + |\x|^2)^{d/2 - s}$. \par
Now let $j, n, r$ be as in \rref{assu1}, and
write the inequality $\| g \|_{\LL^{\rs}} \leqs \Ga(s) \| \Di{s} g \|_{\LL^2}$ with $s := n - j$ and
$g := \Di{j} f$, $f \in \Hh{n}$; this gives the inequality \rref{gaggl}
with $\Ga(j,n) = \Ga(n-j)$, which has the explicit
expression \rref{sharpo}. Due to Proposition \ref{prolieb}, \rref{gaggl}
becomes an equality if we consider the element $f \in \Hh{n}$ defined
by Eq.\,\rref{djef}, which is equivalent to
\beq \FF^{-1} (|\k|^j \FF f) = {1 \over (1 + |\x|^2)^{d/2 - n + j}}~; \label{fulf} \feq
but
\beq {1 \over (1 +|\x|^2)^{d/2 - n + j}} = {1 \over 2^{d/2 - n + j -1}\, \Gamma(d/2-n+j)}\,
\FF^{-1}\Big({K_{n-j}(|\k|) \over |\k|^{n-j}}\Big) \label{firstff} \feq
and ${K_{n-j}(|\k|)/ |\k|^{n-j}} \in \LL^1$; to prove these
statements, use Lemma \ref{lembes} in Appendix \ref{appe} with $\mu = \sigma = n-j$. In view
of \rref{firstff}, Eq.\,\rref{fulf} is equivalent to Eq.\,\rref{prefhyp}
$$ f  = {1 \over 2^{d/2  - n + j -1}\Gamma(d/2-n+j)} \,\FF^{-1}
\Big({K_{n-j}(|\k|) \over |\k|^{n}}\Big)~. $$
Now, using Lemma \ref{lembes} with $\mu=n-j$ and $\sigma=n$
we obtain the remaining statements to be proved: $K_{n-j}(|\k|)/|\k|^{n} \in \LL^1$
if and only if $n < j/2 + d/2$, in this case  $f \in \LL^\infty$
and Eq.\,\rref{fhyp} holds; $f \in \LL^2$ if and only if $n < j/2 + d/4$.
\fine
\begin{rema}
\textbf{Remark.} Obviously enough, one would like to generalize
Eq.\,\rref{fhyp} to all $j,n$ as in \rref{assu1}, removing the
limitation $n < j/2 + d/2$. To illustrate the related difficulties,
it suffices to consider the right hand side of \rref{fhyp}
for fixed $j$ and $n \vain (j/2 + d/2)^{-}$. In this limit,
the term $\Gamma(d/2 - n + j/2)$ diverges and the
hypergeometric function in \rref{fhyp} becomes
formally ${}_2 F_1( d/2 - j/2, 0; d/2; - |\x|^2) = 1$.
However, a constant function represents the zero
element of $\Phi'$ ({\footnote{Let us recall that $\Phi'$ can be identified with
$\SS'$ modulo the polynomial functions: see the footnote before Eq.\,\rref{lizo}.}}); so,
the right hand side of Eq.\,\rref{fhyp} gives, in the
limit $n \vain (j/2 + d/2)^{-}$, an indeterminate form
in $\Phi'$. We leave to future work
the analysis of this problem and
the discussion of \rref{fhyp} for $n > j/2 + d/2$
(perhaps possible by analytic continuation arguments).
\fine
\end{rema}
Here is another consequence of the previous results.
\begin{prop}
\label{corro}
\textbf{Corollary.} Let $j, n, r$ be as in \rref{assu1}, and
consider the inequalities \rref{gagg} \rref{sobb} :
$$ \HH{n} \subset \ZZ{2}{r}{j}, ~\| \Di{j} f \|_{\LL^r} \leqs \Ga \| \Di{n} f \|_{\LL^2}
\quad \mbox{for some $\Ga \in [0,+\infty)$ and all $f \in \HH{n}$}; $$
$$ \HH{n} \subset \ZZ{2}{r}{j}, ~ \| \Di{j} f \|_{\LL^r}
\leqs \Sa \sqrt{ \| f \|^2_{\LL^2} + \| \Di{n} f \|^2_{\LL^2}}
~\mbox{for some $\Sa \in [0,+\infty)$ and all $f \in \HH{n}$}. $$
These are true. Their sharp constants $\Ga(j,n)$, $\Sa(j,n)$ are both
equal to the right hand side of Eq.\,\rref{sharpo}, if $n \neq 0$; in the subcase $n=0$,
implying $j=0$, the sharp constants are $\Ga(0,0) = 1$ and $\Sa(0,0) = 1/\sqrt{2}$.
\end{prop}
\textbf{Proof.} Due to Proposition \ref{lemga}, the inequality
\rref{gagg} is equivalent to the extended inequality
\rref{gagn21} of Corollary \ref{coronj}; moreover they have
the same sharp constant $\Ga(j,n)$, given by \rref{sharpo}.
On the other hand Proposition \ref{soga} with $p=q=t=2$ and $\te=1$ ensures that,
for $n \neq 0$,
the inequalities \rref{gagg} \rref{sobb} are equivalent
and possess the same sharp constant. The statements on the subcase
$n=0$ (and $j=0$) are obvious (and make reference to already
mentioned facts, see Remark \ref{remobvious}).
\fine
\subsection{An ``almost general'' $\boma{\LL^2}$ case: proof of the
Gagliardo-Nirenberg and Sobolev inequalities and
upper bounds for their sharp constants}
\label{subsdue}
In this subsection we make the assumptions
\beq 0 \leqs \te \leqs 1\,, \quad 0 \leqs n, j < + \infty\,, \quad
0 \leqs \te n - j < {d \over 2}\,; \quad {1 \over r} =  {1 \over 2} - {\te n - j\over d} ~,\label{unoox} \feq
which differ from the general $\LL^2$ conditions
\rref{unoo} since they exclude the case with
$\te n - j = \dd{d/2}$ (occurring if only if
$n \neq 0$ and $\te = j/n + d/2 n$);
in this sense, we are considering an almost general
$\LL^2$ case. We note that $r \in [2,+ \infty)$.
\par
We consider the inequalities \rref{gagg} \rref{sobb} and their
sharp constants $G(j,n,\te)$, $S(j,n,\te)$.
Following the approach of
\cite{Haj} (see Corollary 2.3), we show
how to infer \rref{gagg} and
an upper bound on $G(j,n,\te)$ using
results on the cases $j = \te n$ and $\te=1$
(see our subsections \ref{jten} and \ref{subsuno})
({\footnote{To be precise, Corollary 2.3 of
\cite{Haj} is about the Gagliardo-Nirenberg
inequality \rref{gag} with $p,q$
arbitrary and $j=0$, whereas here
$p=q=2$ and $j$ can be nonzero.}});
this has implications on the Sobolev
inequality, according to Proposition \ref{soga}.
\begin{prop}
\label{questa}
\textbf{Proposition.}
Let $j,n,\te,r$ be as in \rref{unoox};
consider the Gagliardo-Nirenberg inequality \rref{gagg}
$$ \HH{n} \subset \ZZ{2}{r}{j}~, \quad \| \Di{j} f \|_{\LL^r} \leqs G \| f \|^{1 - \te}_{\LL^2} \| \Di{n} f \|^\te_{\LL^2}
\quad \mbox{for all $f \in \HH{n}$}. $$
This statement is true, and the related sharp constant $G(j,n,\te)$ has the upper bound
\beq G(j,n,\te) \leqs \Gp(j,n,\te)~, \label{gpp} \feq
$$ \Gp(j,n,\te)
:= {1 \over (4 \pi)^{(\te n - j)/2}} \sqrt{{\Gamma(d/2- \te n + j) \over \Gamma(d/2+ \te n -j)}}\,
\left( {\Gamma(d) \over \Gamma(d/2)} \right)^{(\te n - j)/d}. $$
The equality $G(j,n,\te) = \Gp(j,n,\te)$ holds in the
following subcases: \parn
(a) $j = \te n$, where $G(j,n,\te)=1$ due to Proposition \ref{propgagbi}; \parn
(b) $\te =1$, where $G(j,n,\te) = \Ga(j,n)$ as in \rref{sharpo}. \parn
\end{prop}
\textbf{Proof.} The validity of Eq.\,\rref{gagg} and the sharp constants
in the subcases (a)(b) are known from subsections \ref{jten}, \ref{subsuno};
the equality $G(j,n,\te) = \Gp(j,n,\te)$ in these
subcases is checked immediately comparing the known
values of the sharp constants with Eq.\,\rref{gpp}.
\par
In the rest of the proof we assume $\te \neq 0,1$.
With this assumption (and with the previous ones
in \rref{unoox} for $j,n,\te,r$), we see
that Eq.\,\rref{assu1} holds with $n$ replaced by $\te n$
and $j,r$ as before; therefore Corollary \ref{corro} gives
\beq \HH{\te n} \subset \ZZ{2}{r}{j}~, \qquad
\| \Di{j} f \|_{\LL^r} \leqs \Ga(j,\te n) \| \Di{\te n} f \|_{\LL^2}
\quad \mbox{for all $f \in \HH{\te n}$,} \label{cogive} \feq
with $\Ga(j,\te n)$ defined following Eq.\,\rref{sharpo}; one checks that
\beq \Ga(j,\te n) = \Gp(j,n,\te) \quad \mbox{as in \rref{gpp}}~. \label{finds} \feq
On the other hand, due to Proposition \ref{propgagbi} we have
\beq \HH{n} \subset \HH{\te n}, \qquad  \| \Di{\te n} f \|_{\LL^2} \leqs
\| f \|^{1 - \te}_{\LL^2} \| \Di{n} f \|^{\te}_{\LL^2} \quad
\mbox{for all $f \in \HH{n}$}. \label{wten} \feq
Eqs.\,\rref{cogive}\,\rref{finds} \rref{wten} give
\beq \HH{n} \subset \ZZ{2}{r}{j}\,, \quad \| \Di{j} f \|_{\LL^r} \leqs \Gp(j,n,\te)
\| f \|^{1 - \te}_{\LL^2} \| \Di{n} f \|^{\te}_{\LL^2}
\quad \mbox{for all $f \in \HH{n}$}~; \feq
summing up, \rref{gagg} holds and
$G(j,n,\te)$ is bounded from
above by $\Gp(j,n,\te)$.
\fine
\begin{prop}
\label{corquesta}
\textbf{Corollary.}
Let $j,n,\te,r$ be as in \rref{unoox};
consider the Sobolev inequality \rref{sobb}
$$ \| \Di{j} f \|_{\LL^r} \leqs S \sqrt{\| f \|^2_{\LL^2} + \| \Di{n} f \|^2_{\LL^2}}
\quad \mbox{for all $f \in \HH{n}$}. $$
This statement is true, and the related sharp constant $S(j,n,\te)$ has
the upper bound
\beq S(j,n,\te) \leqs \Sp(j,n,\te)~, \label{spp} \feq
$$ \Sp(j,n,\te) := \sqrt{(1 - \te)^{1 - \te} \te^\te}~ \Gp(j,n,\te),\qquad
\Gp(j,n,\te)\,\mbox{as in \rref{gpp}}. $$
The equality $S(j,n,\te) = \Sp(j,n,\te)$ holds in the
following cases: \parn
(a) $n >0$, $j = \te n$, where $S(j,n,\te)= \sqrt{(1-\te)^{1 - \te} \te^{\te}}$
due to Corollary \ref{corbyeq}; \parn
(b) $n >0$, $\te =1$, where $S(j,n,\te) = \Ga(j,n)$ as in \rref{sharpo}
due to Corollary \ref{corro}.
\end{prop}
\textbf{Proof.} Use Propositions \ref{questa} and \ref{soga}. \fine
Let us remark that, due to one of the Gamma function terms in
\rref{gpp}, the upper bound $\Gp(j,n,\te)$ diverges
if we fix $j,n$ and consider the limit $\te \vain
j/n + d/2 n$; the same can be said of the bound $\Sp(j,n,\te)$
defined by \rref{spp}. The case $\te = j/n + d/2 n$, which is
excluded from the conditions \rref{unoox},
is attacked with a different strategy in the next
subsection, where we even obtain the sharp constants.
\subsection{The $\boma{\LL^\infty}$ subcase}
\label{substre}
In this subsection we assume
\beq 0 < n < + \infty\,, \quad 0 \leqs j < + \infty\,,
\quad \te \equiv \te(j,n) := {j \over n} + {d \over 2 n} < 1 \label{unoos} \feq
and note that, with this choice of $\te$, the general definition of
$r$ in Eq.\,\rref{unoo} gives
\beq r = + \infty~. \feq
Our subsequent consideration will frequently refer to the
space $\LL^\infty = C_0(\reali^d, \complessi)$, see Eq.\,\rref{linfo}.
The inequalities \rref{gagg} \rref{sobb} read
\beq \HH{n} \subset \ZZ{2}{\infty}{j}~,~~\| \Di{j} f \|_{\LL^\infty}
\leqs G \, \| f \|^{1 - \te(j,n)}_{\LL^2} \| \Di{n} f \|^{\te(j,n)}_{\LL^2}
\label{gaggs} \feq
$$ \mbox{for some $G$ and all $f \in \HH{n}$}~, $$
\beq \HH{n} \subset \ZZ{2}{\infty}{j}~,~~ \| \Di{j} f \|_{\LL^\infty}
\leqs S \sqrt{ \| f \|^2_{\LL^2} + \| \Di{n} f \|^2_{\LL^2}} \label{sobbs} \feq
$$ \mbox{for some $S$ and all $f \in \HH{n}$}. $$
Differently from the previous subsections, here
we first give a result for \rref{sobbs} and then
present its implications for \rref{gaggs};
this approach gives the sharp constants for both inequalities.\par
In item (iv) of the Introduction we have already mentioned that the
inequalities \rref{gaggs} \rref{sobbs}
have been analyzed for $j=0$ by Ilyin \cite{Ily}, who determined
the sharp constants and the maximizers.
This author is mainly interested in Eq.\,\rref{gaggs},
so he does not present a specific statement about
\rref{sobbs}; however, he essentially derives Eq.\,\rref{sobbs}
with its sharp constant and uses this result as
an intermediate step towards Eq.\,\rref{gaggs}
({\footnote{For completeness, let us mention
that Eq.\,\rref{sobbs} for $j=0$ has a structure
very similar to the inequality
$\| f \|_{\LL^\infty} \leqs K \, \| \sqrt{1 + D^2}^{~n} f \|^2_{\LL^2}$,
for which we have given the sharp constant and
a maximizer in our previous work \cite{imb}\,.}}).
In this paragraph we generalize the results of Ilyin
to the case of arbitrary $j$; after deriving
the sharp constant and a maximizer for Eq.\,\rref{sobbs},
we use Proposition \ref{soga} to obtain the analogous results on \rref{gaggs}.
\par
Let us repeat another fact anticipated in item (iv) of the Introduction.
The maximizer for \rref{gaggs} \rref{sobbs} that we present coincides
for $j=0$ with the one determined by Ilyin
({\footnote{This is not at all surprising, since
Ilyin proves uniqueness of the maximizer
for \rref{gaggs} with $j=0$ up to
translation, rescaling and multiplication
by a constant.}}); besides giving its Fourier
representation (as in \cite{Ily} for $j=0$), we also
express it in terms of the space variables
using the Fox $H$-function or the Meijer $G$-function
(for arbitrary $j$).
The definitions of the above mentioned special functions are summarized in a
specific subsection of Appendix \ref{appe}, where we also give some basic
references about them. We remark that the $G$-function is a special
case of the $H$-function, more frequently implemented in standard packages
for symbolic or numerical computations
with special functions; due to this fact, we emphasize the use of $G$
whenever possible. For certain choices of the parameters, the $H$- or $G$-functions
considered hereafter are in fact
elementary functions; some examples are given in Section \ref{seces}. \par
\par
After all these preliminaries, we focus the attention on \rref{sobbs}.
\begin{prop}
\label{propsharp}
\textbf{Proposition.} For $j, n$ and $\te(j,n)$ as in \rref{unoos},
we have the following. \parn
(i) The Sobolev inequality \rref{sobbs} holds and its sharp constant is
\beq S(j,n) = {1 \over 2^{d/2} \pi^{d/4 - 1/2}
\sqrt{\Gamma(d/2) \, n \sin(\pi \, \te(j,n))}}~. \label{sharp} \feq
A maximizer for \rref{sobbs} is the function
\beq f := \FF^{-1}\Big( {|\k|^j \over 1 + |\k|^{2 n}} \Big)~. \label{maxim} \feq
The above $f$, being the inverse Fourier transform of an $\LL^1$ function,
is in $\LL^\infty$; moreover, $f$ can be expressed as follows using the Bessel function $J_{d/2 - 1}$:
\parn
\vbox{
\beq f = F_{j n}(|\x|)~,\quad ~F_{j n} \in C([0,+\infty),\reali)~, \label{bessel} \feq
$$ F_{j n}(\rho) :=  \int_{0}^{+\infty} \hspace{-0.5cm} d \xi\,
{J_{d/2 - 1}(\rho \,\xi) \over (\rho \,\xi)^{d/2 - 1}}
{ \xi^{d + j - 1} \over 1 + \xi^{2 n}}~\mbox{for $\rho>0$},~~
F_{j n}(0) = {\pi \over 2^{d/2} \Gamma(d/2) n \sin(\pi {j + d \over 2 n})}~. $$
}
(ii) The function $F_{j n}$ of Eq.\,\rref{bessel} can be expressed as a Fox $H$-function,
in the following way: for all $\rho \in [0,+\infty)$,
\beq F_{j n}(\rho) = {1 \over 2^{d/2} n}\,
H \left(\left. \barray{l} (1 - {j + d \over 2 n},{1 \over n});  \\
(0, 1),(1 - {j + d \over 2 n}, {1 \over n}); (1 - {d \over 2},1)  \farray \right|
\left({\rho \over 2}\right)^2\right). \label{hnj} \feq
$F_{j n}$ can also be expressed in terms of the Meijer $G$-function if
$n$ is rational. More precisely, if
\beq n = {N \over M} \qquad N, M \in \{1,2,3,...\} \feq
we have the following:
\beq F_{j n}(\rho) = {M \over 2^{d/2 + M -1} \, \pi^{M - 1} N^{d/2}}\,\,
G \left( \left. \barray{l} a_1,...,a_N ;  \\
b_1,...,b_{N+M}; \bs_1,...,\bs_N  \farray \right| \left(\rho \over 2 N\right)^{2 N} \right),
\label{gnj} \feq
where the parameters labeling the $G$-function are defined as follows:
\parn
\vbox{
\beq a_\ell := 1 - {j + d \over 2 N} - {\ell-1 \over M}~~ \mbox{for $\ell=1,...,M$}~; \label{gnjpa} \feq
$$ b_{h} := {h-1 \over N}~ \mbox{for $h=1,...,N$}~,~~
b_{N+h} := - {j+d \over 2 N} + {h \over M}~~ \mbox{for $h=1,...,M$}~; $$
$$ \bs_\ell := 1 - {d \over 2 N} - {\ell - 1 \over N}~~\mbox{for $\ell=1,...,N$}~. $$
}
\end{prop}
\vskip 0.1cm \noindent
\textbf{Proof.} In the sequel we frequently refer to the integral
\beq I(j,n) := \int_{\reali^d} dk \, {|k|^{2 j} \over 1 + |k|^{2 n}}~, \label{definj} \feq
which is finite due to the assumptions on $j, n$ in \rref{unoos} and given by
\beq I(j,n)  = {\pi^{d/2 + 1} \over
\Gamma(d/2) n \sin(\pi \te(j,n) )} \label{inj} \feq
(see Lemma \ref{lemap} in Appendix \ref{appe}
({\footnote{Use item (ii)
of this lemma with $a=j + d/2$ and $b= n$, so that (by \rref{unoos})
$b > a > 0$ and $a/b = \te(j,n)$.}})).
Our argument is divided in two steps. \parn
\textsl{Step 1. The Sobolev inequality \rref{sobbs} holds and its sharp constant $S(j,n)$ fulfills}
\beq S(j,n) \leqs {\sqrt{I(j,n)} \over (2 \pi)^{d/2}}~. \label{sharpp} \feq
To prove this we note that
\beq g \in \LL^1~ \Rightarrow \FF^{-1} g \in \LL^{\infty}, ~\| \FF^{-1} g \|_{\LL^{\infty}} \leqs
{1 \over (2 \pi)^{d/2}} \| g \|_{\LL^1}~; \label{note} \feq
this is the special case $p=1$ of the Hausdorff-Young inequality \rref{knowy}
(which is in fact derived by elementary means, since the relation
$\FF^{-1} g (x) = (2 \pi)^{-d/2} \int_{\reali^d} dk \, e^{i k x} g(k)$ implies
$|\FF^{-1} g (x)| \leqs (2 \pi)^{-d/2} \int_{\reali^d} dk \,|g(k)|$ for all $x$). \par
Let us fix $f \in \HH{n}$. Writing $\Di{j} f = \FF^{-1} (|\k|^j \FF f)$ and using \rref{note}
with $g = |\k|^j \FF f$ we obtain
\beq \Di{j} f \in \LL^{\infty}, \qquad\| \Di{j} f \|_{\LL^{\infty}} \leqs
{1 \over (2 \pi)^{d/2}} \| |\k|^j \FF f \|_{\LL^1}~, \label{into} \feq
provided that $|\k|^j \FF f \in \LL^1$. In order to
check this statement, we write
$$ |\k|^j \FF f = {|\k|^j \over \sqrt{1 + |\k|^{2 n}}} \, \sqrt{1 + |\k|^{2 n}} \, \FF f $$
and use the H\"older inequality; this ensures that $|\k|^j \FF f$ is actually in $\LL^1$, with
\parn
\vbox{
\beq \|  |\k|^j \FF f \|_{\LL^1} \leqs \Big\| {|\k|^j \over \sqrt{1 + |\k|^{2 n}}} \Big\|_{\LL^2}
\| \sqrt{1 + |\k|^{2 n}} \FF f \|_{\LL^2}  \label{hold} \feq
$$ =  \sqrt{I(j,n)} \sqrt{\| f \|^2_{\LL^2} + \| \Di{n} f \|^2_{\LL^2}} $$
}
(as for the last equality, see \rref{norm} and \rref{definj}).
Inserting \rref{hold} into \rref{into}, we conclude
\beq \Di{j} f \in \LL^\infty,\qquad\| \Di{j} f \|_{\LL^\infty} \leqs {\sqrt{I(j,n)} \over (2 \pi)^{d/2}}
\sqrt{\| f \|^2_{\LL^2} + \| \Di{n} f \|^2_{\LL^2}}~. \feq
This proves the inequality \rref{sobbs} and gives the bound \rref{sharpp} on its sharp constant.
\parn
\textsl{Step 2. Let $f := \FF^{-1}\left( {|\k|^j \over 1 + |\k|^{2 n}} \right)$, as
in \rref{maxim}. Then
$f$, being the Fourier transforms of an $\LL^1$ function, is in $\LL^\infty$;
$f$ possesses the features described by Eqs.\,\rref{bessel}-\rref{gnjpa}, and is in $\HH{n}$. Moreover,}
\beq S(j,n) \geqs {\| \Di{j} f \|_{\LL^\infty} \over \sqrt{\| f \|^2_{\LL^2} + \| \Di{n} f \|^2_{\LL^2}}}
= {\sqrt{I(j,n)} \over (2 \pi)^{d/2}}  \label{sharpm} ~.\feq
Indeed, $f$ fits the framework of Lemma \ref{lemmamej} in Appendix \ref{appe}
({\footnote{This lemma must be used with $a=j/2 + d/2$ and $b= n$.
It is $a,b > 0$ and $b >a$, since $a/b = \te(j,n) - j/2n < \te(j,n) < 1$.}}).
Being the Fourier transform of an $\LL^1$ function, $f$ is in $\LL^\infty$
and can be described via Eqs.\,\rref{bessel}-\rref{gnjpa} due to
the cited lemma. To go on, we note that Eq.\,\rref{norm}
and the convergence of the integral $I(j,n)$ in \rref{definj} yield the following statements
\beq \sqrt{1 + |\k|^{2 n}} \FF f = {|\k|^j \over \sqrt{1 + |\k|^{2 n}}} \in \LL^2, \quad
\mbox{whence $f \in \HH{n}$}~; \feq
\beq \sqrt{\| f \|^2_{\LL^2} + \| \Di{n} f \|^2_{\LL^2}} =
\| \sqrt{1 + |\k|^{2 n}} \FF f \|_{\LL^2} = \sqrt{I(j,n)}~. \label{imp1} \feq
Let us pass to
\beq D^j f =  \FF^{-1}\Big( {|\k|^{2 j} \over 1 + |\k|^{2 n}} \Big) \label{djf} \feq
(automatically in $\LL^\infty$, due to \rref{sobbs}).
This function also fits Lemma \ref{lemmamej} which ensures, amongst else,
\beq \| D^j f \|_{\LL^\infty} = D^j f(0) = {1 \over (2 \pi)^{d/2}} \int_{\reali^d} d k
{|k|^{2 j} \over 1 + |k|^{2 n}}~ = {I(j,n) \over (2 \pi)^{d/2}}~. \label{supf} \feq
Eqs.\,\rref{imp1}\,\rref{supf} imply
$$ {\| \Di{j} f \|_{\LL^\infty} \over \sqrt{\| f \|^2_{\LL^2} + \| \Di{n} f \|^2_{\LL^2}}}~
= {\sqrt{I(j,n)} \over (2 \pi)^{d/2}}~; $$
this ratio is obviously bounded from above by $S(j,n)$, so we have the thesis \rref{sharpm}.
\parn
\textsl{Conclusion of the proof.} Steps 1 and 2 indicate that
\beq S(j,n) = {\sqrt{I(j,n)} \over (2 \pi)^{d/2}} =
{\| \Di{j} f \|_{\LL^\infty} \over \sqrt{\| f \|^2_{\LL^2} + \| \Di{n} f \|^2_{\LL^2}}}~, \feq
with $f$ as in Eqs.\,\rref{maxim}\,\rref{bessel}. Thus $f$ is
a maximizer; finally, expressing $I(j,n)$ via \rref{inj} we obtain Eq.\,\rref{sharp}
for $S(j,n)$.
\fine
\begin{prop}
\textbf{Corollary.}
\label{corgaggs}
For $j, n$ and $\te(j,n)  \equiv \te$ as in \rref{unoos}, the Gagliardo-Nirenberg
inequality \rref{gaggs} holds and its sharp constant $G(j,n)$ is given by
\beq G(j,n)= {S(j,n) \over \sqrt{(1 -\te)^{1 - \te} \te^\te}}
= {1 \over 2^{d/2} \pi^{d/4 - 1/2} \sqrt{\Gamma(d/2)\,(1 -\te)^{1 - \te} \te^\te
\, n \sin(\pi \te)}}~. \label{sharpga} \feq
The function $f$ defined by Eq.\,\rref{maxim} is a maximizer for \rref{gaggs}.
\end{prop}
\textbf{Proof.} Use Propositions \ref{propsharp} and \ref{soga}. \fine
\begin{rema}
\label{refrem}
\textbf{Remark.} From \rref{sharp} \rref{sharpga} it is evident that
$S(j,n) = S(j',n)$ and $G(j,n) = G(j',n)$ if $(j, n), (j', n)$ fulfill conditions
\rref{unoos} and $\te(j,n) + \te(j',n) = 1$.
\fine
\end{rema}
\subsection{Another ``almost general'' $\boma{\LL^2}$ case: alternative
upper bounds for the Gagliardo-Nirenberg and Sobolev sharp constants}
\label{subsqua}
In this subsection we make the assumptions
\beq 0 \leqs \te < 1\,, \quad 0 \leqs n, j < + \infty\,, \quad
0 \leqs \te n - j \leqs {d \over 2}~~; \quad {1 \over r} = {1 \over 2} - {\te n - j \over d}~,
\label{unooy} \feq
which differ from the general $\LL^2$ conditions
\rref{unoo} since they exclude the case $\te=1$;
in this sense we are considering an almost general
$\LL^2$ case, slightly different from the one
of subsection \ref{subsdue}.
\par
Hereafter we use the strategy
of subsection \ref{substre}, i.e., we first derive a result on the Sobolev
inequality and then we point out its implications
for the Gagliardo-Nirenberg inequality; this approach
produces upper bounds for the sharp constants of both inequalities.
Let us note that the conditions \rref{unooy} differ from
the conditions \rref{unoox} of subsection \ref{subsdue} only at
boundary values; in the intersection
of \rref{unoox} with \rref{unooy}, we have two alternative
upper bounds for the Gagliardo-Nirenberg and Sobolev sharp constants.
A unified view of both types of upper bounds is proposed in subsection
\ref{subsqui}\,.
\begin{prop}
\textbf{Proposition.}
\label{propgen}
Let $j, n, \te, r$ be as in \rref{unooy}; consider
the Sobolev inequality \rref{sobb}
$$ \HH{n} \subset \ZZ{2}{r}{j},~\| \Di{j} f \|_{\LL^r} \leqs S \sqrt{\| f \|^2_{\LL^2} + \| \Di{n} f \|^2_{\LL^2}}~. $$
This is true and its sharp constant $S(j,n,\te)$ has the upper bound
\beq S(j,n,\te) \leqs \Spp(j,n,\te), \quad
\Spp(j,n,\te) := {E(j,n,\te) F(j,n,\te) \over \pi^{\te n/2 - j/2}}~; \label{ubound} \feq
here
\beq E(j,n,\te) := {(1 + 2 j/d - 2 \te n/d)^{d/4 + j/2 - \te n/2}
 \over (1 - 2 j/d + 2 \te n/d)^{d/4 - j/2 + \te n/2}}~;
\label{defe} \feq
\beq F(j,n,\te) :=
\left( {\Gamma\left( {d( 1 - \te) \over 2 (\te n - j)} \right) \Gamma\left( {d \te \over 2 (\te n - j)} \right)
\over n \,\Gamma\left({d \over 2} \right) \Gamma \left({d \over 2(\te n - j)}\right) } \right)^{{\te n - j \over d}}
\quad \mbox{if $j < \te n$}~, \label{eqF} \feq
$$ F(\te n,n,\te) := \sqrt{ (1 - \te)^{1 - \te} \te^{\te} } $$
(note that, if $n \neq 0$, $F(\te n,n,\te) = \lim_{j \vain (\te n)^{-}} F(j,n,\te)$). \parn
The equality $\Spp(j,n,\te)= S(j,n,\te)$ holds in the following subcases: \parn
(a) $\te=j/n$ with $n \neq 0$, where $S(j,n,\te)= \sqrt{(1-\te)^{1 - \te} \te^{\te}}$
due to Corollary \ref{corbyeq}; \parn
(b) $\te = j/n + d/2 n$ with $n \neq 0$, where $S(j,n,\te) =$ right hand side of Eq.\,\rref{sharp}.
\end{prop}
Before proving the above statements, we would like to point out their
connections with earlier works. In our previous paper \cite{imb}
we have considered the inequality $\| f \|_{\LL^r} \leqs K \|
\sqrt{1 + D^2}^{~n} f \|_{\LL^2}$, very similar to the case $j=0$ of
\rref{sobb}, and we have determined upper (and lower) bounds for the sharp constant
$K \equiv K(n,r)$. Hereafter we adapt the arguments
of \cite{imb}, replacing the norm
$\| \sqrt{1 + D^2}^{~n} f \|_{\LL^2}$ employed therein with
the equivalent norm $\sqrt{ \| f \|^2_{\LL^2} + \| \Di{n} f \|^2_{\LL^2}}$.
The difference between these norms has some implications from
the computational viewpoint; apart from this (and from the extension to the
case $j \neq 0$), the forthcoming proof of Proposition \ref{propgen} follows
the main ideas of \cite{imb}, e.g., a combination of the H\"older and
Hausdorff-Young inequalities to derive the upper bound in Eq.\,\rref{ubound}.
Let us also mention that, for $d=1$ and $j=0$, the upper bound of \rref{ubound}
has been (announced in \cite{Cot0} and) derived in \cite{Cot} along similar lines.
\vskip 0.1cm \noindent
\textbf{Proof of Proposition \ref{propgen}.} For $n=0$, all statements to be
proved are trivial; in fact, in this case we have $j=0$, $r=2$,
$\Spp(0,0,\te) = \sqrt{(1-\te)^{1 - \te} \te^\te} \geqs 1/\sqrt{2}$
and $1/\sqrt{2}$ is just the sharp constant of the Sobolev inequality, see
Remark \ref{remobvious}.
\par
From here to the end of the proof we assume $n \neq 0$: our argument
is divided in several steps.
\parn
\textsl{Step 1. Defining $\t,s$.}
In the sequel we denote with $\t$, $s$ the solutions of the equations
\beq {1\over r} + {1 \over \t} = 1~, \qquad {1 \over \t} = {1 \over 2} + {1 \over s} ~.
\label{st} \feq
Recalling that $r \in [2,+\infty]$, see \rref{remr}, we infer
from \rref{st} that
\beq \t \in [1,2]~, \qquad s \in [2,+ \infty]~. \feq
From the
explicit expression of $r$ in \rref{unooy} we readily obtain
\beq \t = {2 d \over d + 2 (\te n - j)}~, \qquad s = {d \over \te n - j}~. \feq
\textsl{Step 2. One has}
\beq {|\k|^j \over \sqrt{1 + |\k|^{2 n}}} \in \LL^s~, \qquad \Big\| {|\k|^j
\over \sqrt{1 + |\k|^{2 n}}} \Big\|_{\LL^s} = \pi^{\te n/2 - j/2} F(j,n,\te)
\label{tesf} \feq
\textsl{with $F$ as in \rref{eqF}; moreover,
$F(\te n,n,\te) = \lim_{j \vain (\te n)^{-}} F(j,n,\te)$}.
If $j < \te n$, we have $s\neq +\infty$ and one obtains \rref{tesf} using
Lemma \ref{lemap} in Appendix
\ref{appe} to evaluate $\int_{\reali^d} d k ({|k|^j \over \sqrt{1 + |k|^{2 n}}})^s$
({\footnote{To be precise, use item (i) of the cited lemma
with $a=j + d/s = \te n$, $b=n$ and $u=s/2 = d/2 (\te n - j)$.}}).
If $j = \te n$ we have $s = +\infty$;  ${|\k|^{\te n} \over
\sqrt{1 + |\k|^{2 n}}}$ is clearly continuous and vanishing at infinity,
hence in $\LL^\infty$, and $\| {|\k|^{\te n} \over
\sqrt{1 + |\k|^{2 n}}} \|_{\LL^\infty}$ $ = \sup_{\eta \in [0,+\infty)} {\eta^\te \over \sqrt{1 + \eta^{2}}}$ $
= \sqrt{ (1 - \te)^{1 - \te} \te^{\te} } = F(j,n,\te)$ (the sup is attained at
$\eta = \sqrt{\te/(1 - \te)}$).
Finally, the statement $F(\te n,n,\te) = \lim_{j \vain (\te n)^{-}}
F(j,n,\te)$ is checked expressing $\ln F(j,n,\te)$ via the Stirling formula
$\ln \Gamma(z)$ $= (z - 1/2) \ln z$ $- z + (1/2) \ln (2 \pi)$ $+ {1/(12 z)} +
O(1/z^2)$ for $z \vain + \infty$. \parn
\textsl{Step 3. The Sobolev inequality
\rref{sobb} holds, and its sharp constant $S(j,n,\te)$ fulfills}
\beq S(j,n,\te)
\leqs \Spp(j,n,\te) \label{thes3} \feq
\textsl{with $\Spp(j,n,\te)$ as in Eq.\,\rref{ubound}.}
Let us keep in mind the definitions of $\t,s$
in \hbox{Step 1}, and write down the Hausdorff-Young inequality \rref{knowy} with $p$, $p'$
replaced by $\t$, $r$; this gives (expressing $r$, $\t$ in terms of $j,n,\te$),
\beq g \in \LL^\t \quad\Rightarrow\quad \FF^{-1} g \in \LL^{\t}, ~\| \FF^{-1} g
\|_{\LL^{r}} \leqs {E(j,n,\te) \over \pi^{\te n  - j}} \, \| g
\|_{\LL^\t}\label{knowu} \feq
with $E$ as in \rref{defe}. Let us consider a function $f \in
\HH{n}$; writing $\Di{j} f = \FF^{-1} (|\k|^j \FF
f)$ and using \rref{knowu} with $g = |\k|^j \FF f$, we obtain
\beq \Di{j} f \in
\LL^{r},~~\| \Di{j} f \|_{\LL^{r}} \leqs {E(j,n,\te) \over \pi^{\te n - j}} \,  \|
|\k|^j \FF f \|_{\LL^\t}~, \label{intoo} \feq
provided that $|\k|^j \FF f \in \LL^\t$. In order to check this last statement, we write
$$ |\k|^j \FF f = {|\k|^j \over \sqrt{1 + |\k|^{2 n}}} \, \sqrt{1 + |\k|^{2 n}} \, \FF f $$
and recall the second relation \rref{st} connecting $\t, s$\,; this yields that
$|\k|^j \FF f$ is actually in $\LL^\t$, with
\beq \|  |\k|^j \FF f \|_{\LL^\t} \leqs \Big\| {|\k|^j \over \sqrt{1 + |\k|^{2 n}}} \Big\|_{\LL^s}
\| \sqrt{1 + |\k|^{2 n}} \FF f \|_{\LL^2}  \label{holdd} \feq
$$ = \pi^{\te n/2 - j/2} F(j,n,\te) \sqrt{\| f \|^2_{\LL^2} + \| \Di{n} f \|^2_{\LL^2}} $$
(as for the last equality, recall \rref{norm} and \rref{tesf}).
Inserting \rref{holdd} into \rref{intoo}, we conclude
\beq \Di{j} f \in \LL^r,~\| \Di{j} f \|_{\LL^r} \leqs
\Spp(j,n,\te) \sqrt{\| f \|^2_{\LL^2} + \| \Di{n} f \|^2_{\LL^2}}~, \feq
with $\Spp(j,n,\te)$ as in \rref{ubound}. This proves
the inequality \rref{sobb}
and the upper bound \rref{thes3} on its sharp constant. \parn
\textsl{Step 4. In the subcase $\te = j/n$,
where $S(j,n,\te)= \sqrt{(1-\te)^{1 - \te} \te^{\te}}$
due to Corollary \ref{corbyeq}, one has $\Spp(j,n,\te) = S(j,n,\te)$.}
This is readily checked using the definition
\rref{ubound} of $\Spp(j,n,\te)$.
\parn
\textsl{Step 5. In the subcase $\te = j/n + d/2 n$,
where $S(j,n,\te)$ equals the
right hand side of Eq.\,\rref{sharp}, one has
$\Spp(j,n,\te) = S(j,n,\te)$.}
This statement is checked using the definition
\rref{ubound} of $\Spp(j,n,\te)$, and recalling that
$\Gamma(z) \Gamma(1-z) = \pi/ \sin (\pi z)$.
\fine
\begin{prop}
\textbf{Corollary.}
\label{corgen}
For $j,n,\te,r$ as in \rref{unooy}, consider the Gagliardo-Nirenberg inequality
\rref{gagg}
$$
\HH{n} \subset \ZZ{2}{r}{j},\quad \| \Di{j} f \|_{\LL^r} \leqs
G \| f \|^{1 - \te}_{\LL^2} \| \Di{n} f \|^\te_{\LL^2}
\quad \mbox{for all $f \in \HH{n}$}~.
$$
This holds and its sharp constant $G(j,n,\te)$
has the upper bound
\beq G(j,n,\te) \leqs \Gpp(j,n,\te)~,\qquad
\Gpp(j,n,\te) := {\Spp(j,n,\te)\over \sqrt{(1 -\te)^{1 - \te} \te^\te}}~,
\label{bouggap} \feq
with $\Spp(j,n,\te)$ as in Eq.\,\rref{ubound}. \parn
The equality $\Gpp(j,n,\te) = G(j,n,\te)$ holds in
the cases: \parn
(a) $j = \te n$, where $G(j,n,\te)=1$ due to Proposition \ref{propgagbi}; \parn
(b) $\te = j/n + d/2 n$ (with $n \neq 0$), where $G(j,n,\te) = {1\over \sqrt{(1 -\te)^{1 - \te} \te^\te}}
\times$ right hand side of Eq.\,\rref{sharp}.
\end{prop}
\textbf{Proof.} For $n \neq 0$, everything follows from Propositions \ref{propgen} and \ref{soga}.
Let us pass to the case $n=0$, implying $j=0$ due to \rref{unooy}; then $j =\te n$ and
\rref{gaggs} holds with $G(j,n,\te)=1$ due to Proposition \ref{propgagbi}. On the
other hand, in this case $\Gpp(j,n,\te) =1$ due to the definitions \rref{bouggap} and
\rref{ubound}.
\fine
Let us remark that, due to the first $\Gamma$ function term
in \rref{eqF}, the upper bound $\Spp(j,n,\te)$
of \rref{ubound} diverges in the limit $\te \vain 1$;
the same happens for the bound $\Gpp(j,n,\te)$ in
\rref{bouggap}. On the other hand the special case $\te=1$,
excluded from the present conditions \rref{unooy},
has been already studied in subsection \ref{subsdue}\,.
\subsection{Putting together the results of subsections \ref{subsdue}, \ref{subsqua}:
upper bounds for the general $\LL^2$ case}
\label{subsqui}
In the cited subsections we have proved that the
Gagliardo-Nirenberg and Sobolev inequalities hold
for $j,n,\te$ as in \rref{unoox} or for $j,n,\te$ as
in \rref{unooy}. The union of these two cases
is just the ``general $\LL^2$ case'' described
by Eq.\,\rref{unoo}, i.e.,
\parn
\vbox{
$$ 0 \leqs \te \leqs 1\,, \quad 0 \leqs n, j < + \infty\,,
\quad 0 \leqs \te n - j   \leqs {d \over 2}\,,
\quad \te \neq 1~
\mbox{if $n = j + \dd{d \over 2}$}\,; $$
$$
{1 \over r} = {1 \over 2} - {\te n - j \over d}\,.
$$}
Summing up, we have the following result.
\begin{prop}
\textbf{Proposition.} The Gagliardo-Nirenberg and Sobolev
inequalities hold in the general $\LL^2$ case \rref{unoo}.
\end{prop}
Let us consider the sharp constants $G(j,n,\te)$, $S(j,n,\te)$.
In the ``almost general'' case \rref{unoox}
we have for them the upper bounds $\Gp(j,n,\te)$, $\Sp(j,n,\te)$ of
Proposition \ref{questa} and Corollary \ref{corquesta},
that coincide with the sharp constants
for $\te=1$ and diverge for $\te$ approaching the limit value
$j/n + d/2n$ (excluded from \rref{unoox}).
In the other ``almost general'' case
\rref{unooy}, we have the upper bounds $\Gpp(j,n,\te)$, $\Spp(j,n,\te)$
of Proposition \ref{propgen} and Corollary \ref{corgen},
that coincide with the sharp constants for $\te = j/n + d/2n$ and
diverge for $\te$ approaching the limit value $1$ (excluded from \rref{unooy}).
\par
Due to the above features, one expects $\Gp(j,n,\te)$, $\Sp(j,n,\te)$
to be better (i.e., smaller) than $\Gpp(j,n,\te)$, $\Spp(j,n,\te)$ for
$\te$ close to $1$, and the contrary to happen for $\te$ close
to $j/ n + d/2 n$. This is confirmed by the numerical values reported in
the next section for some sample choices of $j, n, \te$
(in space dimension $d=1,2,3$: see pages \pageref{pata3}-\pageref{pata4}).
Indeed, in these tests the $+$ bounds are
better than the $++$ bounds only for $\te$ very close to $1$.
\par
To conclude these considerations, let us recall that the $+$ and $++$ bounds
agree with the sharp constants even in the elementary case $j = \te n$
(with $n \neq 0$, if one considers the Sobolev inequality).
\subsection{Lower bounds for the sharp constants in
the general $\boma{\LL^2}$ case}
\label{subssei}
Let us still refer to the general $\LL^2$ framework
of Eq.\,\rref{unoo}. There is an obvious
strategy to obtain lower bounds on the sharp constants
$G(j,n,\te)$ and $S(j,n,\te)$:
one chooses any function $h \in \HH{n} \setminus \{ 0 \}$,
hereafter referred to as a ``trial function'', and notes that
\beq
G(j,n,\te) \geqs
{\| \Di{j} h \|_{\LL^r} \over \| h \|^{1- \te}_{\LL^2} \| \Di{n} h \|^\te_{\LL^2}}~, \label{gh1} \feq
\beq
S(j,n,\te) \geqs {\| \Di{j} h \|_{\LL^r} \over \sqrt{ \| h \|^2_{\LL^2}
+ \| \Di{n} h \|^2_{\LL^2}}} ~; \label{sh1} \feq
of course, in the choice of $h$ one should try to make the right hand side
of Eq.\,\rref{gh1} or \rref{sh1} as large as possible.
Hereafter we present alternative lower bounds which cover
the general case \rref{unoo}, and become very accurate for
$\te$ close to $1$. These are obtained using as a trial
function for the Gagliardo-Nirenberg inequality an
approximant of the maximizer given by \rref{fhyp} for
the case $\te=1$. From Corollary \ref{coronj},
we know that the maximizer in \rref{fhyp}
is in $\Hh{n}$, but it may fail to be in $\HH{n}$; the approximant introduced hereafter
is a regularization, depending on a parameter
$\ep >0$, that certainly belongs to $\HH{n}$.
Here is the statement implementing these ideas.
\begin{prop}
\label{progmm} \textbf{Proposition.} Let $j,n,\te,r$ be as in
Eq.\,\rref{unoo}. The sharp constant $G(j,n,\te)$ in the
Gagliardo-Nirenberg inequality \rref{gagg} has the lower bound
\beq G(j,n,\te) \geqs \Gm(j,n,\te \,|\,\ep)\quad \mbox{for all $\ep \in
(0,+\infty)$}, \label{lowboundep} \feq
which obviously implies
\beq G(j,n,\te) \geqs \Gm(j,n,\te) := \sup_{\ep \in (0,+\infty)}
\Gm(j,n,\te\,|\,\ep)~. \label{lowbound} \feq
Here we have put
\beq \Gm(j,n,\te \,|\,\ep)
:= {Y(j,n,\te \,|\,\ep) \over U(j,n \,|\, \ep)^{1-\te} V(j,n \,|\, \ep)^\te}~,
\label{gmm} \feq
where
\beq U(j,n \,|\,\ep) :=\left( {2 \pi^{d/2} \over \Gamma(d/2)}
\int_{0}^{+\infty} \hspace{-0.2cm} d \xi \, { \xi^{2 n- 2 j + d - 1} K^2_{n -
j}(\xi) \over (\xi^2 + \ep^2)^{2 n-j}} \right)^{1/2}, \label{znjte} \feq
\beq V(j,n \,|\, \ep) := \left( {2 \pi^{d/2} \over \Gamma(d/2)} \int_{0}^{+\infty}
\hspace{-0.2cm} d \xi \, { \xi^{4 n- 2 j + d - 1} K^2_{n - j}(\xi) \over (\xi^2
+ \ep^2)^{2 n-j}} \right)^{1/2} \label{vnjte} \feq
($K$ indicates, as usual, the Macdonald function); moreover
\beq Y(j,n,\te \,|\, \ep) := \left({2 \pi^{d/2}
\over \Gamma(d/2)} \int_{0}^{+\infty} \hspace{-0.2cm} d \rho \, \rho^{d-1}
|M_{j n \ep}(\rho)|^{r} \right)^{1/r}, \label{ynjep} \feq
where $M_{j n \ep} \in C([0,+\infty), \reali)$ is defined by
\beq M_{j n \ep}(\rho) :=
\int_{0}^{+\infty} \!\! d\xi~{J_{d/2 - 1}(\rho \,\xi) \over (\rho \,\xi)^{d/2-1}}\, {
\xi^{n + d -1} K_{n - j}(\xi) \over (\xi^2 + \ep^2)^{n-j/2}}
\label{defhnjep} \feq
(intending $ \dd{\left. {J_{d/2-1}(\rho\,\xi) \over
(\rho \,\xi)^{d/2-1}} \right|_{\rho=0} := \lim_{s \vain 0^{+}}
{J_{d/2-1}(s) \over s^{d/2-1}} = {1 \over 2^{d/2 - 1} \Gamma(d/2)}}$; note
that $Y$ depends on $\te$ through $r$\,).
\end{prop}
\textbf{Proof.}
The idea is to apply the inequality \rref{gh1}
choosing for $h$ the trial function
\beq h_{j n \ep} := \FF^{-1} g_{j n \ep}~, \
\qquad g_{j n \ep} := { |\k|^{n-j} K_{n - j}(|\k|) \over (|\k|^2 + \ep^2)^{n-j/2}}
\qquad (\ep > 0)~. \feq
In the limit $\ep \vain 0$ this function becomes (up to
a multiplicative constant) the maximizer
\rref{fhyp} for the case $\te=1$ of the
Gagliardo-Nirenberg inequality. In the sequel
we fix any $\ep >0$ and analyze the features
of the above function. \par
First of all, we note that
$g_{j n \ep}$ is in the space
of continuous, rapidly decreasing
functions on $\reali^d$; to check
this, one should recall that
$\lim_{\xi \vain 0^{+}} \xi^\mu K_{\mu}(\xi)$
is finite for any $\mu >0$, while $K_{\mu}(\xi) = O(e^{-\xi}/\sqrt{\xi})$
for $\xi \vain + \infty$ \cite{Nist}.
It is clear that $g_{j n \ep}, |\k|^n g_{j n \ep} \in \LL^2$, whence
$h_{j n \ep} \in \HH{n}$; moreover, by
the elementary rules for radial integrals
(Eq.\,\rref{then} of Appendix \ref{appe}),
\beq \| h_{j n \ep} \|_{\LL^2} = \| g_{j n \ep} \|_{\LL^2} =
U(j,n \,|\, \ep) \quad \mbox{as in \rref{znjte}}, \label{u1} \feq
\beq \| \Di{n} h_{j n \ep} \|_{\LL^2} = \| |\k|^n g_{j n \ep} \|_{\LL^2} =
V(j,n \,|\, \ep) \quad \mbox{as in \rref{vnjte}}~. \label{u2} \feq
Let us pass to
\beq \Di{j} h_{j n \ep} = \FF^{-1} (|\k|^j g_{j n \ep}) =
\FF^{-1} { |\k|^{n} K_{n - j}(|\k|) \over (|\k|^2 + \ep^2)^{n-j/2}}~.
\feq
The general theory of radial Fourier transforms,
see Appendix \ref{appe}
({\footnote{In particular,
Proposition \ref{furgen}.}}), gives
\beq \Di{j} h_{j n \ep} = M_{j n \ep}(|\x|) \qquad
M_{j n \ep}~\mbox{as in \rref{defhnjep}} \feq
and also ensures that $M_{j n \ep}$ is
continuous on $[0,+\infty)$.
The function $\Di{j} h_{j n \ep}$ is certainly in $\LL^r$
(this is established invoking the Gagliardo-Nirenberg inequality,
or more directly via the Hausdorff-Young inequality: in fact
$\FF \Di{j} h_{j n \ep}$ is continuous and rapidly decreasing,
hence in $\LL^{r'}$ with $1/r + 1/r'= 1$). We have
\beq  \| h_{j n \ep} \|_{\LL^r} = Y(j,n \,|\,\ep) \quad \mbox{as in
\rref{ynjep}}~. \label{u3} \feq
In conclusion, the inequality \rref{gh1} for the
trial function $h_{j n \ep}$ and Eqs.\,\rref{gmm} \rref{u1} \rref{u2} \rref{u3}  give
$$ G(j,n,\te) \geqs
{\| \Di{j} h_{j n \ep} \|_{\LL^r} \over
\| h_{j n \ep} \|^{1- \te}_{\LL^2} \| \Di{n} h_{j n \ep} \|^\te_{\LL^2}} =
{Y(j,n,\te \,|\,\ep) \over U(j,n \,|\, \ep)^{1-\te} V(j,n \,|\,\ep)^\te} =
\Gm(j,n,\te \,|\,\ep)~, $$
which is just the thesis \rref{lowboundep}.
\fine
\begin{prop}
\textbf{Corollary.} Let $j,n,\te,r$ be
as in Eqs.\,\rref{unoo}, and $n \neq 0$. $S(j,n,\te)$ has the lower bound
\beq S(j,n,\te)\geqs \Sm(j,n,\te \,|\,\ep) \quad \mbox{for all $\ep >0$}, \label{lowboundepso} \feq
which obviously implies
\beq S(j,n,\te) \geqs \Sm(j,n,\te) := \sup_{\ep \in (0,+\infty)} \Sm(j,n,\te \,|\,\ep)~;
\label{lowboundso} \feq
here we have put
\beq \Sm(j,n,\te \,|\,\ep) := \sqrt{(1 -\te)^{1 - \te} \te^\te}~\Gm(j,n,\te \,|\,\ep),
\quad \mbox{$\Gm$ as in \rref{gmm}}. \label{gss} \feq
\end{prop}
\textbf{Proof.} Use Propositions \ref{progmm} and \ref{soga}.
\fine
The previous Corollary excludes the trivial case $n=0$
of Eq.\,\rref{unoo}, implying $j=0$; let recall
that $S(0,0,\te)=1/\sqrt{2}$ due to Remark \ref{remobvious}.
\vskip 0.2cm
The forthcoming Proposition \ref{propgenlow} (with its Corollary
\ref{corgen2}) presents an alternative lower bound
for $S(j,n,\te)$ (and its equivalent for $G(j,n,\te)$),
holding under certain conditions on $\te$;
in this case we use \rref{sh1} choosing $h = f_{\lambda}$, where
$f$ is the maximizer of the special case $\te = j/n + d/2n$
(see Eq.\,\rref{maxim}) and the scaling parameter $\lambda$ is determined so as
to maximize the right hand side of \rref{sh1}. By construction,
this lower bound is accurate when $\te$ is close to the
special value $j/n + d/ 2 n$.
\begin{prop}
\textbf{Proposition.}
\label{propgenlow}
Let
\beq 0 < n < + \infty\,, \qquad 0 \leqs j \leqs n\,, \qquad
{j \over n} \leqs \te \leqs {j \over n} + {d \over 2 n} < 1\,;
\label{unoi} \feq
define $r$ as usual, via Eq.\,\rref{unoo}.
$S(j,n,\te)$ in Eq.\,\rref{sobb} has the lower bound
\parn
\vbox{
\beq S(j,n,\te) \geqs \Smm(j,n,\te)~, \label{esm} \feq
$$ \Smm(j,n,\te) := {I(j,n,\te) \over \pi^{d/4 + 1/2}}
\sqrt{\Gamma(d/2) n \sin( \pi(j/n + d/2 n)) {(1 - \te)^{1 - \te} \te^\te} \over
(1 - j/n - d/2 n)^{1-\te} (j/n + d/2 n)^\te }~; $$}
here
\beq I(j,n,\te) := \left\{ \barray{cc}
\dd{\pi \over 2^{d/2} \Gamma(d/2) \,n \,\sin(\pi \te )} & \mbox{if $\te = j/n + d/2 n$} \\
\left({2 \pi^{d/2} \over \Gamma(d/2)} \int_{0}^{+\infty} d\rho\, \rho^{d-1} |L_{j n}(\rho)|^r
\right)^{1/r} & \mbox{if $\te \neq j/n + d/2 n$} \farray \right. \label{definjte} \feq
and $L_{j n} \in C([0,+\infty),\reali)$ is defined in terms of the Bessel function $J_{d/2-1}$ by
\parn
\vbox{
\beq {~} \hspace{-0.4cm} L_{j n}(\rho) :=  \int_{0}^{+\infty} \hspace{-0.5cm} d \xi\,
{J_{d/2 - 1}(\rho \,\xi) \over (\rho \,\xi)^{d/2 - 1}}{ \xi^{d + 2 j -1} \over
1 + \xi^{2 n}}\quad\mbox{for $\rho>0$}, \label{defel} \feq
$$ L_{j n}(0) = {\pi \over 2^{d/2} \Gamma(d/2) \, n \sin(\pi {2 j + d \over 2 n})}~.$$
}
$L_{j n}$ can be expressed in terms of the Fox $H$-function or, if $n$ is rational,
in terms of the Meijer $G$-function; the expressions for $L_{j n}$ via $H$ or $G$
are equal to the ones given for the function $F_{j n}$ in Eqs.\,\rref{hnj}-\rref{gnjpa}, with $j$
replaced systematically by $2 j$. \parn
The equality $S(j,n,\te) = \Smm(j,n,\te)$ holds in the following cases: \parn
(a) $\te=0$, implying $j=0$, $r=2$ (due to \rref{unoi}) and $S(j,n,\te)=1$ (due to
Corollary \ref{corbyeq}); \parn
(b) $\te = j/n + d/2n$, where $S(j,n,\te) =$ right hand side of Eq.\,\rref{sharp}.
\end{prop}
\textbf{Proof.}
\textsl{Step 1. One has}
\beq S(j,n,\te) \geqs \Smm(j,n,\te) \label{thes4} \feq
\textsl{with $\Smm(j,n,\te)$ as in Eq.\,\rref{esm}; moreover
the function $L_{j n}$ of Eq.\,\rref{defel}, related to the definition of $\Smm(j,n,\te)$, has
all the features claimed in the proposition.}
As anticipated, this result is obtained using the
maximizer of the special case $\te = j/n + d/2 n$, i.e.,
the $\HH{n}$ function
$$ f := \FF^{-1}\left( {|\k|^j \over 1 + |\k|^{2 n}} \right) $$
of Eqs.\,\rref{maxim}\,\rref{bessel}. For any $\te$ as in \rref{unoi},
we apply \rref{sh1} with the trial function $h := f_{\lambda}$,
with $\lambda > 0$ unspecified for the moment.
This gives
\beq S(j,n,\te) \geqs {\| \Di{j} f_{\lambda} \|_{\LL^r} \over
\sqrt{ \| f_{\lambda}\|^2_{\LL^2} + \| \Di{n} f_{\lambda} \|^2_{\LL^2}} } \label{scas2} \feq
$$ =
{ \| \Di{j} f \|_{\LL^r} \over
\sqrt{\lambda^{-2 \te n} \| f \|^2_{\LL^2} + \lambda^{2 (1-\te) n} \| \Di{n} f \|^2_{\LL^2} } }
\quad \mbox{for all $\lambda > 0$}~, $$
the last equality following from \rref{scas}. To go on, we note that
\parn
\vbox{
\beq \| f \|_{\LL^2} = \left\| {|\k|^j \over 1 + |\k|^{2 n}} \right\|_{\LL^2} = \sqrt{A(j,n)}~, \label{451} \feq
$$ A(j,n) := {\pi^{d/2 + 1} \,(1 -j/n - d/2 n) \over \Gamma(d/2) \, n \sin( \pi (j/n + d/2 n))}~; $$
}
\parn
\vbox{
\beq \| \Di{n} f \|_{\LL^2} = \left\| {|\k|^{j+n} \over 1 + |\k|^{2 n}} \right\|_{\LL^2} = \sqrt{B(j,n)}~,
\label{452} \feq
$$ B(j,n) := {\pi^{d/2 + 1} \,(j/n + d/2 n) \over \Gamma(d/2) \, n \sin( \pi (j/n + d/2 n))}~. $$
}
For the expression of
the above norms in terms of $A(j,n)$ and $B(j,n)$, see Lemma \ref{lemap} in Appendix \ref{appe}
({\footnote{Use this lemma with the following choices:
$a=d/4 + j/2$, $b=n$, $u=2$ to get Eq.\,\rref{451},
$a=d/4 + j/2 + n/2$, $b=n$, $u=2$ to get Eq.\,\rref{452}.}}). Moreover
$\Di{j}f = \FF^{-1} \left({|\k|^{2 j} \over 1 + |\k|^{2 n}} \right)$ is in
$\LL^\infty$ and admits the representation
\beq \Di{j} f(x) = L_{j n}(|x|) \quad \mbox{for $x \in \reali^d$},\quad
L_{j n} \in C([0,+\infty),\reali)\quad\mbox{as in \rref{defel}}~; \label{besseldj} \feq
for this statement we refer to Lemma \ref{lemmamej}, also allowing to represent
$L_{j n }$ like $F_{j n}$ in Eqs.\,\rref{hnj}-\rref{gnjpa} with $j$ replaced by $2 j$
({\footnote{Use this lemma with $a=j + d/2$ and $b=n$.}}).
If $\te \neq j/n + d/2 n$ we have $r \neq \infty$, so
\beq \| \Di{j} f \|_{\LL^r} = \left({2 \pi^{d/2} \over \Gamma(d/2)} \int_{0}^{+\infty} \!\!\!d\rho \,\rho^{d-1}
|L_{j n}(\rho)|^r
\right)^{1/r}; \feq
if $\te = j/n + d/2 n$ we have $r = \infty$ and, due to
Lemma \ref{lemmamej}, Eq.\,\rref{eff0},
\beq \| \Di{j} f \|_{\LL^\infty} = (\Di{j} f)(0) =
{\pi \over 2^{d/2} \Gamma(d/2) \,n \,\sin(\pi \te )}~. \feq
In both cases,
we can write
\beq \| \Di{j} f \|_{\LL^r} = I(j,n,\te)\quad\mbox{as in Eq.\,\rref{definjte}}. \label{453} \feq
Inserting Eqs.\,\rref{451}\,\rref{452}\,\rref{453} into \rref{scas2} we obtain
\beq S(j,n,\te) \geqs { I(j,n,\te) \over
\sqrt{\lambda^{-2 \te n} A(j,n) + \lambda^{2 (1-\te) n} B(j,n) } }
\qquad \mbox{for all $\lambda > 0$}~, \feq
which implies
\beq S(j,n,\te) \geqs { I(j,n,\te) \over
\sqrt{\inf_{\lambda \in (0,+\infty)}
\big(\lambda^{-2 \te n} A(j,n) + \lambda^{2 (1-\te) n} B(j,n)\big) } }~.  \feq
The above $\inf$ equals ${A(j,n)^{1-\te} B(j,n)^\te \over (1 - \te)^{1 - \te} \te^\te}$
and is attained at $\lambda = \Big({\te\, A(j,n) \over (1-\te)\,B(j,n)}\Big)^{1/2n}$
(for $\te=0$, the last statement means that the $\inf$ is the
$\lambda \vain 0^{+}$ limit of the function). In conclusion,
\beq S(j,n,\te) \geqs I(j,n,\te) \sqrt{ {(1 - \te)^{1 - \te} \te^\te} \over
A(j,n)^{1-\te} B(j,n)^\te }  = \Smm(j,n,\te)~, \label{concl} \feq
where the last equality follows expressing $A(j,n), B(j,n)$ via
Eqs.\,\rref{451} \rref{452} and comparing with the definition
\rref{esm} of $\Smm(j,n,\te)$. \parn
\textsl{Step 2. If $\te=0$, Eq.\,\rref{unoo} implies $j=0$, $r=2$ due to
\rref{unoo} and one has $S(j,n,\te) = 1$ due to Corollary
\ref{corbyeq}; moreover $\Smm(j,n,\te) =1$.} Let $\te=0$. It is obvious that
\rref{unoo} implies $j=0$, $r=2$ and that the cited corollary gives
 $S(0,n,0)=1$. To go on, we note that
Eq.\,\rref{concl} gives $\Smm(0,n,0) := I(0,n,0)/\sqrt{A(0,n)}$;
on the other hand, Eqs.\,\rref{451}\,\rref{453}
imply $\sqrt{A(0,n)} = \| f \|_{\LL^2} = I(0,n,0)$,
so $\Smm(0,n,0)=1$.
\parn
\textsl{Step 3. If $\te = j/n + d/2 n$, one has $S(j,n,\te) = \Smm(j,n,\te) =$
the right hand side of Eq.\,\rref{sharp}.}
This statement is checked by direct comparison between the definitions
\rref{esm} of $\Smm(j,n,\te)$ and Eq.\,\rref{sharp}. An alternative
proof is as follows: the derivation of the bound
$S(j,n,\te) \geqs \Smm(j,n,\te)$ given in Step 1
is reduced, in the special case $\te = j/n + d/2 n$,
to the derivation of the bound
$S(j,n,j/n + d/2 n) \geqs$ right hand side of \rref{sharp}, given in Step 2
within the proof of Proposition \ref{propsharp}. (Note that
if $f$ is the function defined by Eq.\,\rref{maxim}
the fact that $f$ is a maximizer for $\te = j/n + d/2 n$
automatically ensures that in this subcase the ratio \rref{scas2} attains
its maximum at $\lambda =1$.) Thus $\Smm(j,n,j/n + d/2 n) =$
right hand side of \rref{sharp}; on the other hand, from Proposition \ref{propsharp}
we already know that $S(j,n,j/n + d/2 n) =$
right hand side of \rref{sharp}, whence the thesis.
\fine
\begin{prop}
\textbf{Corollary.}
\label{corgen2}
For $j$, $n$, $\te$, $r$ as in \rref{unoi}, $G(j,n,\te)$ fulfills the bounds
\beq G(j,n,\te) \geqs \Gmm(j,n,\te)~, \qquad
\Gmm(j,n,\te) := {\Smm(j,n,\te)
\over \sqrt{(1 -\te)^{1 - \te} \te^\te}}~,\label{bouggapm} \feq
with $\Smm(j,n,\te)$ as in Eq.\,\rref{esm}. \parn
The equality $G(j,n,\te) = \Gmm(j,n,\te)$ holds in
the cases: \parn
(a) $\te=0$, where $j=0$, $r=2$ (due to \rref{unoi}) and $G(j,n,\te)=1$ (due to
Proposition \ref{propgagbi}); \parn
(b) $\te = j/n + d/2 n$, where $G(j,n,\te) = {1\over \sqrt{(1 -\te)^{1 - \te} \te^\te}}
\times$ right hand side of Eq.\,\rref{sharp}.
\end{prop}
\textbf{Proof.} Use Propositions \ref{propgenlow} and \ref{soga}. \fine
\textbf{Comparison between the previous lower bounds.}
Let us repeat that the lower bounds $\Gmm(j,n,\te)$,
$\Smm(j,n,\te)$ are defined under the
condition \rref{unoi} (containing the
essential limitation $j/n + d/2 n < 1$),
while $\Gm(j,n,\te)$, $\Sm(j,n,\te)$ are defined in
the general $\LL^2$ case
\rref{unoo}. When both lower bounds $-$ and $--$
are defined, one naturally wonders which one is better (i.e., larger).
This depends on the choice of $j,n,\te$: see, for example,
the cases considered in Table III, page \pageref{pata4}
(and the introduction to this table on page \pageref{pata3}).
\section{Some examples}
\label{seces}
\textbf{Sharp constants and maximizers for the Sobolev
and Gagliardo-Nirenberg inequalities, in the $\boma{\LL^\infty}$ subcase
\rref{unoos}\,-\,\rref{sobbs}.}
The forthcoming Table I reports the values of
$\te(j,n)$ and of the sharp constants $S(j,n)$, $G(j,n)$ obtained
from Eqs.\,\rref{unoos}\,\rref{sharp}\,\rref{sharpga}, for $d=1,2,3$
and some choices of $j, n$. In this table, many lines
consider pairs $(j,n)$ and $(j',n)$ such that
$\te(j,n)+\te(j',n) = 1$, so that $S(j,n) = S(j',n)$ and $G(j,n) = G(j',n)$
(see Remark \ref{refrem}). For example, the second line
refers to the pairs $\dd{(0,{3/2})}$ and $\dd{({1/2},{3/2})}$,
such that $\dd{\te(0,{3/2})={1/3}}$ and $\dd{\te({1/2},{3/2})={2/3}}$.
\par
We know that, for $j,n$ as in \rref{unoos}, the inequality
\rref{gaggs} or \rref{sobbs} admits a maximizer $f = F_{j n}(|\x|)$,
where $F_{j n} : [0,+\infty) \vain \reali$, $\rho
\mapsto F_{j n}(\rho)$ is defined by Eq.\,\rref{bessel} in terms
of an integral involving a Bessel function. In a few cases
the integral is elementary, otherwise $F_{j n}$ can be expressed
in terms of a Fox $H$-function or Meijer $G$-function as in
Eqs.\,\rref{hnj}-\rref{gnjpa}. Table II reports, as examples, the functions
$F_{j n}$ for $d=1,2,3$ and some choices of $(j,n)$. The parameters
of the $G$-function for the cases appearing in Table II have
been determined using Eqs.\,\rref{gnj}\,\rref{gnjpa} and also taking into account
the symmetry properties mentioned in Appendix \ref{appe};
these allow to rearrange the sequences of parameters
$(a_{\ell})$, $(b_{\ell})$, $(\bs_{\ell})$ of Eq.\,\rref{gnjpa}
in increasing order and to eliminate
from the parameters any pair $(a_{\ell_1}, \bs_{\ell_2})$ with
$a_{\ell_1} = \bs_{\ell_2}$.
\begin{table}
\hskip -1.7cm
\textbf{Table I.} Values of $\te(j,n)$ and sharp constants
of the inequalities \rref{gaggs}\rref{sobbs}, given by \rref{unoos}\rref{sharp}\rref{sharpga}: some examples.
\vskip 1.0cm
{~}
\hskip -2cm
\footnotesize{
\begin{tabular}{|c||l||l|l|l|}
\hline
$d$ & $(j,n)$ \pb & $\te(j,n)$ & $S(j,n)$ & $G(j,n)$ \\
\hline
\hline
\multirow{4}{0.6cm}{$\barray{ccccc} ~\\~\\~\\~\\ 1 \farray$}
& $(0,1)$ \pb
& $\dd{1 \over 2}$
& $\dd{1 \over \sqrt{2}} = 0.707...$
& $1$ \\
\cline{2-5}
& $\dd{(0,{3 \over 2})}$ or $\dd{({1 \over 2},{3 \over 2})}$ \pb
& $\dd{1 \over 3}$ or $\dd{2 \over 3}$
& $\dd{{\sqrt{2} \over 3^{3/4}}} = 0.620... $
& $\dd{{2^{1/6} \over 3^{1/4}}} = 0.852...$ \\
\cline{2-5}
& $(0,2)$ or $(1,2)$\pb
& $\dd{1 \over 4}$ or $\dd{3 \over 4}$
& $\dd{{1 \over 2^{3/4}}} = 0.594...$
& $\dd{{2^{1/4} \over 3^{3/8}}} = 0.787...$ \\
\cline{2-5}
& $(4,10)$ or $(5,10)$ \pb
& $\dd{9 \over 20}$ or $\dd{11 \over 20}$
& $\dd{1 \over 2 \sqrt{5 \sin(9 \pi/20)}} = 0.224...$
& $\dd{1 \over 3^{9/20} 11^{11/40} \sqrt{\sin(9 \pi/20)}} = 0.317...$  \\
\hline
\multirow{4}{0.6cm}{$\barray{ccccc} ~\\~\\~\\~\\ 2 \farray$}
& $\dd{(0, {3 \over 2})}$ \pb
& $\dd{2 \over 3}$
& $\dd{1 \over 3^{3/4}} = 0.438...$
& $\dd{1 \over 2^{1/3} 3^{1/4}} = 0.603...$ \\
\cline{2-5}
& $(0,2)$ \pb
&  $\dd{1 \over 2}$
&  $\dd{{1 \over 2^{3/2}}} = 0.353... $
&  $\dd{{1 \over 2}}$ \\
\cline{2-5}
& $\dd{(0,{5 \over 2})}$ or $\dd{ ({1 \over 2},{5 \over 2}) }$ \pb
& $\dd{2 \over 5}$ or $\dd{3 \over 5}$
& $\dd{{1 \over \sqrt{5}} \left( {2 \over 5 + \sqrt{5}} \right)^{1/4}} = 0.324...$
& $\dd{{2^{1/20} \over 3^{3/10} (5 + \sqrt{5})^{1/4}}} = 0.453...$ \\
\cline{2-5}
& $(0,3)$ or $(1,3)$ \pb
& $\dd{1 \over 3}$ or $\dd{2 \over 3}$
& $\dd{{1 \over \sqrt{2} \, 3^{3/4}}} = 0.310...$
& $\dd{{1 \over 2^{5/6} 3^{1/4}}} = 0.426...$  \\
\hline
\multirow{4}{0.6cm}{$\barray{ccccc} ~\\~\\~\\~\\ 3 \farray$}
& $\dd{(0,2)}$ \pb
& $\dd{3 \over 4}$
& $\dd{1 \over 2^{5/4} \sqrt{\pi}} = 0.237...$
& $\dd{1 \over 2^{1/4} 3^{3/8} \sqrt{\pi}} = 0.314...$ \\
\cline{2-5}
& $\dd{(0,{5 \over 2})}$ \pb
&  $\dd{3 \over 5}$
& $\dd{{1 \over \sqrt{5} \, \pi} \left( {2 \over 5 + \sqrt{5}} \right)^{1/4}} = 0.182... $
& $\dd{{2^{1/20} \over 3^{3/10} (5 + \sqrt{5})^{1/4} \sqrt{\pi}}} = 0.256...$ \\
\cline{2-5}
& $\dd{(0,3)}$ \pb
& $\dd{1 \over 2}$
& $\dd{{1 \over 2 \sqrt{3} \, \pi}} = 0.162...$
& $\dd{{1 \over \sqrt{6 \pi}}} = 0.230...$ \\
\cline{2-5}
& $(1,3)$ \pb
& $\dd{5 \over 6}$
& $\dd{{1 \over \sqrt{6 \pi}}} = 0.230...$
& $\dd{{1 \over 5^{5/12} \sqrt{\pi}}} = 0.288...$  \\
\hline
\end{tabular}
}
\end{table}
\begin{table}
\textbf{Table II.} Maximizers for the inequalities \rref{gaggs} \rref{sobbs},
computed via \rref{bessel}\,-\,\rref{gnjpa}: \\ some examples.
\vskip 1.0cm
{~}
\footnotesize{
\begin{tabular}{|c||l||l|}
\hline
$d$ & $(j,n)$ \pb & $F_{j n}(\rho)$  \\
\hline
\hline
\multirow{3}{0.6cm}{$\barray{cccc} ~\\~\\~\\ 1 \farray$}
& $(0,1)$ \pb
& $\dd{\sqrt{\pi \over 2} \, e^{-\rho}}$ \\
\cline{2-3}
& $\dd{(0,2)}$ \pb & $\dd{ {\sqrt{\pi} \over 2}
\Big( \cos{\rho \over \sqrt{2}} - \sin{\rho \over \sqrt{2}} \Big)
\Big( \cosh{\rho \over \sqrt{2}} - \sinh{\rho \over \sqrt{2}} \Big)}$ \\
\cline{2-3}
& $(1,2)$ \pb
& $\dd{{1 \over 2} \, G \left( \left. \barray{l} 1/2 \,;  \\
0,1/2,1/2\,; \, 1/4,3/4  \farray \right| \left({\rho \over 4}\right)^{4} \right)}$ \\
\hline
\multirow{3}{0.6cm}{$\barray{cccc} ~\\~\\~\\ 2 \farray$}
& $\dd{(0,{3 \over 2})}$ \pb
& $\dd{{1 \over 6 \pi} \, G \left( \left. \barray{l} 1/6 \,;  \\
0, 1/6, 1/3, 2/3, 2/3\,; \, 0,1/3  \farray \right| \left({\rho \over 6}\right)^{6} \right)}$
\\
\cline{2-3}
& $(0,2)$ \pb
&  $\dd{{1 \over 4} \, G \left( \left. \barray{l}~  \,;  \\
0, 1/2, 1/2\,; \, 0  \farray \right| \left({\rho \over 4}\right)^{4} \right)}$\\
\cline{2-3}
& $\dd{(1,3)}$  \pb
& $\dd{{1 \over 6} \, G \left( \left. \barray{l} 1/2 \,;  \\
0, 1/3, 1/2, 2/3\,; \, 0,1/3, 2/3  \farray \right| \left({\rho \over 6}\right)^{6} \right)}$ \\
\hline
\multirow{3}{0.6cm}{$\barray{cccc} ~\\~\\~\\ 3 \farray$}
& $\dd{(0,2)}$ \pb
& $\dd{\sqrt{\pi \over 2} \,  \, {e^{-\rho/\sqrt{2}} \over \rho} \, \sin {\rho \over \sqrt{2} }}$  \\
\cline{2-3}
& $\dd{(0,3)}$ \pb
& $\dd{{1 \over 6 \sqrt{6}} \, G \left( \left. \barray{l}~  \,;  \\
0, 1/3, 1/2, 2/3 \,; \, -1/6, 1/6  \farray \right| \left({\rho \over 6}\right)^{6} \right)}$\\
\cline{2-3}
& $\dd{(1,3)}$ \pb
& $\dd{{1 \over 6 \sqrt{6}} \, G \left( \left. \barray{l} 1/3  \,;  \\
{0, 1/3, 1/3, 2/3} \,; \, {-1/6, 1/6, 1/2}  \farray \right| \left({\rho \over 6}\right)^{6} \right)}$ \\
\hline
\end{tabular}
}
\end{table}
\vfill \eject \noindent
{~}
\vskip -1.4cm \noindent
\textbf{Upper and lower bounds for the sharp constant in the Gagliardo-Nirenberg inequality
\rref{gagg}: numerical values in some examples.}
\label{pata3}
In the forthcoming Table III we present, for $d=1,2,3$: \parn
(i) some choices of $(j,n)$; \parn
(ii) the interval $\Theta_{j n} \equiv \Theta$
of the values of $\te$ fulfilling the conditions \rref{unoo}
for the chosen pair $(j,n)$; \parn
(iii) one or two sample choices for $\te$ in $\Theta$; \parn
(iv) the value of $r(j,n,\te) \equiv r$, see Eq.\,\rref{unoo}; \parn
(v) the numerical values of the lower bounds $\Gm$,
$\Gmm$ (Eqs.\,\rref{lowbound}\,\rref{bouggapm})
and of the upper bounds $\Gp$, $\Gpp$ (Eqs.\,\rref{gpp}\,\rref{bouggap}) for the sharp constants
($\Gm \equiv \Gm(j,n,\te),\, ...\, ,$ $\Gpp \equiv \Gpp(j,n,\te)$).
To be more precise, the quantities indicated in Table III as $\Gm,\,...\,, \Gpp$
are lower or upper approximants obtained numerically for the
theoretical bounds defined by \rref{lowbound} \rref{bouggapm} \rref{gpp} \rref{bouggap}.
\par
In the column about $\Gmm(j,n,\te)$, blank boxes refer to cases
where this bound is undefined because $j/n + d/2 n \geqs 1$.
For each choice of $(j,n,\te)$,
boldface is used to indicate the best
available lower bound (i.e., the maximum between $\Gm$ and $\Gmm$, when
both of them are defined)
and the best available upper bound (i.e., the minimum between $\Gp$ and $\Gpp$).\par
Let us add some information about the computation of the upper and lower bounds.
Concerning $\Gp$ and $\Gpp$, which have explicit expressions
in terms of Gamma and elementary functions, we have calculated their
numerical values, rounding up the results to the number of digits reported in the table.
\par
To compute the lower bound $\Gm$, which is
defined maximizing with respect to a parameter $\ep >0$,
we have operated in this way: firstly we have computed numerically all
integrals in Eqs.\,\rref{gmm}-\rref{defhnjep} for a grid of
sample values of $\ep$ (with spacing $1/100$), and then
we have taken the maximum over the grid (rounding down
to the digits reported in the table). This is, in any case,
a lower bound for the sharp constant (in fact smaller,
but not too much smaller than the sup in the definition \rref{lowbound} of $\Gm$).
\par
The lower bound $\Gmm$ depends on the integral
of a certain function $L_{j n} \equiv L$, see Eq.\,\rref{defel}. For the choices
of $j,n$ under consideration, this function
has been expressed as a Meijer G-function following
the indications after Eq.\,\rref{defel}; subsequently, its integral
has been computed numerically (and the final
value for $\Gmm$ obtained in this way has been
rounded down to the digits reported in the table).
\par
All the numerical calculations mentioned before have been
performed using {\tt{Mathematica}}.
For the choices of $j,n,\te$ in the table, we
have the following indications: \parn
(a) The upper bound $\Gpp$ is generally better than $\Gp$;
$\Gp$ is better only for $\te$ very close to $1$. This is
in agreement with the anticipations given on page \pageref{panti}. \parn
(b) The best lower and upper bounds are generally close (or even very
close), thus confining the sharp constant to a narrow interval.
Less satisfactory results are obtained when $j,n$ are large
and close (see, e.g., the case $d=1, j=9, n=10, \te = 37/40$, in which
the ratio lower bound/upper bound is, approximately, $0.53$). \par
To conclude, we point out that it would be easy to produce an analogue
of Table III for the Sobolev inequality \rref{sobb}, reporting
the lower or upper bounds
$S_{\sigma}(j,n,\te) \equiv S_{\sigma}$
of the previous section for the sharp constants ($\sigma = -,- -,+,++$); let us
recall that $S_{\sigma}(j,n,\te) = \sqrt{(1 -\te)^{1 - \te} \te^\te} \, G_{\sigma}(j,n,\te)$
for all $\sigma$. The previous comments about $G_{\sigma}$
(in particular, statements (a) (b)) could be repeated for the bounds $S_{\sigma}$.
\begin{table}[h]
\textbf{Table III.}\label{pata4} On the Gagliardo-Nirenberg inequality \rref{gagg}, in
some examples with $d=1,2,3$. $\Gm, \Gmm$ are the lower bounds
\rref{lowbound} \rref{bouggapm} on the sharp constants,
$\Gp$, $\Gpp$ are the upper bounds \rref{gpp} \rref{bouggap}
(all of them depending on $j,n,\te$).
Boldface is used to indicate the best
lower and upper bounds; see page \pageref{pata3} for more indications.
\vskip 0.5cm
{~}
\hskip -0.3cm
\footnotesize{
\begin{tabular}{|c||l|l|l|l||l|l||l|l|}
\hline
$d$ & $(j,n)$ \pb &~~~~~~~~~~~~ $\Theta$ & $\te$ & $r$ &
$\Gm$ & $\Gmm$ & $\Gp$ & $\Gpp$ \\
\hline
\hline
\multirow{6}{0.6cm}{$\barray{ccccc} ~\\~\\~\\1 \\ ~~ \farray$}
& $(0,1)$ \pbb
& $\{0 \leqs \te \leqs 1/2\}$
& 1/3
& 6
& \textbf{0.849}
& 0.832
& 1.204
& \textbf{0.873}
\\
\cline{2-9}
& $(3/4,1)$  \pbb
& $\{3/4 \leqs \te \leqs 1\}$
& 9/10
& 20/7
& \textbf{0.867}
& ~
& 1.030
& \textbf{0.944}
\\
\cline{2-9}
&  $(3/4,1)$ \pbb
& $\{3/4 \leqs \te \leqs 1\}$
& 99/100
&  50/13
& \textbf{0.950}
& ~
& \textbf{1.078}
& 1.564
\\
\cline{2-9}
&  $(1,2)$ \pbb
& $\{1/2 \leqs \te \leqs 3/4\}$
& 5/8
&  4
&  0.608
&  \textbf{0.633}
&  1.087
&  \textbf{0.711}
\\
\cline{2-9}
&  $(5,10)$ \pbb
& $\{1/2 \leqs \te \leqs 11/20\}$
& 21/40
&  4
&  0.080
&  \textbf{0.421}
&  1.087
&  \textbf{0.471}
\\
\cline{2-9}
&  $(9,10)$ \pbb
& $\{9/10 \leqs \te \leqs 19/20\}$
& 37/40
& 4
&  \textbf{0.317}
&  0.00894
&  1.087
&  \textbf{0.592}
\\
\hline
\multirow{6}{0.6cm}{$\barray{ccccc} ~\\~\\~\\2\\~ \farray$}
&  $(0,2)$ \pbb
& $\{0 \leqs \te \leqs 1/2\}$
& 1/3
&  6
&  \textbf{0.504}
&  0.498
&  0.741
&  \textbf{0.511}
\\
\cline{2-9}
&  $(0,3)$\pbb
& $\{0 \leqs \te \leqs 1/3\}$
& 1/6
& 4
& 0.533
& \textbf{0.547}
& 0.752
& \textbf{0.554}
\\
\cline{2-9}
&  $(1/2,1)$ \pbb
& $\{1/2 \leqs \te \leqs 1\}$
& 3/4
& 8/3
& \textbf{0.766 }
& ~
& 0.848
& \textbf{0.782 }
\\
\cline{2-9}
&  $(1/2,1)$ \pbb
& $\{1/2 \leqs \te \leqs 1\}$
& 9/10
& 10/3
& \textbf{0.714}
& ~
& \textbf{0.781}
& 0.795
\\
\cline{2-9}
&  $(1,3)$ \pbb
& $\{1/3 \leqs \te \leqs 2/3\}$
& 5/9
& 6
& 0.387
& \textbf{0.414}
& 0.741
& \textbf{0.436}
\\
\cline{2-9}
&  $(9,10)$ \pbb
& $\{9/10 \leqs \te < 1\}$
& 19/20
& 4
& \textbf{0.359}
& ~
& 0.752
& \textbf{0.504}
\\
\hline
\multirow{6}{0.6cm}{$\barray{ccccc} ~\\~\\~\\3\\~ \farray$}
& $(0,2)$\pbb
& $\{0 \leqs \te \leqs 3/4\}$
& 3/8
& 4
& \textbf{0.389}
& 0.359
& 0.494
& \textbf{0.394}
\\
\cline{2-9}
& $(0,3)$\pbb
& $\{0 \leqs \te \leqs 1/2\}$
& 1/3
& 6
& 0.273
& \textbf{0.278}
& 0.428
& \textbf{0.284}
\\
\cline{2-9}
& $(1,3)$ \pbb
& $\{1/3 \leqs \te \leqs 5/6\}$
& 2/3
& 6
& \textbf{0.264}
& 0.250
& 0.428
& \textbf{0.284}
\\
\cline{2-9}
& $(2,3)$ \pbb
& $\{2/3 \leqs \te \leqs 1\}$
& 95/100
& 60/13
& \textbf{0.385}
& ~
& 0.461
& \textbf{0.453}
\\
\cline{2-9}
& $(2,3)$\pbb
& $\{2/3 \leqs \te \leqs 1\}$
& 99/100
& 300/53
& \textbf{0.396}
& ~
& \textbf{0.433}
& 0.677
\\
\cline{2-9}
& $(9,10)$\pbb
& $\{9/10 \leqs \te \leqs 1\}$
& 19/20
& 3
& \textbf{0.321}
& ~
& 0.609
& \textbf{0.469}
\\
\hline
\end{tabular}
}
\end{table}
\vfill \eject \noindent
\textbf{Acknowledgments.} We are grateful to Alberto Cialdea and
Enrico Valdinoci for useful
bibliographical indications.
This work was partly supported by INdAM, INFN and by MIUR, PRIN 2010
Research Project  ``Geometric and analytic theory of Hamiltonian systems in finite and infinite dimensions''.
\appendix
\section{Appendix. Proof of Proposition \ref{lemga}}
\label{appega}
Let $p,q \in [1,+\infty]$, $j, n \in [0,+\infty)$, and assume
Eq.\,\rref{req1}
$$ 0 \leqs {1 \over r} = {1 \over q} - {n - j \over d} \leqs 1 ~. $$
The Gagliardo-Nirenberg inequality of
order $(p,q;j,n,\te=1)$ is statement \rref{gag1}, reported hereafter:
$$ \ZZ{p}{q}{n} \subset \ZZ{p}{r}{j},~~\| \Di{j} f \|_{\LL^r} \leqs G \| \Di{n} f \|_{\LL^q}
\quad \mbox{for some $G \in [0,+\infty)$ and all $f \in \ZZ{p}{q}{n}$}. $$
The corresponding extended Gagliardo-Nirenberg inequality is statement \rref{gagn1}, that we also report:
$$ \WW{q}{n} \subset \WW{r}{j},~~ \| \Di{j} f \|_{\LL^r} \leqs G \| \Di{n} f \|_{\LL^q}
\quad \mbox{for some $G \in [0,+\infty)$ and all $f \in \WW{q}{n}$}. $$
Let us recall the notations $G(p,q;j,n)$ and $G(q;j,n)$ for
the corresponding sharp constants. Proposition \ref{lemga} states the
equivalence of these inequalities and the equality of
their sharp constants for $q \neq 1,+\infty$.
\salto
\textbf{Proof of Proposition \ref{lemga}.}
\textsl{Step 1. The extended inequality \rref{gagn1}
implies the inequality \rref{gag1},
with $G(p,q; j,n) \leqs G(q;j,n)$.}
This follows immediately by comparing the two statements
under analysis and taking into account that
$\ZZ{p}{q}{n} = \LL^p \cap \WW{q}{n}$,
$\ZZ{p}{r}{j} = \LL^p \cap \WW{r}{j}$.
\parn
\textsl{Step 2. The inequality \rref{gag1} implies
the extended inequality \rref{gagn1}, with
$G(q;j,n) \leqs G(p,q;j,n)$.}
Let us assume the inequality \rref{gag1}, that we
write hereafter using its sharp constant:
\beq \ZZ{p}{q}{n} \subset \ZZ{p}{r}{j}~,~~\| \Di{j} f' \|_{\LL^r} \leqs G(p,q;j,n) \| \Di{n} f' \|_{\LL^q}
\quad \mbox{for all $f' \in \ZZ{p}{q}{n}$}. \label{asfo} \feq
Let us consider a function
\beq f \in \WW{q}{n}, \feq
which is fixed in the sequel. By the density of $\ZZ{p}{q}{n}$
in $\WW{q}{n}$ (see the final lines of Section
\ref{prelim}), there is a sequence
\beq f_\ell \in \ZZ{p}{q}{n} \qquad (\ell \in \naturali) \feq
such that, in the limit $\ell \vain \infty$,
\beq
\Di{n} f_{\ell} \vain \Di{n} f ~~ \mbox{in}~ \LL^q~. \label{a4} \feq
Since $\LL^q \hookrightarrow \Phi'$ we also have
$\Di{n} f_{\ell} \vain \Di{n} f$ in $\Phi'$; applying to
this relation the continuous map $\Di{-n} : \Phi' \vain \Phi'$ we obtain
\beq f_{\ell} \vain f ~~\mbox{in $\Phi'$}~. \label{a33} \feq
To go on, let us write the inequality in \rref{asfo} with
$f' = f_\ell$ or $f' = f_{\ell} - f_{\ell'}$, for $\ell, \ell' \in \naturali$;
this gives the following:
\beq \| \Di{j} f_{\ell} \|_{\LL^r} \leqs G(p,q;j,n) \| \Di{n} f_{\ell} \|_{\LL^q} \label{a5} \feq
\beq \| \Di{j} f_{\ell} - \Di{j} f_{\ell'} \|_{\LL^r}
\leqs G(p,q;j,n) \| \Di{n} f_{\ell} - \Di{n} f_{\ell'} \|_{\LL^q}~.
\label{a6} \feq
Eqs.\,\rref{a4} and \rref{a6} imply
\beq \| \Di{j} f_\ell - \Di{j} f_{\ell'} \|_{\LL^r} \vain 0
\qquad \mbox{for $\ell, \ell' \vain +\infty$}~, \feq
yielding the existence of
\beq f^j := \lim_{\ell \vain + \infty} \Di{j} f_\ell~~ \mbox{in}~ \LL^r~. \label{lr} \feq
On the other hand, Eq.\,\rref{a33} and the continuity of $\Di{j}: \Phi' \vain \Phi'$ give
\beq \Di{j} f_\ell \vain \Di{j} f~\mbox{in}~\Phi'~; \label{tcomp} \feq
noting that \rref{lr} implies $\Di{j} f_\ell \vain f^j$ in $\Phi'$, by comparison
with \rref{tcomp} we obtain $f^j = \Di{j} f$ and, returning to \rref{lr},
\beq \Di{j} f = \lim_{\ell \vain + \infty} \Di{j} f_\ell~~ \mbox{in $\LL^r$}~. \label{llr} \feq
Now, sending $\ell$ to $+\infty$
in Eq.\,\rref{a5} and recalling Eqs.\,\rref{a4}\,\rref{llr}, we obtain
\beq \| \Di{j} f \|_{\LL^r} \leqs G(p,q;j,n) \| \Di{n} f \|_{\LL^q}~. \label{a55} \feq
These conclusions hold for an arbitrary $f \in \WW{q}{n}$; thus,
the extended inequality \rref{gagn1} is true and its
sharp constant is such that $G(q;j,n) \leqs G(p,q;j,n)$. \parn
\textsl{Step 3. Conclusion of the proof.} Steps 1 and 2 clearly give the thesis. \fine
\section{Appendix. On certain integrals and inverse Fourier transforms.
The $\boma{H}$- and $\boma{G}$-functions}
\label{appe}
\textbf{Radial integrals.}
Let us consider a measurable function $g$ on $\reali^d$ of the form
\beq g(k) = G(|k|) \qquad \mbox{for $k \in \reali^d \setminus \{ 0 \}$}~, \label{a1} \feq
where
\beq G : (0,+\infty) \vain \complessi ~\mbox{measurable and}~
\int_{0}^{+\infty} \hspace{-0.2cm} d \xi \, \xi^{d-1}
\, |G(\xi)| < + \infty~. \label{a3} \feq
Then, as well known, $g \in \LL^1$ and
\beq \int_{\reali^d} \hspace{-0.2cm} d k \, g(k) = {2 \pi^{d/2} \over \Gamma(d/2)}
\int_{0}^{+\infty} \hspace{-0.2cm} d \xi \, \xi^{d-1} G(\xi)~. \label{then} \feq
(More generally, if $g$ has the form \rref{a1} with a meaurable
$G : (0,+\infty) \vain [0,+\infty)$, Eq. \rref{then} holds with both
sides possibly equal to $+\infty$). \par
From the above statement one infers the following lemma,
used in the main text; this refers to the Pochhammer symbol $(z)_{\ell}$
(see, e.g., \cite{Nist}), which is defined
as follows for $z \in \complessi$ and $\ell \in \{0,1,2,...\}$:
\beq (z)_0 := 1~, \qquad (z)_{\ell} := \prod_{i=0}^{\ell-1} (z + i)~
\mbox{for $\ell=1,2,...$}~. \feq
\begin{prop}
\label{lemap}
\textbf{Lemma.}
(i) Let $a, b,u \in \reali$, $b > a >0$ and $u > 0$; one has
\beq \int_{\reali^d}
\hspace{-0.2cm}  d k {|k|^{2 a u - d} \over (1 + |k|^{2 b})^u} = {2 \pi^{d/2}
\over \Gamma(d/2)} \int_{0}^{+\infty} \hspace{-0.2cm} d \xi {\xi^{2 a u - 1}
\over (1 + \xi^{2 b})^{u}}  \label{simi0} \feq
and
\beq \int_{0}^{+\infty} d \xi {\xi^{2 a u - 1} \over (1 + \xi^{2 b})^u}
= {\Gamma(a u/b) \Gamma(u - a u/b) \over 2 b \,\Gamma(u)}~. \label{simi}
\feq
(ii) In particular, let $a, b \in \reali$ with $b > a > 0$ and $u \in \{1,2,3,...\}$; then
\beq \int_{0}^{+\infty} d \xi {\xi^{2 a u - 1} \over (1 + \xi^{2 b})^u} =
{\pi (1 - au/b)_{u-1} \over 2 b (u-1)! \sin(\pi a u/b)} \label{int}~. \feq
In the special case $a u /b= m$ $\in \{1,2,...,u-1\}$,
the last formula should be applied understanding
\beq {\pi (1 - m)_{u-1} \over \sin(\pi m)} := \lim_{\ep \vain 0}
{\pi (1 - m - \ep)_{u-1} \over \sin(\pi(m + \ep))} = (-1)^{m+1} \prod^{u-2}_{i=0, i \neq m-1} (1 - m + i)~.
\label{limit} \feq
\end{prop}
\textbf{Proof.} (i) Eq.\,\rref{simi0} follows from Eq.\,\rref{then}.
Moreover, with a change of variable $\xi = t^{1/2 b}$ we get
$$ \int_{0}^{+\infty} d \xi {\xi^{2 a u - 1} \over (1 +
\xi^{2 b})^u} = {1 \over 2 b} \int_{0}^{+\infty} d t {t^{a u/b - 1} \over (1 + t)^u} =
{\Gamma(a u/b) \Gamma(u - a u/b) \over 2 b \,\Gamma(u)} $$
(for the last statement see, e.g., \cite{Nist}, Eq.\,(5.12.3)); this proves Eq.\,\rref{simi}. \parn
(ii) Let us apply Eq.\,\rref{simi} with $u \in \{1,2,3,...\}$, assuming
provisionally that $a u/b$ is noninteger.
The known relations
$\Gamma(z + \ell) = (z)_{\ell} \Gamma(z)$ (for $\ell \in \naturali$) and
$\Gamma(z) \Gamma(1-z) = \pi/\sin(\pi z)$ (for $z$ noninteger) give
\beq \Gamma(a u/b) \Gamma(u - a u/b) = \Gamma(a u/b) \Gamma(1 - a u/b + (u-1)) \label{eqsub} \feq
$$ = (1 - a u/b)_{u-1} \Gamma(a u/b) \Gamma(1 - a u/b) = {\pi (1 - a u/b)_{u-1} \over \sin(\pi a u/b)}~. $$
Substituting Eq.\,\rref{eqsub} into \rref{simi}, and writing $\Gamma(u) = (u-1)!$,
we obtain the thesis \rref{int}. \par
It remains to analyze the special case $u \in \{1,2,3,...\}$,
$a u/b = m \in \hbox{$\{1,2,...,m-1\}$}$. Elementary considerations of continuity
show that Eq.\,\rref{int} holds again if one defines
${\pi (1 - m)_{u-1} \over \sin(\pi m)}$ as the $\ep \vain 0$ limit appearing
in Eq.\,\rref{limit}; on the other hand,
\beq {\pi (1 - m - \ep)_{u-1} \over \sin(\pi(m + \ep))} =
{\pi \prod_{i=0}^{u-2} (1 - m - \ep + i)  \over \sin(\pi(m + \ep))} \feq
$$ = - {\pi \ep \over \sin(\pi(m + \ep))} \prod_{i=0, i \neq m-1}^{u-2} (1 - m - \ep + i)
\underset{\ep \vain 0}{\longrightarrow} (-1)^{m+1} \prod_{i=0, i \neq m-1}^{u-2} (1 - m + i)~, $$
which justifies the second equality in \rref{limit}. \fine
\textbf{General results on radial, inverse Fourier transforms.} Let us summarize some standard facts.
\begin{prop}
\label{furgen}
\textbf{Proposition.} Let $g$ be as in \rref{a1} \rref{a3}, hence in $\LL^1$, and
$f := \FF^{-1} g \in \LL^\infty$.
$f$ has the radial structure
\beq f(x) = F(|x|)~\mbox{for $x \in \reali^d$},~~ \label{f1} \feq
where $F \in C([0,+\infty),\complessi)$ is given by
\beq F(\rho) := \int_{0}^{+\infty}
d \xi~\xi^{d-1} {J_{d/2 - 1}(\rho\, \xi) \over (\rho\, \xi)^{d/2-1}} \, G(\xi) \quad
\label{ff1} \feq
$$ \mbox{(intending $ \dd{\left. {J_{d/2-1}(\rho\, \xi) \over (\rho \,\xi)^{d/2-1}} \right|_{\rho = 0} :=
\lim_{s \vain 0^{+}} {J_{d/2-1}(s) \over s^{d/2-1}} = {1 \over 2^{d/2 - 1} \Gamma(d/2)}}$)}~. $$
Moreover, if \hbox{$G(\xi) \geqs 0$}
for all $\xi \in (0,+\infty)$, then
\beq \| f \|_{\LL^{\infty}} = f(0) = F(0)~. \label{ofco} \feq
\end{prop}
\textbf{Proof.} The representation
\rref{f1}\rref{ff1} for $f$ is well known, see, e.g., \cite{Boc}.
Of course, we have
\beq f(x) = {1 \over (2 \pi)^{d/2}} \int_{\reali^d} dk \, e^{ik \sc x} G(|k|) \label{fx} \feq
for each $x \in \reali^d$ and, in particular,
\beq f(0) = {1 \over (2 \pi)^{d/2}} \int_{\reali^d} d k \,G(|k|)~; \label{f0} \feq
if $G(\xi) \geqs 0$ for all $\xi$, Eqs.\,\rref{fx} \rref{f0} imply
\beq |f(x)| \leqs {1 \over (2 \pi)^{d/2}} \int_{\reali^d} d k \, G(|k|) = f(0)~
\mbox{for all $x \in \reali^d$}, \feq
yielding Eq.\,\rref{ofco}.
\fine
\textbf{The Fox and Meijer functions $\boma{H,G}$.}
Before computing some inverse Fourier transforms, it is
necessary to say a few words on the above cited functions. \par
The Fox $H$-function \cite{Mat} \cite{Wolf} depends
on one complex variable $z$ and on a set of
parameters $a_\ell, A_{\ell}$ ($\ell=1,...,n$), $\as_\ell, \As_{\ell}$ ($\ell=1,...,\ns$),
$b_\ell, B_{\ell}$ ($\ell=1,...,m$), $\bs_\ell, \Bs_{\ell}$ ($\ell=1,...,\ms$)
where $n,\ns,m,\ms$ are nonnegative integers, $a_\ell,\as_{\ell}, b_{\ell}, \bs_{\ell}$
are complex numbers and $A_\ell, \As_\ell, B_\ell, \Bs_\ell$ are real and positive.
The definition reads
\vskip 0.1cm \noindent
\vbox{
\beq H \left(\!\!\! \left.\barray{l} (a_1,A_1),...,(a_n,A_n) ; (\as_1,\As_1)...,(\as_\ns, \As_\ns) \\
(b_1, B_1),...,(b_m, B_m); (\bs_1,\Bs_1),...,(\bs_\ms, \Bs_\ms)  \farray \right| z \! \right)  \label{hfox} \feq
$$ :={1\over 2 \pi i} \int_{\L} d s  \, z^s
{ \prod_{\ell=1}^n \Gamma(1 - a_\ell + A_\ell s) \prod_{\ell=1}^m \Gamma(b_\ell - B_\ell s)  \over
\prod_{\ell=1}^\ns \Gamma(\as_\ell - \As_\ell s) \prod_{\ell=1}^\ms \Gamma(1 - \bs_\ell + \Bs_\ell s) }~, $$
}
where $\L$ is an oriented path in $\complessi$ such that the ``left'' and the ``right'' poles of the integrand
in \rref{hfox} are on the left and on the right of $\L$, respectively. By definition, the sets
of the left and right poles are
\parn
\vbox{
\beq \LL \equiv \LL((a_1,A_1),...,(a_n,A_n)) := \label{eqleft} \feq
$$ \{ \sigma \in \complessi~|~\sigma~\mbox{is a pole of the function $s \mapsto \Gamma(1 - a_\ell + A_\ell s)$ for
some $\ell \in \{1,...,n\}$} \} $$
$$ = \{ {- 1 + a_\ell - k \over A_{\ell}}~|~\ell=1,...,n;~ k =0,1,2,...\}~; $$
}
\vbox{
\beq \RR \equiv \RR((b_1,B_1),...,(b_m,B_m)) := \label{eqright} \feq
$$ \{ \sigma \in \complessi~|~\sigma~\mbox{is a pole of the function $s \mapsto \Gamma(b_\ell - B_\ell s)$ for
some $\ell \in \{1,...,m\}$} \} $$
$$ = \{ {b_\ell + k \over B_{\ell}} ~|~\ell=1,...,m;~ k =0,1,2,...\}~; $$
}
the parameters $a_\ell, A_\ell, b_\ell, B_\ell$, must be such that $\LL$ and $\RR$ do not intersect.
\par
The path $\L$ must be conveniently specified, and suitable conditions must be
put on the parameters, on $z$ and on $\arg z$ to ensure a nonambiguous
definition of $z^s$ and the convergence of the integral in Eq.\,\rref{hfox}.
For the sake of the present paper, it suffices to consider the following choices and
conditions
({\footnote{Here and in the sequel the expression
``$\L$ is a path from $c - i \infty$ to $c + i \infty$'',
with $c \in \reali$, must be understood as follows:
$\L$ has a parametrization $s = \param(t)$, where $\param
\in C^1(\reali, \complessi)$ is such that
$\param(t) = c + i t + O(1/t)$ and $\param'(t) = i + O(1/t)$ for $t \vain \pm \infty$.
If this happens and $\phi, \psi$ are two complex functions
defined on the support of $\L$, we write ``$\phi(s) \sim \psi(s)$
for $s \vain c \pm i \infty$ along $\L$'' to mean the following:
for any parametrization $\param$ as before, one has $\phi(\param(t))/\psi(\param(t))
\vain 1$ for $t \vain \pm \infty$.}}):
\beq \L := \mbox{any path from $c - i \infty$ to $c + i \infty$ ($c \in \reali)$} \label{con1} \feq
$$ \mbox{with $\LL$ on its left, $\RR$ on its right}~; $$
\beq z \in \complessi \setminus (-\infty,0]~, \qquad |\arg z| < \pi~; \label{argz} \feq
\beq \alf >0~\mbox{and}~ |\arg z| < {\alf \pi \over 2} ~\mbox{if $\alf < 2$}~,
\label{c3a} \feq
\beq \alf := \sum_{\ell=1}^n A_{\ell} - \sum_{\ell=1}^\ns \As_{\ell}
+ \sum_{\ell=1}^m B_{\ell} - \sum_{\ell=1}^\ms \Bs_{\ell}~. \label{c3b} \feq
The second condition in Eq.\,\rref{argz} is understood as the definition of $\arg z$
to be used in Eq.\,\rref{hfox}, where one intends $z^s := e^{s (\log |z| + i \arg z)}$.
Due to condition \rref{c3a}, the integrand of Eq.\,\rref{hfox}
decays exponentially ({\footnote{In fact, let $\phi(s)$ denote the function
under the sign of integral in Eq.\,\rref{hfox}, also depending on $z$
and on the parameters $a_\ell,...,\Bs_{\ell}$. Using the Stirling's
formula for the Gamma function, one readily checks that,
for $s \vain c \pm i \infty$ along $\L$, $ |\phi(s)| \sim \Phi_{\pm}
|\mbox{Im} s|^{\dd{- \omega + c \,\Omega}} \,
e^{\dd{- (\pi \alf/2 \pm \arg z) |\mbox{Im} s|}}$ \,
where $\alf$ is defined by \rref{c3b},
$\omega :=
\mbox{Re}(\sum_{\ell=1}^n a_{\ell} + \sum_{\ell=1}^\ns \as_{\ell}
- \sum_{\ell=1}^m b_{\ell} - \sum_{\ell=1}^\ms \bs_{\ell}$),~
$\Omega :=
\sum_{\ell=1}^n A_{\ell} + \sum_{\ell=1}^\ns \As_{\ell}
- \sum_{\ell=1}^m B_{\ell} - \sum_{\ell=1}^\ms \Bs_{\ell}$
and $\Phi_{\pm}$ are constants that could  be written as well in
terms of $z,a_{\ell},...,\Bs_{\ell}$. The coefficients
$\pi \alf/2 \pm \arg z$ in the exponential are both positive if
and only if $\alf$ and $\arg z$ fulfill the inequalities \rref{c3a}.}}).
\par
One could implement the definition \rref{hfox} of $H$
with other choices of the path $\L$ and other conditions on the parameters and on $z$
\cite{Mat}, but these alternatives are not considered in this paper.
\par
The Meijer G-function \cite{Fik} \cite{Luke} \cite{Nist}  is an $H$-function with
$A_\ell, B_\ell, \As_\ell, \Bs_\ell=1$ for all
$\ell$; its definition reads
\parn
\vbox{
\beq G \left( \!\!\!\left. \barray{l} a_1,...,a_n ; \as_1,...,\as_\ns \\
b_1,...,b_m; \bs_1,...,\bs_\ms  \farray \right| z \!\right) \label{meij} \feq
$$ := {1\over 2 \pi i} \int_{\L} d s \, z^s
{ \prod_{\ell=1}^n \Gamma(1 - a_\ell + s) \prod_{\ell=1}^m \Gamma(b_\ell - s)  \over
\prod_{\ell=1}^\ns \Gamma(\as_\ell - s) \prod_{\ell=1}^\ms \Gamma(1 - \bs_\ell + s) }~. $$
}
The path $\L$ in the above integral is chosen as in Eq.\rref{con1},
and conditions \rref{argz} \rref{c3a} are prescribed.
Of course, in the present case the left and right poles $\LL,\RR$
and the parameter $\alf$ are described by Eqs.\,\rref{eqleft}\,\rref{eqright}\,\rref{c3b}
with $A_\ell, ...,\Bs_\ell=1$; in particular,
the last of these equations becomes
\beq \alf := n - \ns + m - \ms~. \label{c3bg} \feq
In the previous considerations it has been assumed $z \neq 0$;
the $H$- or $G$- functions at $z =0$ are defined taking the $z \vain 0$ limit,
if this exists.
\par
The $H$- and $G$-functions have certain symmetry properties with respect
to their parameters, which appear from the definitions \rref{hfox}
\rref{meij}. For example, from \rref{meij} it is evident that the $G$-function
is invariant under a permutation of anyone of the four
sequences $(a_1,...,a_n)$, $(\as_1,...,\as_\ns)$,
$(b_1,...,b_m)$, $(\bs_1,...,\bs_\ns)$. Moreover,
a $G$-function does not change if we remove from
the list of its parameters any pair $(a_{\ell_1}, \bs_{\ell_2})$
with $a_{\ell_1} = \bs_{\ell_2}$, or any pair $(\as_{\ell_1}, b_{\ell_2})$
with $\as_{\ell_1} = b_{\ell_2}$; in fact, for any such pair
there is a mutual cancellation of the corresponding Gamma functions in \rref{meij}.
\par
Let us repeat a comment appearing also in the main
text: in comparison with the $H$-function, the $G$-function is
more frequently implemented in standard packages for symbolic or numerical computations
with special functions. We also mention that $G$ can be expressed as a linear combination
of generalized hypergeometric functions \cite{Luke} \cite{Nist}.
As a final comment, the notation employed here for $G$
is suggested by the {\tt{Mathematica}} command for this function, and the notation
for $H$ arises from a natural generalization of the style used for $G$
({\footnote{In the literature, the function $H$ of
Eq.\,\rref{hfox} is more frequently written as
$$ H^{m n}_{p q} \left( z \left| \barray{l} (a_1,A_1),...,(a_p,A_p) \\
(b_1, B_1),...,(b_q, B_q)  \farray \right. \!\!\!\right)~ $$
where $p := n + \ns$,
$((a_1,A_1),...,(a_p,A_p)) :=$
$((a_1,A_1),...,(a_n,A_n), (\as_1,\As_1)...,(\as_\ns, \As_\ns))$,
$q := m + \ms$ and
$((b_1,B_1),...,(b_1,B_q)) :=$
$((b_1, B_1),...,(b_m, B_m),(\bs_1,\Bs_1),...,(\bs_\ms, \Bs_\ms)$;
see, e.g., \cite{Wolf}.\par
The function $G$ of Eq.\,\rref{meij} is more frequently written as
$$ G^{m n}_{p q} \left( z \left| \barray{l} a_1,...,a_p \\
b_1,...,b_q  \farray \right. \!\!\!\right) $$
where $p := n + \ns$, $(a_1,...,a_p):=$ $(a_1,...,a_n,a^{*}_1,...,a^{*}_{\ns})$,
$q := m + \ms$, $(b_1,...,b_q) := (b_1,...,b_m,b^{*}_1,...,b^{*}_{\ms})$.
}}).
In the sequel we are interested in
the $H$- or $G$- function in some special cases, where all the parameters
and the variable are real.
\salto
\textbf{Some inverse Fourier transforms.}
\begin{prop}
\label{lemmamej}
\textbf{Lemma.} Let
\beq g(k) := { |k|^{2 a-d} \over 1 + |k|^{2 b} }~~\mbox{for $k \in \reali^d$}~~(b>a>0),
\qquad f := \FF^{-1} g~; \feq
then the following holds. \parn
(i) $g \in \LL^1$, so that $f \in \LL^\infty$. $f$ has the radial structure \rref{f1}
$f(x) = F(|x|)$
where
\beq F(\rho) := \int_{0}^{+\infty}
d \xi~{J_{d/2 - 1}(\rho\, \xi) \over (\rho\, \xi)^{d/2-1}} {\xi^{2 a - 1} \over 1 + \xi^{2 b}}
\quad \mbox{for $\rho \in (0,+\infty)$}, \label{efro} \feq
\beq F(0) := {\pi \over 2^{d/2} \Gamma(d/2) \,b\, \sin(\pi a/b)} \label{ef0}~. \feq
Moreover
\beq \| f \|_{\LL^\infty} = f(0) = F(0)~. \label{eff0} \feq
(ii) For all $\rho \in [0,+\infty)$,
\beq F(\rho) = {1 \over 2^{d/2} b}\,
H \left( \!\!\!\left. \barray{l} (1 - {a \over b},{1 \over b});  \\
(0, 1),(1 - {a \over b}, {1 \over b}); (1 - {d \over 2},1)  \farray \right|
\left({\rho \over 2}\right)^2 \right) \label{fh} \feq
(note that the right hand side contains an $H$-function as in
\rref{hfox} with $n=1$, $\ns=0$, $m=2$, $\ms=1$; in this
case Eq.\,\rref{c3b} gives $\alf = 2/b >0$). \parn
(iii) In the rational case
\beq b = {N \over M}~, \qquad N, M \in \{1,2,3,...\} \feq
$F$ can be expressed as follows:
\beq F(\rho) = {M \over 2^{d/2 + M -1} \pi^{M - 1} N^{d/2}}\,
G \left( \left. \barray{l} a_1,...,a_N;  \\
b_1,...,b_{N+M}; \bs_1,...,\bs_N  \farray \right| \left(\rho \over 2 N\right)^{2 N} \right),
\label{fg} \feq
where:
\beq a_\ell := 1 - {a \over N} - {\ell-1 \over M} \quad \mbox{for $\ell=1,...,M$}~; \feq
$$ b_{h} := {h-1 \over N} \quad\mbox{for $h=1,...,N$}~,\quad
b_{N+h} := - {a \over N} + {h \over M}\quad \mbox{for $h=1,...,M$}~; $$
$$ \bs_\ell := 1 - {d \over 2 N} - {\ell - 1 \over N}~~\mbox{for $\ell=1,...,N$} $$
(here we are considering a Meijer $G$-function as in Eq.\,\rref{meij}
with $n=N$, $\ns=0$, $m = N+M$, $\ms=N$; Eq.\,\rref{c3bg} gives $\alf = 2 N >0$).
\end{prop}
\textbf{Proof.}
(i) All statements in this item are readily proved using
Proposition \ref{furgen} with $G(\xi) = \dd{ \xi^{2 a-d} \over 1 + \xi^{2 b} }$;
in particular, this proposition gives
$$ F(0) = {1 \over 2^{d/2 - 1} \Gamma(d/2)}
\int_{0}^{+\infty} d \xi~{ \xi^{2 a -1} \over 1 + \xi^{2 b} } =
{\pi \over 2^{d/2} \Gamma(d/2)\, b \, \sin(\pi a/b)}~, $$
thus justifying \rref{ef0} (as for the integral in the above line, use Lemma \ref{lemap} with $u=1$). \parn
(ii) Here (and in the subsequent proof of (iii)) we assume $\rho>0$;
once Eq.\,\rref{fh} is proved for $\rho >0$, its extension
to $\rho=0$ follows from the continuity of $F$ at zero (and the same can
be said about Eq.\,\rref{fg} of item (iii)). \par
Let us start from the Mellin-Barnes representation of the Bessel functions, that
can be written as follows:
\beq J_{\nu}(x) = {1 \over 2 \pi i} \int_{\L} \!\! ds \left({x \over 2}\right)^{\nu + 2 s}
{\Gamma(-s) \over \Gamma(1 + \nu + s) } \label{mebes} \feq
under the following conditions:
\beq \mbox{$\L$ a path from $c - i \infty$ to $c + i \infty$ ($c \in \reali$)} \label{condel} \feq
$$ \mbox{such that the poles $s=0,1,2,...$ of $\Gamma(-s)$ are on the right of $\L$}~; $$
\beq x \in (0,+\infty),~ \nu \in \complessi~,~\mbox{Re}\, \nu + 2 c > 0 \label{condnuc} \feq
(see \cite{Par}, page 115, Eq.\,(3.4.21) with the related comments,
and \cite{Wat}, page 192, Eq.(7))
({\footnote{The following remark is
a variation of some considerations
from page 115 of \cite{Par}.
For $x \in (0,+\infty)$, $\nu
\in \complessi$, $c \in \reali$, one has
$\left| \left({x \over 2}\right)^{\nu + 2 s} \,
\Gamma(-s) / \Gamma(1 + \nu + s) \right|$
$\sim H_{\pm}/|\mbox{Im} s |^{\tiny{\mbox{Re}}\,\nu + 2 c + 1}$
for $s \vain c \pm i \infty$ along $\L$, for suitable constants
$H_{\pm}$ depending on $x,\nu,c$; therefore
the integral in \rref{mebes} is absolutely
convergent under the condition $\mbox{Re}\, \nu + 2 c > 0$,
that appears for this reason in Eq.\,\rref{condnuc}.}}).
In the sequel we use the representation \rref{mebes} with
\beq \nu = {d \over 2} -1 \geqs - {1 \over 2}~,~c > {1 \over 4} \feq
(so that $\nu + 2 c > 0$); expressing in this way the term $J_{d/2 - 1}(\rho \, \xi)$ in
Eq.\,\rref{efro}, and reversing therein the order of the integrations in $s$ and $\xi$ we conclude
\beq F(\rho) = {1\over 2 \pi i}\,{1 \over 2^{d/2-1}}  \int_{\L} ds
\left( {\rho \over 2} \right)^{2 s}
{\Gamma(-s) \over \Gamma(d/2 + s)} \int_{0}^{+\infty} \!\!\! d\xi \, {\xi^{2 a + 2 s - 1} \over 1 + \xi^{2 b}}~.
\label{fro} \feq
On the other hand, Eq.\,\rref{simi} with $u=1$
and $a$ replaced by $a + s$ gives
\beq \int_{0}^{+\infty} \!\!\! d\xi \, {\xi^{2 a + 2 s  - 1} \over 1 + \xi^{2 b}}
= {\Gamma(a / b + s / b) \, \Gamma(1 - a / b - s / b) \over 2 b}~,
\feq
so from \rref{fro} we get
\beq F(\rho) = {1\over 2 \pi i}\,{1 \over 2^{d/2} b}  \int_{\L} ds
\left( {\rho \over 2} \right)^{2 s}
{\Gamma({a / b} + {s / b}) \, \Gamma(-s) \,\Gamma(1 - {a / b} - {s / b}) \over \Gamma({d / 2} + s)}~.
\label{weca} \feq
Comparing this result with the general definition \rref{hfox} of the Fox $H$-function
we obtain the thesis \rref{fh}, provided that the path of integration separates
as required the left and right poles of the integrand. Indeed, in the case
of \rref{weca} the sets of left and right poles are
\parn
\vbox{
\beq \LL =
\{ \sigma \in \complessi~|~\sigma~\mbox{is a pole of the function $s \mapsto
\Gamma({a  \over b} + {s \over b})$} \}  \feq
$$ = \{ - a - b k~|~k =0,1,2,...\}~; $$}
\parn
\vbox{
\beq \RR =
\{ \sigma \in \complessi~|~\sigma~\mbox{is a pole of the function}~ s \mapsto \Gamma(-s)~\mbox{or} \feq
$$ \mbox{of the function}~ s \mapsto \Gamma(1 - {a \over b} - {s \over b}) \}
= \{ 0,1,2,...\} \cup \{b - a + b k~ |~k=0,1,2,...\}~. $$}
The requirement that the sets $\LL,\RR$ are on the left and on the right of $\L$
is a strengthening of the conditions \rref{condel} initially given on this path;
we can in any case choose $\L$ so as to fulfill these stronger requirements,
thus getting the thesis \rref{fh}. \parn
(iii) We maintain for $\L$ the choice employed in the proof of (ii), and
use Eq.\,\rref{weca} with $b=N/M$. After a change of variable
$s \vain N s$ in the integral therein, we get
\beq F(\rho) = {1\over 2 \pi i}\,{M \over 2^{d/2}}  \int_{\L/N} \!\!\!\!\!\! ds
\left( {\rho \over 2} \right)^{2 N s}
{\Gamma(M {a \over N} + M s) \, \Gamma(-N s) \, \Gamma(1 - M {a \over N} - M s) \over \Gamma({d \over 2} + N s)}~.
\label{wecan} \feq
In the above $\L/N$ indicates the path defined as follows:
if $s = \param(t)$ is a parametrization of $\L$, then
by definition $s = \param(t)/N$ is a parametrization of
$\L/N$. Of course, $\L/N$ goes from $c/N - i \infty$ to $c/N + i \infty$. \par
Now, in the integral representation of $F(\rho)$ each Gamma function
contains the term $s$ or $-s$ multiplied by a positive integer. On the other
hand, the known Gauss multiplication rule for the Gamma function
\cite{Nist}
states that, for $n \in \{1,2,3,...\}$,
\beq \Gamma(n z) = {n^{n z - 1/2} \over (2 \pi)^{n/2 - 1/2}} \prod_{h=1}^n \Gamma(z + {h-1 \over n})~; \feq
we use this rule with $n = M$ and $z = a/N + s$, or $n = N$ and $z = -s$, or $n = M$ and $z= 1/M - a/N -s$,
or $n = N$ and $z = d/(2 N) + s $, which allows to reformulate \rref{wecan} as
\parn
\vbox{
\beq F(\rho) = {1\over 2\pi i}\,{M \over 2^{d/2 + M -1} \pi^{M - 1} N^{d/2}} \feq
$$\times \int_{\L/N} \!\!\! ds
\left( {\rho \over 2 N} \right)^{2 N s} {
\prod_{h=1}^M \Gamma({a \over N} + {h-1 \over M} +s)
\prod_{h=1}^N \Gamma({h-1 \over N} - s)
\prod_{h=1}^M \Gamma(- {a \over N} + {h \over M} -s)
\over \prod_{h=1}^N \Gamma({d \over 2 N} + {h-1 \over N} + s )}~.
$$}
Now, comparing with the general definition \rref{meij} of the
$G$-function we obtain the thesis \rref{fg} (again, the path
separates correctly the left and right poles of the integrand).
\fine
\salto
\textbf{More on inverse Fourier transforms.} The forthcoming statement
relies on the modified Bessel function
of the second kind
(Macdonald function) $K_{\mu}$ and on the hypergeometric function ${}_2 F_1$
\cite{Nist} \cite{Wat}.
\begin{prop}
\label{lembes}
\textbf{Lemma.} Let $\mu, \si \in \reali$, and define
\beq g(k) := {K_{\mu}(|k|) \over |k|^{\si}} \qquad (k \in \reali^d\setminus\{0\})~. \feq
Then $g \in \LL^1$ if and only if
\beq |\mu| + \si < d~. \label{musi} \feq
Assuming this,
\beq f := \FF^{-1} g \feq
is in $\LL^\infty$ and possesses the radial structure \rref{f1} $f(x) = F(|x|)$, where
\beq F(\rho) :=
{\Gamma(d/2 + \mu/2 - \si/2) \, \Gamma(d/2 -\mu/2 - \si/2) \over 2^{\si + 1 - d/2} \,\Gamma(d/2)}\, \label{a40} \feq
$$ \times {\,}_2 F_1(d/2 +  \mu/2 - \si/2, d/2 -\mu/2 - \si/2; d/2 ; -\rho^2)~
\quad \mbox{for $\rho \in [0,+\infty)$}~. $$
In particular, if $\si = \mu$,
\beq F(\rho) = {2^{d/2-\mu-1} \Gamma(d/2 - \mu) \over (1 + \rho^2)^{d/2 - \mu}}
\qquad \mbox{~for $\rho \in [0,+\infty)$}~. \label{a41} \feq
Moreover
\beq f \in \LL^2 ~\Leftrightarrow~ g \in \LL^2 ~\Leftrightarrow~ 2 (|\mu| + \si) < d~. \label{punti}\feq
\end{prop}
\textbf{Proof.} We apply Proposition \ref{furgen} with $G(\xi) := K_{\mu}(\xi)/\xi^{\si}$.
Recalling that $K_{\mu}(\xi) = O(e^{-\xi}/\sqrt{\xi})$ for $\xi \vain + \infty$ and that
$K_{\mu}(\xi) \sim \mbox{const.} \xi^{-|\mu|}$ for $\mu \neq 0$ and $\xi \vain 0^{+}$,
$K_{0}(\xi) \sim - \ln \xi $ for $\xi \vain 0^{+}$, we see that
$\int_{0}^{+\infty} d\xi \,\xi^{d-1} |G(\xi)| < + \infty$, i.e., $g \in \LL^1$,
if and only if \rref{musi} holds. Assuming this
the general rule \rref{f1} for radial, inverse Fourier transforms gives in the present case
$\FF^{-1} g (x) = F(|x|)$ with
\beq F(\rho) := \int_{0}^{+\infty}
\!\!d\xi ~\xi^{d - \si - 1} {J_{d/2 - 1}(\rho \,\xi) \over (\rho \,\xi)^{d/2-1}} K_{\mu}(\xi)
\mbox{~~for $\rho>0$}, \quad F(0) := \lim_{\rho \vain 0^{+}} F(\rho)~. \label{ineq} \feq
The integral in \rref{ineq} can be computed using results
from \cite{Wat}, \S 13.45, page 410; Eq.(1) therein yields our Eq.\,\,\rref{a40} and Eq.(2) therein (in a variant taking into account
the identity $K_{-\mu} = K_{\mu}$) yields our Eq.\,\rref{a41}. \par
To conclude, let us justify the statements in \rref{punti}.
Since $f = \FF^{-1} g$, it is obvious that $f \in \LL^2 \Leftrightarrow g \in \LL^2$.
On the other hand,
\beq g \in \LL^2 ~\Leftrightarrow~ \int_{0}^{+\infty} \!\! d\xi \,\xi^{d-1}
\Big( {K_{\mu}(\xi) \over \xi^{\si}} \Big)^2 < +\infty ~\Leftrightarrow~ 2 (|\mu| + \si) < d~. \feq
The first equivalence is obvious; to obtain the second equivalence note that,
by the previously mentioned asymptotic behavior of $K_{\mu}$,
the function of $\xi$ in the above integral decays exponentially
for $\xi \vain + \infty$ and behaves like $1/\xi^{2(|\mu| + \si) - d + 1}$ for $\mu
\neq 0$ and $\xi \vain 0^{+}$, or like $(\ln \xi)^2/\xi^{2 \si - d + 1}$ for $\mu=0$ and
$\xi \vain 0^{+}$.
\fine

\end{document}